\documentclass{article}[12pt]              

\usepackage{graphicx} % Required for inserting images
\usepackage{srcltx,epsfig}
\usepackage{amsmath, amssymb, amsthm}
\usepackage{color}
\usepackage{float}
\usepackage{xcolor}
\usepackage{multirow}
\usepackage[margin=1in]{geometry}
\usepackage{subfig}
\usepackage{verbatim}
\usepackage{pgfplots}
\pgfplotsset{compat=1.18} 
\usepackage{mathrsfs}
\usepackage{array}
\usepackage{booktabs} 
\usepackage{placeins} 
\newcommand\pd[2]{\dfrac{\partial {#1}}{\partial {#2}}}
\newcommand\ppd[1]{\dfrac{\partial }{\partial {#1}}}

\newcommand{\hS}{\hat{S}}
\newcommand{\mH}{\mathcal{H}}
\newtheorem{remark}{Remark}[section]

\usepackage{bm}
\usepackage[colorlinks,linkcolor=blue,anchorcolor=blue,citecolor=red]{hyperref}
\usepackage[T1]{fontenc}
\usepackage{tikz}
\usetikzlibrary{shapes.geometric, arrows}
\usetikzlibrary{matrix}
\usepackage{algorithm}
\usepackage{algpseudocode}
\usepackage[numbers]{natbib}
\usepackage{makecell}
\usepackage{siunitx}
\newcommand{\bmx}{\bm{x}}

\newcommand{\dd}{\mathrm{d}}
\graphicspath{{images/}}

\providecommand{\keywords}[1]
{
  \small	
  \textbf{Keywords:} #1
}
\title{The unified gas kinetic wave-particle method for the neutron transport equation}
\author{Guangwei Liu\thanks{Beijing Computational Science Research Center,
		Beijing 100094, China, email: \href{mailto:guangweiliu@csrc.ac.cn}{guangweiliu@csrc.ac.cn}.}, 
		~~Shuang Tan\thanks{Institute of Applied Physics and Computational Mathematics, Beijing 100094, China, email: \href{mailto:tan_shuang@iapcm.ac.cn}{tan\_shuang@iapcm.ac.cn}.}, 
	~~Yanli Wang\thanks{Beijing Computational Science Research Center,
		Beijing 100094, China, email: \href{mailto:ylwang@csrc.ac.cn}{ylwang@csrc.ac.cn}.}}
\date{\today}
\numberwithin{equation}{section} 
\begin{document}

\maketitle
\tableofcontents

\newpage
%-------------------------------------------------
\newpage
\begin{abstract}
   The unified gas-kinetic wave-particle (UGKWP) method is proposed for the neutron transport equation, addressing the inherent multiscale nature of neutron propagation in both optically thin and thick regimes. UGKWP couples macroscopic diffusion and microscopic transport processes within a unified time-dependent framework, allowing a smooth transition between the free transport and diffusion regimes. This method is readily extended to multi-group neutron transport models and is applicable to both steady-state and eigenvalue problems. Several numerical examples, including the 1D and 3D single-group and 3D multi-group problems, are studied, indicating UGKWP a promising framework for scalable and accurate simulation of multigroup neutron transport in complex geometries.
 \end{abstract}
 

\section{Introduction}
The neutron transport equation (NTE), a derivative of the Boltzmann equation, is a cornerstone of nuclear engineering, governing the phase-space distribution of neutrons in applications ranging from reactor physics to radiation shielding. Numerically solving NTE is notoriously challenging due to its high dimensionality and the multi-scale nature of neutron interactions, which span from diffusion-dominated to free-streaming regimes.

These challenges have spurred the development of two principal classes of numerical methods as the deterministic and stochastic methods. For the deterministic methods, the discrete ordinate method ($S_N$) \cite{carlson1968transport, lewis1984computational} is widely used, but may suffer from ray effects and demands fine angular resolution for anisotropic scattering \cite{Lathrop1968RayEffects, lewis1984computational}. The spherical harmonics method $P_N$ \cite{marshak1947note, sanchez1982review}, where the orthogonal polynomial is utilized as the basis function to approximate the angular flux, is efficient for nearly isotropic problems but may produce nonphysical oscillations where angular gradients are strong \cite{marshak1947note}. For the method of characteristics (MOC) \cite{askew1972moc, liu2011new}, it is accurate but becomes complex and computationally demanding in unstructured 3D geometries \cite{Gentry2021MOC3D}. To address the multi-scale challenges inherent in NTE, a unified gas-kinetic scheme (UGKS) is proposed first for the rarefied gas kinetic problem \cite{xu2010unified}, and then to the radiative transfer \cite{sun2015asymptotic} and neutron transport equations \cite{shuang2019parallel, tan2020time}, etc. UGKS, as a finite-volume method, constructs numerical fluxes by evolving an integral solution of the kinetic model equation across a time step, capturing the fluid dynamics from the rarefied to the continuum limit. 

For the stochastic methods, the Monte Carlo method (MC) \cite{lux1991monte} is one of the most popular methods, but its convergence can be prohibitively slow, and it often suffers from high statistical noise in problems characterized by low-flux regions or thick shielding. Based on UGKS, a particle-based unified gas-kinetic particle (UGKP) method was subsequently developed, which couples the deterministic evolution of the macroscopic variables with a stochastic particle representation for the non-equilibrium components \cite{shi2020asymptotic, hu2024implicit}. Moreover, to further reduce the computational cost, the unified gas-kinetic wave-particle method (UGKWP) refines this approach by decomposing the full distribution function into a background "wave" component and a "particle" component \cite{zhu2019unified, liu2020unified}. The wave, representing the near-equilibrium part, is solved deterministically, while the particle, representing the non-equilibrium deviation, is tracked stochastically. This wave-particle decomposition has demonstrated remarkable success in the transport problems with similar multi-scale characteristics, such as rarefied gas flow \cite{zhu2019unified, liu2020unified}, radiative transfer \cite{li2020unified}, and plasma physics \cite{pu2025unified}.

Given its robust physical foundation and proven multi-scale capabilities, UGKWP is expected to efficiently and accurately solve the neutron transport equation from the transport-dominated to free-streaming regimes. In this paper, UGKWP is proposed for the neutron transport equation, where the complex interaction effect is systematically constructed, including absorption, scattering, fission, and the source term. Within its multi-scale framework, both the macroscopic moment equations and the microscopic kinetic equation are consistently coupled within each time step, where an integral solution is utilized to construct the numerical flux across cell interfaces, bridging the different physical scales. The key advantage of this design is that the spatial and temporal resolution of the scheme is no longer restricted by the neutron's mean free path and collision time, enabling efficient and unified simulations across all physical regimes, from free transport to diffusion. 

For the detailed application, UGKP is first presented based on UGKS \cite{shuang2019parallel, tan2020time} and then UGKPW is constructed based on UGKP by adding the concept of the wave-particle for both the single and multi-group models. For UGKP, it is a particle implementation of UGKS, where the evolution of the particles is divided into a free transport step and a collision step. For the free transport step, the particles are tracked accurately along the characteristic lines. For the collision step, the macroscopic variables are updated by the numerical flux contributed from the particles' motion, and the microscopic distribution function is updated by resampling the updated macroscopic variables. For UGKWP, the free transport particles in UGKP are split into collisionless and collisional ones, where the collisional particles, which can be expressed in terms of macroscopic variables, will no longer be resampled, further reducing the computational cost. Several numerical examples are studied to validate UGKWP for neutrons, including the 1D and 3D single-group examples and 3D multi-group problems.

The following paper is arranged as follows. The neutron transport equation and related properties are introduced in Sec. \ref{Sec:NTE}. In Sec. \ref{sec:ugkwp_single_group}, the construction of UGKWP is presented in detail for the single-group model, including the numerical discretization and the particle simulation. Then, UGKWP is extended to the multi-group model in Sec.\ref{sec:multi-group}, with the properties of UGKWP in the free-streaming and diffusion-dominated regimes as well. Several numerical examples are studied in Sec. \ref{sec:num} to validate UGKWP, and the conclusion is provided in Sec. \ref{sec:conclusion} with the boundary conditions and related parameters listed in App. \ref{sec:app}.

%-------------------------------------------------
\section{Neutron transport equation}
\label{Sec:NTE}
In this section, the model equations in the neutron transport simulation are first introduced. If the delayed neutron is not considered, the neutron transport equation \cite{davison1958neutron} has the form as 
\begin{equation}
    \label{eq:neu_eq}
        \begin{array}{l}
            \displaystyle{\frac{1}{v} \frac{\partial \phi(t, \boldsymbol{x}, \Omega, E)}{\partial t}}+\Omega \cdot \nabla \phi(t, \boldsymbol{x}, \Omega, E)+\Sigma(\boldsymbol{x}, E) \phi(t, \boldsymbol{x}, \Omega, E) \\[12pt]
            \qquad =\displaystyle{\frac{1}{4 \pi}} \int_{4\pi}\int_0^{+\infty} \Sigma_{s}\left(\boldsymbol{x}, \Omega' \cdot \Omega, E' \rightarrow E\right) \phi\left(t, \boldsymbol{x}, \Omega', E'\right) \mathrm{d} E' \mathrm{d} \Omega' \\[12pt]
            \qquad +\displaystyle{\frac{\chi(E)}{4 \pi}} \int_0^{+\infty} \nu \Sigma_{f}\left(\boldsymbol{x}, E'\right) \int_{4\pi} \phi\left(t, \boldsymbol{x}, \Omega', E'\right) \mathrm{d} \Omega' \mathrm{d} E'+q(t, \boldsymbol{x}, \Omega, E),
        \end{array}
\end{equation}
where $\phi$ is the neutron angular flux who depends on  $(\bm{x}, t, \Omega, E)$. Here, $\Omega$ is the angular variable, $\bm{x}$ is the spatial variable and $E$ is the neutron energy. $v$ is the neutron velocity which is defined by the neutron energy $E$ as
\begin{equation}
    \label{eq:velocity_energy}
    E = \frac{1}{2}mv^{2},
\end{equation}
where $m$ is the neutron mass. $\Sigma(\bm{x}, E)$ is the total cross section, and $\Sigma_{s}\left(\boldsymbol{x}, \Omega' \cdot \Omega, E' \rightarrow E\right)$ is the scattering cross section form $(\Omega', E')$ to $(\Omega, E)$. $\Sigma_{f}\left(\boldsymbol{x}, E'\right)$ is the fission cross section and $\nu$ is the mean number of fission neutrons produced in a fission. $\chi(E)$ is the fission spectrum. $q(t, \boldsymbol{x}, \Omega, E)$ is the external source term.
In the single-group neutron transport equation, the neutron velocity is commonly set to $v=1$ after non-dimensionalization. For multi-group neutron transport equations, the velocity $v$ must be specified according to the discrete energy group divisions.
In the numerical simulation, the neutron angular flux $\phi$ is discretized in the high-dimensional phase $(t, \bm{x}, \Omega, E)$. For the energy discretization, the multi-energy-group approximation is utilized \cite{lewis1984computational}, where the energy $E$ is divided into $G$ intervals, with the group angular flux in each $g (1\leqslant g \leqslant G)$-th group defined as 
\begin{equation}
    \phi_{g}(t, \boldsymbol{x},\Omega) = \int_{E_{g-1}}^{E_{g}}\phi(t, \boldsymbol{x},\Omega, E)\mathrm{d}E.
\end{equation}
Introducing the energy function $h(E)$, which satisfies 
\begin{equation}
    \label{eq:energy_fun}
    \phi(t, \boldsymbol{x}, \Omega, E)\approx h(E)\phi_{g}(t, \boldsymbol{x}, \Omega),\qquad E_{g-1}\leqslant E \leqslant E_{g}, \qquad \int_{E_{g-1}}^{E_{g}}h(E)\mathrm{d}E = 1,
\end{equation}
the multi-group transport equation can be derived as 
\begin{equation}
\label{eq:multi-group}
        \begin{aligned}
        &\frac{1}{v_{g}} \frac{\partial \phi_{g}(t, \boldsymbol{x}, \Omega)}{\partial t}  +\Omega \cdot \nabla \phi_{g}(t, \boldsymbol{x}, \Omega)+\Sigma_{g}(\boldsymbol{x}) \phi_{g}(t, \boldsymbol{x}, \Omega)  = q_{g}(t, \boldsymbol{x}, \Omega)\\
         &+\frac{1}{4 \pi} \sum_{g'=1}^{G} \int_{4\pi} \Sigma_{s, g' \rightarrow g}\left(\boldsymbol{x}, \Omega' \cdot \Omega\right) \phi_{g'}\left(t, \boldsymbol{x}, \Omega'\right) \mathrm{d} \Omega' 
         +\frac{\chi_{g}}{4 \pi} \sum_{g'=1}^{G} \nu \Sigma_{f, g'}(\boldsymbol{x}) \int_{4\pi} \phi_{g'}\left(t, \boldsymbol{x}, \Omega'\right) \mathrm{d} \Omega'
        \end{aligned}
\end{equation}
by integrating the transport equation \eqref{eq:neu_eq} with energy $E$ in each interval, where
\begin{equation}
\label{eq:group_var}
    \begin{gathered}
        \frac{1}{v_{g}}  =\int_{E_{g-1}}^{E_{g}} \frac{1}{v} h(E) \mathrm{d} E, \qquad 
        \Sigma_{g}(\boldsymbol{x})  =\int_{E_{g-1}}^{E_{g}} \Sigma(\boldsymbol{x}, E) h(E) \mathrm{d} E,\qquad  \chi_{g}= \int_{E_{g-1}}^{E_{g}} \chi(E) \mathrm{d} E,\\       
        \Sigma_{s, g' \rightarrow g}\left(\boldsymbol{x}, \Omega' \cdot \Omega\right)=  \int_{E_{g-1}}^{E_{g}} \int_{E_{g}-1}^{E_{g}} \Sigma_{s}\left(\boldsymbol{x}, \Omega' \cdot \Omega, E' \rightarrow E\right) h\left(E'\right) \mathrm{d} E' \mathrm{d} E,\\
         \nu \Sigma_{f, g}(\boldsymbol{x})=  \int_{E_{g-1}}^{E_{g}} \nu \Sigma_{f}(\boldsymbol{x}, E) h(E) \mathrm{d} E,  \qquad  
        q_{g}(t, \boldsymbol{x}, \Omega)=  \int_{E_{g-1}}^{E_{g}} q(t, \boldsymbol{x}, \Omega, E) \mathrm{d}E.
        \end{gathered}
\end{equation}
With some rearrangement, the relaxation form of the multi-energy group neutron transport equation can be  derived as 
\begin{equation}
    \label{eq:neu_eq_relax}
        \displaystyle{\frac{\partial \phi_g(t, \boldsymbol{x},\Omega)}{\partial t}} + \xi_{g}\cdot\nabla\phi_{g}(t, \boldsymbol{x},\Omega) = \displaystyle{\frac{S_{g}(t, \boldsymbol{x},\Omega) - \phi_{g}(t, \boldsymbol{x},\Omega)}{\tau_{g}(\boldsymbol{x})}},
\end{equation}
where $S_{g}$ is the generalized source term whose form is
\begin{equation}
    \label{eq:source_term}
    \begin{aligned}
            S_{g}(t, \boldsymbol{x},\Omega) &= \frac{1}{4 \pi \Sigma_{g}(\boldsymbol{x})}\left( \sum_{g'=1}^{G} \int_{4\pi} \Sigma_{s, g' \rightarrow g}\left(\boldsymbol{x}, \Omega' \cdot \Omega\right) \phi_{g'}\left(t, \boldsymbol{x}, \Omega'\right) \mathrm{d} \Omega' \right.\\
        &\qquad  \left. +\chi_{g}\sum_{g'=1}^{G} \nu \Sigma_{f, g'}(\boldsymbol{x}) \int_{4\pi} \phi_{g'}\left(t, \boldsymbol{x}, \Omega'\right) \mathrm{d} \Omega' \right) + \frac{q_{g}(t, \boldsymbol{x}, \Omega)}{\Sigma_{g}(\boldsymbol{x})},
    \end{aligned}
\end{equation}
with 
\begin{equation}
    \label{eq:neutron_vel}
    \xi_g = v_g  \Omega, \qquad \tau_g = \frac{1}{v_g \Sigma_{g}(\boldsymbol{x})}
\end{equation}
the neutron transport velocity for the $g$-th group, and the characteristic collision time, respectively. 

Specially, when $G = 1$, the multi-group model  \eqref{eq:neu_eq_relax} is reduced into the single group model whose form is 
\begin{equation}
    \label{eq:single_group}
     \displaystyle{\frac{\partial \phi(t, \boldsymbol{x},\Omega)}{\partial t}} + \xi\cdot\nabla\phi(t, \boldsymbol{x},\Omega) = \displaystyle{\frac{S(t, \boldsymbol{x},\Omega) - \phi(t, \boldsymbol{x},\Omega)}{\tau(\boldsymbol{x})}},
\end{equation}
and considering an isotropic scattering assumption for neutrons,  $S_g(t, \bm x, \Omega)$ in \eqref{eq:source_term} is reduced into 
 \begin{equation}
     \label{eq:single_group_S}
     S(t, \bm x, \Omega) = \frac{\sum_s(\bm x) + \chi \nu \sum_f(\bm x)}{4 \pi \sum(\bm x)} \int_{4 \pi} \phi(t, \bm x, \Omega) \mathrm{d} \Omega + \frac{q(t, \bm x, \Omega)}{\sum(\bm x)},\qquad \tau=\frac{1}{v \sum(\bm x)}, \qquad \xi = v \Omega. 
 \end{equation}

 The neutron transport equation is challenging to solve due to its high dimensionality and strong couplings in angle and energy. The angular flux depends on time, space, angle, and energy, leading to high computational cost. Scattering and fission terms introduce global nonlocality, while large variations in cross sections and neutron speeds cause stiffness and numerical instability.

Classical methods, including discrete ordinates method ($S_N$) \cite{carlson1968transport, lewis1984computational}, spherical harmonics method ($P_N$) \cite{lewis1984computational, marshak1947note}, method of characteristics (MOC) \cite{liu2011new}, variational nodal method (VNM) \cite{zhang2022variational} and Monte Carlo method \cite{lux1991monte}, have been widely utilized. These methods offer advantages in structured grids, isotropic problems, or complex geometries, but often struggle with multiscale behavior and anisotropic fluxes. In this work, we explore the unified gas-kinetic wave-particle (UGKWP) method \cite{zhu2019unified}, which couples transport and collision processes within a multiscale framework, and UGKWP for the single-group model \eqref{eq:single_group} is first introduced, then extended to the multi-group systems.

%-------------------------------------------------
\section{UGKWP for single-group neutron transport equation}
\label{sec:ugkwp_single_group}
In this section, UGKWP for the single-group neutron transport equation is introduced. UGKWP is a multi-scale method for kinetic theory that combines macroscopic and microscopic dynamics with great technique. It was initially developed for the rarefied gas dynamics \cite{liu2020unified}, and has been successfully applied to plasma \cite{liu2021unified, pu2025unified}, and the radiative transport equation \cite{li2020unified, liu2023implicit}. Without loss of generality, we begin from the spatially 1D single-group model, and introduce the macroscopic scalar flux $\psi(t, x)$ as 
\begin{equation}
    \label{eq:macro_flux}
    \psi(t,x) = \int_{-1}^1 \phi(t, x, \mu) \dd \mu.
\end{equation}
Then the governing equation of the microscopic neutron angular flux $\phi(t,x,\mu)$ for the 1D single-group model has the form below
\begin{gather}
    \label{eq:1D}
     \displaystyle{\frac{\partial \phi(t,x,\mu)}{\partial t}} + \xi \pd{\phi(t,x,\mu)}{x} = \displaystyle{\frac{S(t,x,\mu) - \phi(t,x,\mu)}{\tau(x)}}, \\
     \label{eq:S}
          S(t, x, \mu) = \hat{S}(t, x) + \hat{q}(t, x, \mu),
          \qquad \tau = \frac{1}{v \Sigma(x)}, \quad \xi = v \mu, \quad \mu \in [-1, 1], 
\end{gather}
with 
\begin{equation}
\label{eq:hat_Sq}
     \hat{S}(t, x) = \dfrac{\Sigma_s(x) + \chi \nu \Sigma_f(x)}{2 \Sigma(x)} \psi(t,x), \qquad  \hat{q}(t, x, \mu) = \frac{q(t,x,\mu)}{\Sigma(x)}.
\end{equation}
Integrating \eqref{eq:1D} over $\mu$, we obtain the governing equation of $\psi$ as 
\begin{gather}
    \label{eq:macro}
     \pd{\psi(t, x)}{t} + \ppd{x} \int_{-1}^1 \xi \phi(t, x, \mu) \dd \mu =  \frac{G(t,x) - \psi(t,x)}{\tau(x)}, \\
     \qquad G(t,x) =\hat{G}(t,x) + \hat{Q}(t,x), \qquad \hat{G}(t, x) =  \int_{-1}^1 \hat{S}(t,x) \dd \mu, \qquad \hat{Q}(t, x) = \int_{-1}^1 \hat{q}(t,x, \mu) \dd \mu. 
\end{gather}
Since UGKWP is an enhanced particle method of UGKP by employing the wave-particle decomposition of the equilibrium part, and UGKP \cite{zhu2019unified} is a direct particle extension of UGKS, UGKS and UGKP will be first presented before introducing UGKWP.

\subsection{UGKS for the neutron transport equation}
\label{sec:UGKS}
UGKS was initially developed for rarefied gas dynamics \cite{xu2014direct, xu2010unified}, and has been successfully applied to plasma \cite{liu2017unified, liu2021unified}, radiative transport equation \cite{sun2015asymptotic, sun2017implicit}, and the neutron transport problems \cite{shuang2019parallel, tan2020time} as well. 

UGKS is implemented in the framework of the finite volume method (FVM), and the discretization for both the macroscopic scalar flux $\psi(t,x)$ and the microscopic neutron angular flux $\phi(t,x,\mu)$ is considered. 
Without loss of generality, the source term $q(t, x, \mu)$ is set as zero. For the $i$-th discrete finite volume cell $[x_{i-1/2}, x_{i+1/2}]$, the numerical solution at $t = t^{n}$ is defined as 
\begin{equation}
    \label{eq:num}
     \phi_{i}^{n} = \frac{1}{\Delta x}\int_{x_{i-1/2}}^{x_{i+1/2}} \phi(t^n, x, \mu) \mathrm{d}x, \qquad \psi_{i}^{n} = \frac{1}{\Delta x}\int_{x_{i-1/2}}^{x_{i+1/2}} \psi(t^n, x)\mathrm{~d}x.
\end{equation}
Since the discretization in angular space is not necessary in UGKP and UGKWP, we will also not introduce angular discretization in UGKS. Then, the numerical solutions are updated as 
\begin{align}
   \label{eq:sg_micro_update}
    \phi_{i}^{n+1}& =\phi_{i}^{n}  -\frac{\Delta t}{\Delta x}\left(F^{n}_{i+1 / 2}-F^{n}_{i-1 / 2}\right) +\Delta t\frac{S_{i}^{n+1}-\phi_{i}^{n+1}}{\tau_{i}} ,  \\
     \label{eq:sg_macro_update}
    \psi_{i}^{n+1}& =\psi_{i}^{n} -\frac{\Delta t}{\Delta x}\left(H^{n}_{i+1 / 2}-H^{n}_{i-1/2}\right) +\Delta t\frac{G_i^{n+1}-\psi_{i}^{n+1}}{\tau_{i}},
\end{align}
where the microscopic numerical flux $F^{n}_{i+1/2}$ is calculated by integrating $\phi$ with time $t$ at the cell interfaces $x_{i+1/2}$ as 
\begin{equation}
    \label{eq:sg_micro_flux}
    F^{n}_{i+1/2}=\frac{1}{\Delta t} \int_{t^{n}}^{t^{n+1}} \xi\phi(t, x_{i+1/2}, \mu)\mathrm{~d} t,
\end{equation}
and the macroscopic numerical flux $H^{n}_{i+1/2}$ is calculated by integrating $F^{n}_{i+1/2}$ with velocity $\mu$ as 
\begin{equation}
    \label{eq:sg_macro_flux}
    H^{n}_{i+1/2} = \int_{-1}^{1} F^{n}_{i+1/2} \mathrm{d}\mu. 
\end{equation}
Moreover, the generalized source term $S_{i}^{n+1}$ are 
\begin{gather}
    \label{eq:sg_S}
    S_{i}^{n+1} = \hS_{i}^{n+1} + \hat{q}_i^{n+1}, 
        \qquad  
   \hat{S}_i^{n+1} = \dfrac{\Sigma_s(x_i) + \chi \nu \Sigma_f(x_i)}{2 \Sigma(x_i)} \psi_i^{n+1}, \qquad \hat{q}_i^{n+1} =    
   \frac{1}{\Delta x}\int_{x_{i-1/2}}^{x_{i+1/2}} \hat{q}(t^{n+1}, x, \mu) \mathrm{d}x\\
   G_i^{n+1} = \hat{G}_i^{n+1} + \hat{Q}_i^{n+1},  \qquad   \hat{G}_i^{n+1} = \int_{-1}^1 \hS_i^{n+1} \dd \mu,
    \qquad \hat{Q}_i^{n+1} = \int_{-1}^1 \hat{q}_i^{n+1} \dd \mu, 
 \end{gather}
\begin{remark}
    For the discretization in angular space, the discrete ordinate method, denoted as $S_N$ is utilized, with the Gauss-Legendre quadrature adopted \cite{abramowitz1948handbook}. 
\end{remark}

In UGKS, the key ingredient is the construction of the multiscale numerical flux using the integral solution of the microscopic governing equation \eqref{eq:1D}. The exact form of the integral solution \cite{shuang2019parallel} is  
\begin{equation}
     \label{eq:integral_sol}
    \phi(t, x, \mu) = \int_{t_{0}}^{t} e^{-\frac{t-s}{\tau(x)}}\frac{1}{\tau(x)}S(s, x - \xi(t-s), \mu)\mathrm{d}s + e^{-\frac{t-t_{0}}{\tau(x)}}\phi(t_{0}, x - \xi(t-t_{0}), \mu),
\end{equation}
where $S(t, x, \mu)$ is defined in \eqref{eq:S} with $x-\xi(t-s)$ the particle trajectory, and $s$ the time variable from $t_0$ to $t$. $\phi(t_0, x - \xi(t-t_{0}), \mu)$ is the initial angular flux at $t_0$. To achieve second-order accuracy, $\hat{S}$ is expanded at the spatial space $x$ and time $t_0$ as 
\begin{equation}
    \label{eq:sec_expansion}
    \hat{S}(s, x - \xi(t-s)) \approx \hat{S}(t_0, x) + \xi(s-t)\partial_{x}\hat{S}(t_0, x) +  (s-t_0)\partial_{t}\hat{S}(t_0, x),
\end{equation}
and $\phi$ is expanded at the spatial space $x$ as 
\begin{equation}
    \label{eq:sec_phi}
    \phi(t_0, x-\xi(t-t_0), \mu) \approx \phi(t_0, x, \mu) + \xi(t_0 - t)\partial_x \phi(t_0, x, \mu),
\end{equation}
Substituting \eqref{eq:sec_expansion} and \eqref{eq:sec_phi} into \eqref{eq:integral_sol}, the integral solution becomes 
\begin{equation}
    \label{eq:evolution_sol}
    \phi(t, x, \mu) = \left(1 - e^{-\frac{t-t_0}{\tau}}\right) \left( \bar{S}(t, x, \mu)|_{t_0}  +\bar{q}(t,x,\mu)|_{t_0}\right) + e^{-\frac{t-t_0}{\tau}} \bar{\phi}(t, x, \mu)|_{t_0},    
   \end{equation}
with 
\begin{align}
    \label{eq:energy_2nd}
    \bar{S}(t, x, \mu)|_{t_0} &= \hat{S}(t_0, x) + \left(c_0 - \tau\right)\partial_{t}\hat{S}(t_0, x) + \left(c_0 e^{-\frac{t-t_0}{\tau}} - \tau\right)\xi \partial_{x}\hat{S}(t_0,x) , \\
    \bar{\phi}(t, x, \mu)|_{t_0} &= \phi(t_0, x, \mu) + \xi(t_0 - t)\partial_x \phi(t_0, x, \mu), \\
    \bar{q}(t,x,\mu)|_{t_0} & = \frac{1}{1 - e^{-\frac{t - t_0}{\tau}}}\int_{t_0}^t e^{-\frac{t-s}{\tau(x)}} \frac{\hat{q}(s, x - \xi(t-s), \mu)}{\tau(x)} \dd s,\qquad c_0 = \frac{t - t_0}{1 - e^{-\frac{t - t_0}{\tau}}}, 
\end{align}
Moreover, the first-order approximation gives
\begin{equation}
    \label{eq:energy_1st}
    \bar{S}(t, x, \mu)|_{t_0} = \hat{S}(t_0,x), \qquad \bar{\phi}(t,x, \mu)|_{t_0} = \phi(t_0, x, \mu). 
\end{equation}
With the integral solution \eqref{eq:evolution_sol}, the microscopic numerical flux \eqref{eq:sg_micro_flux} can be split into the macroscopic and microscopic part as 
\begin{equation}
    \label{eq:micro_flux}
    F_{i+1/2}^{n} \triangleq F^{n, \rm{ma}}_{i+1/2} + F^{n, \rm{mi}}_{i+1/2}  + F_{i+1/2}^{n, q},
\end{equation}
where $F^{n, \rm{ma}}_{i+1/2}$, $F^{n, \rm{mi}}_{i+1/2}$ 
and $F_{i+1/2}^{n, q}$ are the numerical flux by $\bar{S}$, $\bar{\phi}$ and $\bar{q}$, respectively.
\begin{align}
\label{eq:micro_flux_S}
        F^{n, \rm{ma}}_{i+1/2} & = \frac{1}{\Delta t} \int_{t^{n}}^{t^{n+1}} (1-e^{-\frac{t-t_n}{\tau}})\xi \bar{S}(t, x_{i+1/2}, \mu)|_{t^n} \dd t \\ \nonumber
         & \hspace{2cm} = 
        c_{1}\xi \hat{S}(t^{n}, x_{i+1/2}) + c_{2}\xi\partial_{t}\hat{S}(t^{n},x_{i+1/2}) + c_{3}\xi^{2}\partial_{x}\hat{S}(t^{n},x_{i+1/2}), \\
        \label{eq:micro_flux_phi}
        F^{n, \rm{mi}}_{i+1/2} & = \frac{1}{\Delta t} \int_{t^{n}}^{t^{n+1}} e^{-\frac{t-t_n}{\tau}} \xi\bar{\phi}(t, x_{i+1/2}, \mu)|_{t^n} \dd t
        = c_{4}\xi\phi(t^{n},x_{i+1/2},\mu) + c_{5}\xi^{2}\partial_{x}\phi(t^{n}, x_{i+1/2},\mu), \\
        F^{n, q}_{i+1/2} & = \frac{1}{\Delta t} \int_{t^{n}}^{t^{n+1}} \xi \bar{q}(t,x,\mu)\dd t.
\end{align}
with 
\begin{align}
\label{eq:ugks_flux_coe}
        &c_{1} = 1 - \frac{\tau a_1}{\Delta t}, \qquad c_{2} = \frac{\Delta t}{2} - \tau + \frac{\tau^{2} a_1}{\Delta t}, \qquad c_{3} = \frac{2\tau^2 a_1}{\Delta t}-\tau a_0-\tau,\\
        & c_{4} = \frac{\tau a_1}{\Delta t}, \qquad c_{5} = \tau a_0 - \frac{\tau^2 a_1}{\Delta t}, \qquad a_0 = e^{-\frac{\Delta t}{\tau}}, \qquad a_1 = 1 - a_0.
\end{align}
 Correspondingly, the macroscopic flux \eqref{eq:sg_macro_flux} can also be split into the macroscopic and microscopic parts as 
\begin{equation}
    \label{eq:macro_flux_H}
    H_{i+1/2}^n \triangleq H_{i+1/2}^{n, \rm{ma}} + H_{i+1/2}^{n, \rm{mi}} + H_{i+1/2}^{n, q}, \qquad H_{i+1/2}^{n, s} = \int_{-1}^1 F^{n, s}_{i+1/2} \dd \mu, \qquad s = {\rm ma}, {\rm mi}, q. 
\end{equation}
\begin{remark}
The integral in $F^{n, q}_{i+1/2}$ and $H^{n, q}_{i+1/2}$ can be easily obtained, since the expression of $q(t,x,\mu)$ is already known. In the simulation, the variables together with the spatial and temporal derivatives are all reconstructed with $S^{n+1}_{i}, S^{n}_{i}, S^{n}_{i+1}$ and $\phi_i^n, \phi_{i+1}^n$, which we will introduce in the following sections.
\end{remark}

As is stated in \cite{zhu2019unified}, the scale-adaptive numerical flux adopted in UGKS can capture the multiscale flow evolution, and the time step length in UGKS is not constrained by the particle mean free path and the collision in the transport modeling as well. However, the stochastic particle methods possess very high efficiency in the kinetic theory, especially for the problem where the mean free path is large in the three-dimensional case. Therefore, the particle implementation of UGKS, which is the unified gas kinetic particle method (UGKP) is introduced in the next section.

\subsection{Unified gas-kinetic particle method (UGKP)}
UGKP is a particle implementation of the unified gas kinetic scheme (UGKS), and has also been proposed for the rare gas dynamics \cite{zhu2019unified}, and the radiative transfer equations \cite{hu2024implicit, shi2020asymptotic, shi2021improved}. Here, we will introduce UGKP for the single-group model of the neutron transport equation. 

\begin{figure}[!hptb]
	\centering
	\includegraphics[width = 0.6\textwidth]{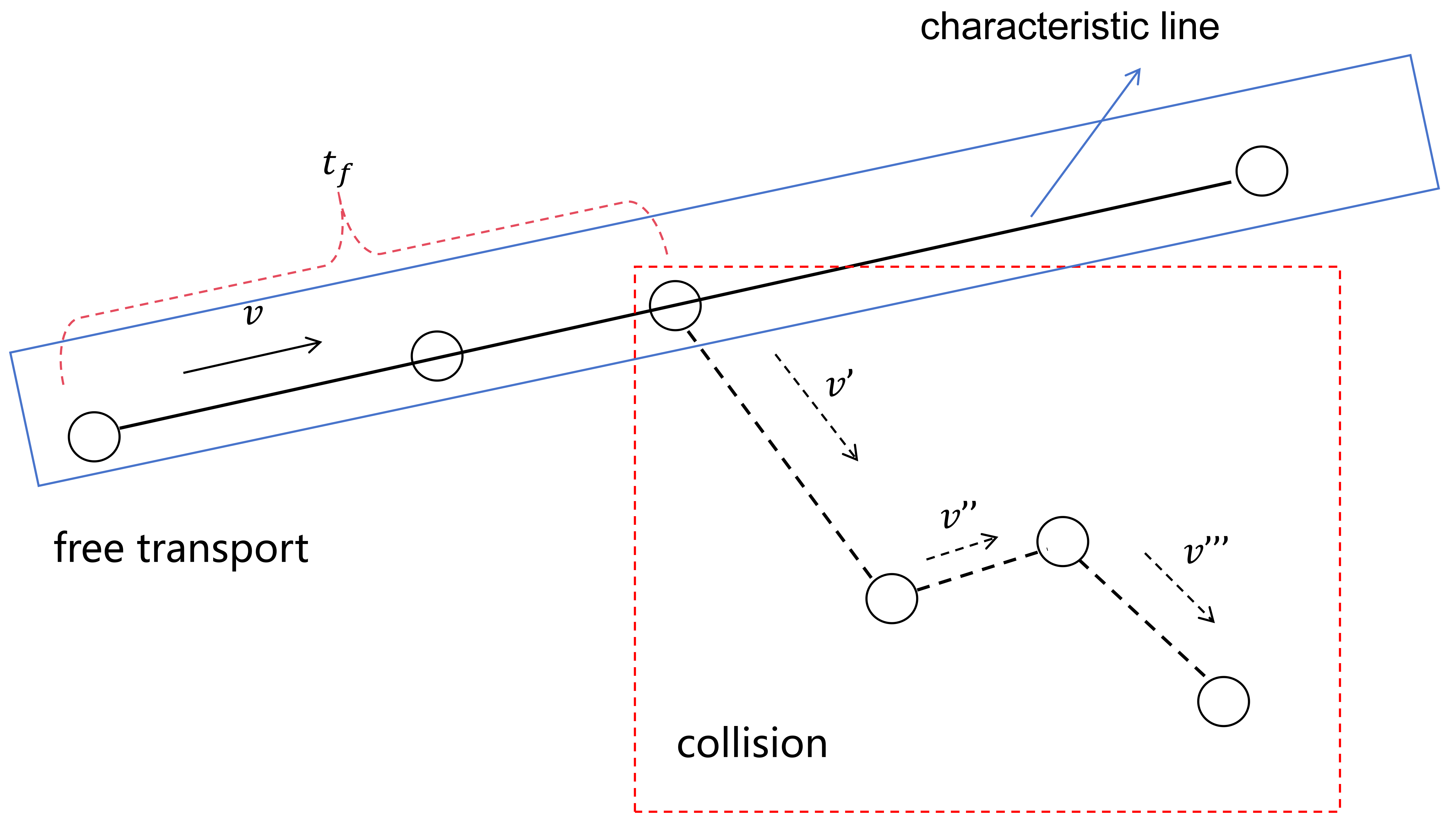}
	\caption{The evolution of one particle in a time step $\Delta t$, where the particle will keep free transport until $t_f$, and first collision occurs.}  
	\label{fig:evolution}
\end{figure}

The evolution of a particle in a time step $\Delta t$ is illustrated in Fig. \ref{fig:evolution}. For simplicity, all the effects in $\hat{S}$ are called collision, including scattering, fission, etc. In a time step, a particle will continue to undergo free transport until a collision occurs. Before the first collision, the particle will keep the characteristic line and retain its initial discrete velocity. Once the collision occurs, the particle's velocity changes, and its exact location and velocity cannot be traced any longer. UKGP has the same modeling scale as UGKS, where the time step length $\Delta t$ could be much larger than the kinetic scale particle collision time. Assuming the free transport time before the first collision as $t_f$, a huge number of collisions can happen during the time $(t_f, \Delta t)$. Therefore, tracking the motion of the particles in UGKP will be split into the free transport part and the collision part. 

We begin from the integral solution \eqref{eq:evolution_sol}, where the probability of maintaining the initial distribution function through the particle's free transport is given by $c_0$. Therefore, the motion of all the particles before their first collision can be tracked statistically. Then, the effect of collision is to annihilate the particles and recover them by resampling from the special distribution $S$. 

UGKP is implemented in the framework of a finite volume method, which models the above physical process on the scale of mesh size and time step. Precisely, UGKP is updated in two processes as the free transport process and the collision process. In the free transport process, the particles are tracked accurately along the characteristic lines, and the macroscopic scalar flux $\psi(t, x)$ is updated by the numerical flux contributed by these free transport particles. In the collision process, the macroscopic scalar flux $\psi(t, x)$ is updated by the numerical flux contributed from the particles' motion, while the microscopic neutron angular flux $\phi(t, x, \mu)$ is updated by resampling according to the updated $\psi(t,x, \mu)$. We will then introduce these processes in detail. 
\paragraph{Free transport process}
As shown in \eqref{eq:evolution_sol}, the cumulative distribution for particles free transport time at time $t^n$ is 
\begin{equation}
    \label{eq:time}
    \mathcal{G}(t) = e^{-\frac{t-t^n}{\tau}}.
\end{equation}
Thus, the free transport time of a particle with a time step $\Delta t = t^{n+1}-t^n$ is determined by inversion of the cumulative distribution function as in \cite{zhu2019unified}
\begin{equation}
    \label{eq:tf}
    t_f = \min(-\tau \ln(r_0), \Delta t),
\end{equation}
where $r_0$ is a random number generated from a uniform distribution $[0, 1]$. For each particle, we assume that each particle has the same mass $m_e$. At time $t= t^n$, a particle at position $x_p^n$ with the velocity $\xi_p$, can be accurately tracked as 
\begin{equation}
    \label{eq:free_trans}
    x_p = x_p^n + \xi_p t_f. 
\end{equation}
If $t_f = \Delta t$, then this particle is collisionless. In this case, $x_p^{n+1} = x_p$, the microscopic states of these particles are all updated. In $t_f < \Delta t$, the collision will occur at least once during $(t_f, \Delta t)$ for this particle, which we will discuss in the collision process. During the free transport process, the contribution to the numerical flux at the cell interface $x_{i+1/2}$ is obtained by counting the particles across the interface as 
\begin{equation}
    \label{eq:micro_flux_1}
    \mH_{i+1/2}^{n, {\rm mi}} =\frac{1}{\Delta t} \sum_{p \in X_{i+1/2}} {\rm sign}(\xi_p) m_e,
\end{equation}
where $X_{i+1/2}$ is the set of the particles moving across the interface $x_{i+1/2}$ during the free transport process. Compared to UGKS, the free transport process recovers the fluxes contributed by the initial distribution function $\bar{\phi}(t, x, \mu)|_{t_n}$, and the microscopic numerical flux $\mH_{i+1/2}^{n, {\rm mi}}$ is corresponding $H_{i+1/2}^{n, {\rm mi}}$ in \eqref{eq:macro_flux}. 

For now, the free transport process is carried out for both the microscopic and macroscopic levels. For the collision process, we will first introduce the implementation of the macroscopic level followed by the microscopic level. 

\paragraph{Collision process in the macroscopic level}
In the collision process, the particles collide and keep on moving and colliding. During this process, once the particles move across the mesh interfaces, they will contribute to the macroscopic flux as well. However, the motion of the particles in the collision process can not be exactly followed, and the related numerical flux can not be directly obtained as in the free transport process. Thus, the method in UGKS is utilized here to obtain the macroscopic numerical flux. The integral solution \eqref{eq:evolution_sol} shows that the macroscopic numerical flux $\mH_{i+1/2}^{n, {\rm ma}}$ is the same as $H_{i+1/2}^{n, {\rm ma}}$ in \eqref{eq:macro_flux_H}.

To obtain $H_{i+1/2}^{n, \rm{ma}}$,  $\hS(t^n, x_{i+1/2})$, $\partial_t \hS(t^n, x_{i+1/2})$, and $\partial_x \hS(t^n, x_{i+1/2})$ in the microscopic numerical flux $F_{i+1/2}^{n, \rm{ma}}$ \eqref{eq:micro_flux_S} should be is reconstructed firstly, and the reconstruction method in UGKS \cite{shuang2019parallel} is utilized here. We begin with the reconstruction of  $\phi_{i+1/2}^n$, where the linear reconstruction is utilized to achieve a second-order scheme. Thus, near the interface $x_{i+1/2}$, it is reconstructed as 
\begin{equation}
    \label{eq:recon_phi}
    \phi_{i+1/2}^n(t^n, x, \mu) = 
    \left\{\begin{aligned}
        \phi_{i+1/2}^{n,L} + \delta_x \phi^n_{i}(x - x_{i+1/2}), \qquad x < x_{i+1/2}, \\
        \phi_{i+1/2}^{n,R} + \delta_x \phi^n_{i+1}(x - x_{i+1/2}), \qquad x >  x_{i+1/2},
    \end{aligned}
    \right. 
\end{equation}
where $(\delta_x \phi^n)|_{i}$ is the reconstructed slope in the cell $i$ in the $x$- direction. Here, the van-Leer limiter \cite{van1977towards} is adopted in the reconstruction, and the reconstructed values at the interfaces are 
\begin{equation}
    \label{eq:rec_phi_1}
    \phi_{i+1/2}^{n,L} = \phi_{i}^{n} + \frac{\Delta x}{2}\delta_x \phi^n_{i}, \qquad 
    \phi_{i+1/2}^{n,R} = \phi_{i+1}^{n} - \frac{\Delta x}{2}\delta_x \phi^n_{i+1}.
\end{equation}
Then, $\hS_{i+1/2}^n$ is reconstructed as 
\begin{equation}
    \label{eq:rec_hs}
    \hS_{i+1/2}^n =\left\{
    \begin{array}{ll}
       \hS_{i+1/2}^{n,L} = \bar{\Sigma}_{i+1/2} \psi_{i+1/2}^{n,L}, & \text{in cell}~i, \\
       \hS_{i+1/2}^{n,R} = \bar{\Sigma}_{i+1/2} \psi_{i+1/2}^{n,R},& \text{in cell}~i+1,
    \end{array}
    \right.
    \qquad 
     \bar{\Sigma}_{i+1/2} = \frac{\Sigma_s(x_{i+1/2}) + \chi \nu \Sigma_f(x_{i+1/2})}{2\Sigma(x_{i+1/2})},
\end{equation}
where $\psi_{i+1/2}^{n,L(R)}$ is obtained from  $\phi_{i+1/2}^{n,L(R)}$, and the related derivatives are reconstructed as 
\begin{equation}
    \label{eq:rec_S}
    \begin{aligned}
     \xi  \partial_x \hS_{i+1/2}^n &=  \xi \frac{\hS_{i+1/2}^{n}- \hS_{i}^{n}}{\Delta x/2}H[\xi] + 
    \xi \frac{\hS_{i+1}^{n}- \hS_{i+1/2}^{n}}{\Delta x/2}(1 - H[\xi]) \qquad 
  \partial_{t} \hS_{i+1/2}^n& =\frac{\hS_{i+1/2}^{n+1}-\hS_{i+1/2}^{n}}{\Delta t},
    \end{aligned}
\end{equation}
where $H[x]$ is the Heaviside function. Substituting \eqref{eq:single_group_S}, \eqref{eq:micro_flux_S} and \eqref{eq:rec_S} into \eqref{eq:macro_flux_H}, $\mH_{i+1/2}^{n, \rm{ma}}$ is calculated as 
\begin{equation}
    \label{eq:macro_flux_1}
    \mH_{i+1/2}^{n, \rm{ma}} = \int_{-1}^{1} c_3 \xi^2 \partial_x \hS^n_{i+1/2} \dd \mu =
    \frac{2v^{2}c_3}{3\Delta x}\left(\hS_{i+1/2}^{n, L} - \hS_i^n + \hS_{i+1}^n - \hS_{i+1/2}^{n,R}\right). 
\end{equation}
Together with \eqref{eq:micro_flux}, \eqref{eq:macro_flux_H},\eqref{eq:micro_flux_1}, and \eqref{eq:macro_flux_1}, the macroscopic scalar flux $\psi(t, x)$ is updated as 
\begin{equation}
    \label{eq:update_UGKP_macro}
    \begin{aligned}
    \psi_{i}^{n+1}&=\alpha_{i}\left(\psi_{i}^{n}+v \Delta t \hat{Q}_{i}^{n+1}-\frac{\Delta t}{\Delta x}\left(\mH_{i+1/2}^{n,  {\rm ma}}-\mH_{i-1/2}^{n, {\rm ma}}\right) - \frac{\Delta t}{\Delta x}\left(\mH_{i+1/2}^{n,  {\rm mi}}-\mH_{i-1/2}^{n, {\rm mi}}\right) - \frac{\Delta t}{\Delta x}\left(\mH_{i+1/2}^{n,  q}-\mH_{i-1/2}^{n, q}\right)
        \right),
    \end{aligned}
\end{equation}
where 
$\mH_{i+1/2}^{n,  q} = H_{i+1/2}^{n,  q}$ and $\hat{Q}_{i}^{n+1}$ are the numerical flux contributed by the source term $q$, which are defined in \eqref{eq:macro_flux_H} and \eqref{eq:sg_S}, respectively. The coefficient $\alpha_i$ is 
\begin{equation}
    \label{eq:alpha}
    \alpha_i =  \left(1 + v \Delta t(\Sigma - \Sigma_s - \chi \nu \Sigma_f)|_{x = x_i}\right)^{-1}. 
\end{equation}
For now, the collision process on the macroscopic level is complete, and the next step is the advancement on the microscopic level during the collision process. 

\paragraph{Collision process at the microscopic level} 
As is discussed in the free transport process, the motion of all the particles at $(0, t_f)$ has been tracked. For the collisionless particles that $t_f = \Delta t$, their movement is already finished. For the collisional particles with $t_f < \Delta t$, at least one collision will occur during $(t_f, \Delta t)$. Here, we call all the effects in $\hS$ as collision. The collective impact of the collisions is to force all the collisional particles in the local region to follow a specific distribution which is given by $\bar{S} + \bar{q}$. In this process, we will first delete the collisional particles and then resample them from their respective distribution functions. 

Following the thought in \cite{zhu2019unified}, we can obtain the macroscopic scalar flux $\psi$ at cell $i$ corresponding to the collisional particles at the end of the time step according to the conservation law as 
\begin{equation}
    \label{eq:psi_Col}
        \psi_i^{n+1, {\rm ma}} = \psi^{n+1}_{i} - \psi_i^{n+1, {\rm mi}},
\end{equation}
where $\psi_i^{n+1, {\rm mi}}$ are the conservative flow variables from the collisionless particles which survived after the free transport process as $t_f = \Delta t$. In this process, the collisional particles are deleted first and re-sampled. Precisely, the number of the re-sampled particles is decided by $\psi_i^{n+1, {\rm ma}}$ in cell $i$, and its distribution is decided by 
\begin{equation}
    \label{eq:phi_Col}
\phi_i^{n+1, {\rm ma}} = \psi_i^{n+1,{\rm ma}}\frac{\hat{\phi}_i^{n+1, {\rm ma}}  }{\int_{-1}^1\hat{\phi}_i^{n+1, {\rm ma}}\dd \mu}  , \qquad 
\hat{\phi}_i^{n+1, {\rm ma}} = \frac{1}{\Delta x}  \int_{x_{i-1/2}}^{x_{i+1/2}} \bar{S}^{n+1}(t^{n+1}, x, \mu)|_{t^n}  + \bar{q}^{n+1}(t^{n+1}, x, \mu)|_{t^n} \dd x. 
\end{equation}
Substituting \eqref{eq:update_UGKP_macro}, and the reconstruction \eqref{eq:rec_S} into \eqref{eq:phi_Col}, we will derive the resampled distribution of $\phi_i^{n+1, {\rm ma}}$. Moreover, if the first-order scheme is utilized, then $\bar{S}$ is approximated by \eqref{eq:energy_1st}, which is independent of $\mu$. Even more, if the source term $q$ is isotropic, it is also independent of $\mu$. In this case,  \eqref{eq:phi_Col} is independent of $\mu$, and will reduce into a uniform distribution on $[-1, 1]$. 
\begin{remark}
    From the integral solution \eqref{eq:evolution_sol}, we can deduce that the re-sampled distribution function \eqref{eq:phi_Col} is the same as 
    \begin{equation}
  \tilde{\phi}_i^{n+1, {\rm ma}} \triangleq  \left(1 - e^{-\frac{\Delta t}{\tau}}\right)  \frac{1}{\Delta x}  \int_{x_{i-1/2}}^{x_{i+1/2}}  \bar{S}(t^{n+1}, x, \mu)|_{t_n}  +\bar{q}(t^{n+1},x,\mu)|_{t_n} \dd x. 
    \end{equation}
\end{remark}

Together with the particles that survived in the free transport process, those re-sampled according to the source term \eqref{eq:phi_Col}, we obtain all the particles at the end of this time step, and the microscopic particles have been updated. 

\begin{figure}[!hptb]
	\centering
	\includegraphics[width = 0.8\textwidth]{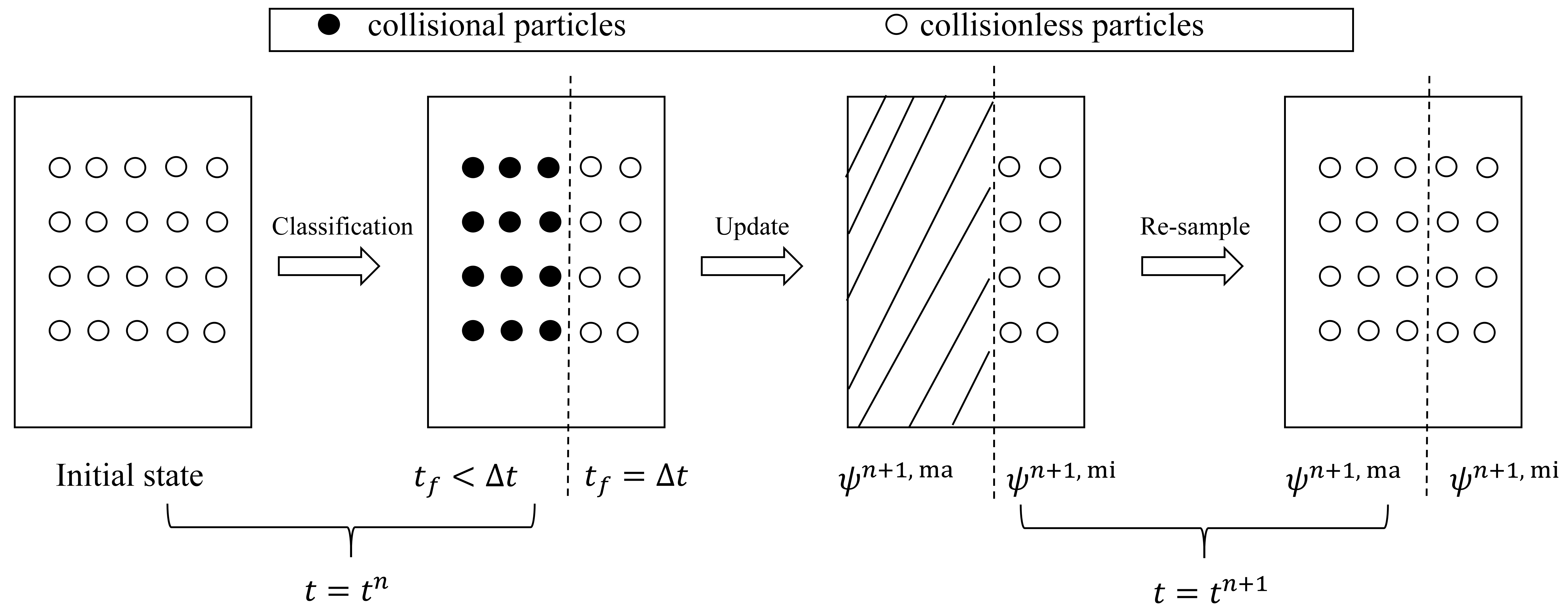}
	\caption{Diagram to illustrate the composition of the particles during time evolution in UGKP.}  
	\label{fig:ugkp}
\end{figure}

\subsection{Unified gas-kinetic wave-particle method (UGKWP)}
In this section, UGKWP is introduced, which is an enhanced particle method based on UGKP, where the concept of the wave-particle is added in UGKWP compared to UGKP. 

 As shown in Fig. \ref{fig:ugkp}, the distribution function is recovered by the re-sampled particles from \eqref{eq:phi_Col}, and the collisionless particles survived after the free transport process. In the next time step $[t^{n+1}, t^{n+2}]$, the re-sampled particles $\psi_i^{n+1, {\rm ma}}$ are reclassified into collisionless particles  $\psi_i^{n+1, {\rm ma}, {\rm free}}$ and collisional particles $\psi_i^{n+1, {\rm ma}, {\rm col}}$ according to $t_f$ again as 
 \begin{equation}
     \label{eq:phi_col_1}
     \psi_i^{n+1, {\rm ma}} = \psi_i^{n+1, {\rm ma}, {\rm free}} + \psi_i^{n+1, {\rm ma}, {\rm col}}.
 \end{equation}
  In the free transport process, only the collisionless particles $\psi_i^{n+1, {\rm ma}, {\rm free}}$ will remain as particles at the end of the time step $[t^{n+1}, t^{n+2}]$, while the collisional particles $\psi_i^{n+1, {\rm ma}, {\rm col}}$ will be deleted first, and re-sampled again. Therefore, only the collisionless particles $\psi_i^{n+1, {\rm ma}, {\rm free}}$ should be re-sampled at $t^{n+1}$, while the contribution of $\psi_i^{n+1, {\rm ma}, {\rm col}}$ can be computed analytically.

 From the cumulative distribution \eqref{eq:time}, we can deduce that the expectation of the proportion of the collisionless particles in each cell is $e^{-\Delta t / \tau}$, and the particles required to re-sample is 
 \begin{equation}
     \label{eq:phi_mi_free}
      \psi_i^{n+1, {\rm ma}, {\rm free}} = \psi_i^{n+1, {\rm ma}} e^{-\Delta t /\tau}.
 \end{equation}
 In this case, the free transport particles are combined with two parts, the first is the collisionless particles that survived in the last time step $\psi_i^{n+1, {\rm mi}}$ and the other is the free transport particles that re-sampled from $\psi_i^{n+1, {\rm ma}, {\rm free}}$. 
Similar to \eqref{eq:micro_flux_1}, the numerical flux contributed by the free transport particles at the cell interface $x_{i+1/2}$ is 
\begin{equation}
    \label{eq:micro_flux_wp_1}
    \mH_{i+1/2}^{n+1, {\rm mi}, {\rm free}} = \frac{1}{\Delta t} \sum_{p \in X_{i+1/2}}{\rm sign}(\xi_p) m_e,
\end{equation}
where $X_{i+1/2}$ is the same as that in \eqref{eq:micro_flux_1}. For the numerical flux contributed from the collisional particles $\psi_i^{n+1, {\rm ma}, {\rm col}}$, it can be computed analytically. 

We first obtain the numerical flux contributed from $\psi_i^{n+1, {\rm ma}}$ in the macroscopic level, where the numerical method in UGKS \eqref{eq:micro_flux_phi} is utilized as 
\begin{equation}
    \label{eq:micro_flux_wp_2} 
    \mH_{i+1/2}^{n+1, {\rm ma, total}} = \int_{-1}^{1} c_4 \xi \phi^{n+1, \rm ma}_{i+1/2} + c_5 \xi^2 \partial_x \phi^{n+1, \rm ma}_{i+1/2} \dd \mu,
\end{equation}
where $c_4$ and $c_5$ are the coefficients in \eqref{eq:ugks_flux_coe}. $\phi^{n+1, \rm ma}_{i+1/2}$ and $\partial_x \phi^{n+1, \rm ma}_{i+1/2}$ are reconstructed with linear reconstruction as \eqref{eq:recon_phi} and \eqref{eq:rec_phi_1}, where $\phi_i^{n+1, \rm ma}$ is the distribution function of the particles $\psi_i^{n+1, {\rm ma}}$ at $t^{n+1}$, whose form is \eqref{eq:phi_Col}. For the numerical flux \eqref{eq:micro_flux_wp_2}, it can be treated as the part that is contributed by the initial distribution in the total numerical flux \eqref{eq:sg_macro_flux}, where the initial distribution function of the particles is $\psi_i^{n+1, {\rm ma}}$ in \eqref{eq:phi_Col}. Then, for the numerical flux contributed from $\psi_i^{n+1, {\rm ma, free}}$, it is obtained by the deterministic discrete velocity method as 
\begin{equation}
    \label{eq:micro_flux_wp_3}
     \mH_{i+1/2}^{n+1, {\rm ma, free}} =\frac{1}{\Delta t}\int_{-1}^1 \int_{t^{n+1}}^{t^{n+2}} \xi \left(\phi_{i+1/2}^{n+1, {\rm ma, free}} + (t^{n+1} - t) \xi \partial_x \phi_{i+1/2}^{n+1, {\rm ma, free}} \right)\dd t \dd \mu,
\end{equation}
where $\phi_{i+1/2}^{n+1, {\rm ma, free}}$ and  $\partial_x \phi_{i+1/2}^{n+1, {\rm ma, free}}$ is reconstructed \eqref{eq:recon_phi} and \eqref{eq:rec_phi_1}, where $\phi_i^{n+1, \rm ma, free}$ is the distribution function of the particles $\psi_i^{n+1, {\rm ma, free}}$ at $t^{n+1}$ as
\begin{equation}
    \label{eq:ugkwp_phi_ma_free}
    \phi_{i+1/2}^{n+1, {\rm ma, free}} = e^{-\frac{\Delta t}{\tau}} \phi_i^{n+1, {\rm ma}}.
\end{equation}
Together with \eqref{eq:micro_flux_wp_2} and \eqref{eq:micro_flux_wp_3}, the numerical flux contributed by the collisional particles $\mH_{i+1/2}^{n+1, {\rm ma, col}}$ is 
\begin{equation}
    \label{eq:micro_flux_wp_4}
    \begin{aligned}
         \mH_{i+1/2}^{n+1, {\rm ma, col}} & =  \mH_{i+1/2}^{n+1, {\rm ma, total}} -  \mH_{i+1/2}^{n+1, {\rm ma, free}} \\
         & = \int_{-1}^1\xi\left(c_4 - e^{-\frac{\Delta t}{\tau}}\right)  \phi_{i+1/2}^{n+1, {\rm ma}} +\xi^2 \left(c_5 + \frac{\Delta t}{2} e^{-\frac{\Delta t}{\tau}}\right)\partial_x \phi_{i+1/2}^{n+1, {\rm ma}} \dd \mu. 
    \end{aligned}
\end{equation}
Therefore, in UGKWP  the advance for $\psi$ is changed from \eqref{eq:update_UGKP_macro} to 
\begin{equation}
    \label{eq:update_UGKWP_macro}
    \begin{aligned}
    \psi_{i}^{n+1}&=\alpha_{i}\left(\psi_{i}^{n}+v \Delta t \hat{Q}_{i}^{n+1}-\frac{\Delta t}{\Delta x}\left(\mH_{i+1/2}^{n,  {\rm ma}}-\mH_{i-1/2}^{n, {\rm ma}}\right)- \frac{\Delta t}{\Delta x}\left(\mH_{i+1/2}^{n,  q}-\mH_{i-1/2}^{n, q}\right) \right. \\
    &\qquad \left. - \frac{\Delta t}{\Delta x}\left(\mH_{i+1/2}^{n,  {\rm mi, free}}-\mH_{i-1/2}^{n, {\rm mi, free}}\right) - \frac{\Delta t}{\Delta x}\left(\mH_{i+1/2}^{n,  {\rm ma, col}}-\mH_{i-1/2}^{n, {\rm ma, col}}\right) 
        \right).
    \end{aligned}
\end{equation}
Compared to UGKP, the combination of the numerical flux $\mH_{i+1/2}^{n,  {\rm mi, free}}$ and $\mH_{i+1/2}^{n,  {\rm ma, col}}$ is the same at the numerical flux $\mH_{i+1/2}^{n,  {\rm mi}}$ in \eqref{eq:update_UGKP_macro}. For now, we have finished the evolution of UGKWP in a time step. The detailed algorithm, and the related properties of UGKWP are presented in the following section

\subsection{Total algorithm and related properties of UGKWP}
\label{sec:AP}
\begin{figure}[!hptb]
	\centering
	\includegraphics[width = 0.8\textwidth]{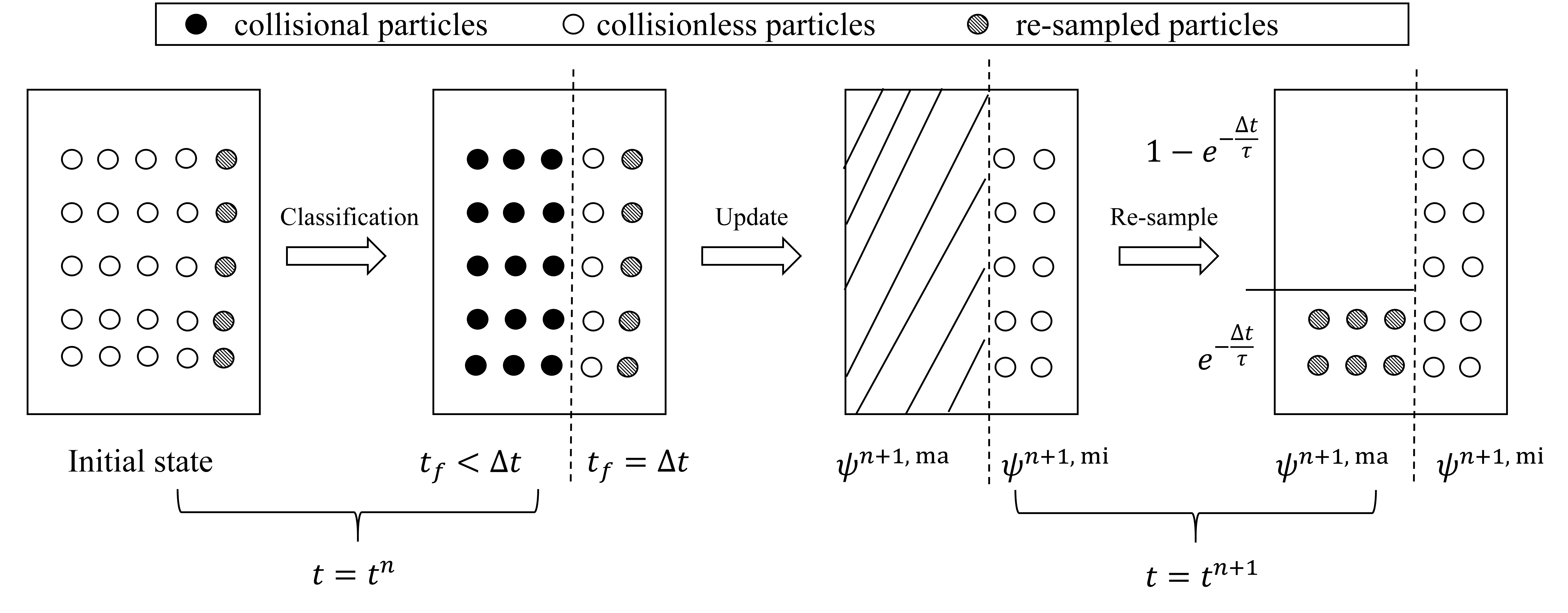}
	\caption{Diagram to illustrate the composition of the particles during time evolution in UGKWP.}  
	\label{fig:ugkwp}
\end{figure}

To make this numerical method clearer, the entire algorithm is listed in Alg. \ref{algo:ugkwp}, and a series of figures is plotted in Fig. \ref{fig:ugkwp} to illustrate this algorithm. Moreover, UGKWP is a multiscale method, which can preserve problems with various cross sections $\Sigma$, and it satisfies the following properties:
\begin{enumerate}
    \item In the free transport limit, the characteristic collision term $\tau(x)$, defined in \eqref{eq:neutron_vel}, goes to $\infty$, and the material is almost empty. In this case, $e^{-\frac{\Delta t}{\tau}}$ goes to $1$, and UGKWP is reduced into a purely particle method, where each particle is traced exactly by free transport.
    \item In the diffusion limit, the characteristic collision $\tau(x)$ goes to zero, and $e^{-\frac{\Delta t}{\tau}}$ goes to $0$. In this case, essentially all the particles are the collision particles, and no particle will be resampled from $\psi_i^{n+1, \text{ma, free}}$ defined in \eqref{eq:phi_mi_free}. UGKWP will become a macroscopic solver, which is very efficient.  
    \item In the transition regime, the macroscopic and microscopic processes are evolved on the wave-particle decomposition, where the update of the macroscopic scalar flux $\psi$ is the combination of macroscopic and microscopic numerical flux. Here, only the fully free transport particles will be resampled to further deduce the computational cost. 
\end{enumerate}

\begin{algorithm}[htbp]
    \caption{Numerical algorithm}
    \label{algo:ugkwp}
    \begin{algorithmic}[1]
        \item Obtain the initial state with the macroscopic scalar flux $\psi^{n}_i$, the collisionless particles $\psi_i^{n, {\rm mi}}$ evolved from the last time step, the collisionless particles $\psi_i^{n, {\rm ma, free}}$ sampled from the updated macroscopic scalar flux $\psi^{n, {\rm ma}}_i$. For the first step, $\psi_i^{0, {\rm mi}} = 0$, and $\psi^{n, {\rm ma}}_i =\psi^{0}_i$;
        \item Update the microscopic variables:
            \begin{enumerate}
                \item[2.1] Obtain the free transport time $t_f$ by the cumulative distribution \eqref{eq:tf}; 
                \item[2.2] Classify the particles $\psi_i^{n, {\rm mi}}$ into the collisionless particles and the collisional particles according to $t_f$;
                \item[2.3] Stream all the particles over the free transport time $t_f$. These particles including all the particles in $\psi_i^{n, {\rm mi}}$, and the free transport particles in $\psi_i^{n, {\rm ma, free}}$ whose free transport time is $t_f = \Delta t$; 
                \item[2.4] Calculate the free transport flux $\mH_{i+1/2}^{n, {\rm mi, free}}$;
                \item[2.5] Update the collisionless particles $\psi_i^{n+1, {\rm mi}}$, which contains two parts, one is survived particles from $\psi_i^{n, {\rm mi}}$ and other is the free transport articles from  $\psi_i^{n, {\rm ma, free}}$;
                \item[2.6] Delete the collisional particles from $\psi_i^{n, {\rm mi}}$;
                \item[2.7] Compute the numerical flux contributed from the unsampled particles $\psi_i^{n, {\rm ma, col}}$ from \eqref{eq:micro_flux_wp_4};
                  \end{enumerate}
    \item Update the macroscopic variables:
    \begin{enumerate}
        \item[3.1] Obtain the numerical flux contributed by the collisional particles \eqref{eq:macro_flux_1};
        \item[3.2] Update the macroscopic scalar flux $\psi_i^{n+1}$ by \eqref{eq:update_UGKWP_macro}, 
        the collisional particles $\phi_i^{n, {\rm ma}}$ by \eqref{eq:psi_Col}, and the collisionless particle in next time step $\psi_i^{n+1, {\rm ma, free}}$ by \eqref{eq:phi_mi_free};
        \item[3.3] Sample the collision particles $\psi_i^{n+1, {\rm ma, free}}$ according to the distribution function \eqref{eq:ugkwp_phi_ma_free};
    \end{enumerate}
\item Determine the computation for the next step: If the finish time is reached, stop the program. Otherwise, go to step 1.
    \end{algorithmic}
\end{algorithm}

\section{UGKWP for multi-group neutron transport simulation}
\label{sec:multi-group}
In Sec. \ref{sec:ugkwp_single_group}, UGKWP for the single-group model of the neutron transport equation is introduced. Here, we will extend it to the spatially multi-group problems. With the isotropic scattering assumption, the generalized source \eqref{eq:source_term} will be reduced into 
\begin{equation}
    \label{eq:source_term_iso}
    \begin{aligned}
            S_{g}(t, x, \mu) &=  \sum_{g'=1}^{G}\frac{\Sigma_{s, g' \rightarrow g}\left(x \right) + \chi_{g}\nu \Sigma_{f, g'}(x)}{2 \Sigma_{g}(x)} \int_{-1}^1  \phi_{g'}\left(t, x, \mu\right) \mathrm{d} \mu  + \frac{q_{g}(t, x, \mu)}{\Sigma_{g}(x)}.
    \end{aligned}
\end{equation}
The macroscopic flux $\psi_g$ for the $g$-group is defined as 
\begin{equation}
    \label{eq:psi_g}
    \psi_g(t, x) = \int_{-1}^1 \phi_g(t,x,\mu)\dd \mu,
\end{equation}
and the related governing equation for $\psi_{g}$ is as follows
\begin{equation}
    \label{eq:eq_psi_g}
    \begin{gathered}
    \pd{\psi_g(t, x)}{t} + \ppd{x} \int_{-1}^1 \xi \phi_g(t, x, \mu) \dd \mu =  \frac{G_g(t,x) - \psi_g(t,x)}{\tau_g(x)}, \\
     \qquad G_g(t,x) =\hat{G}_g(t,x) + \hat{Q}_g(t,x), \qquad \hat{G}_g(t, x) =  \int_{-1}^1 \hat{S}_g(t,x) \dd \mu, \qquad \hat{Q}_g(t, x) = \int_{-1}^1 \hat{q}_g(t,x, \mu) \dd \mu,
     \end{gathered}
\end{equation}
where $\hat{\cdot}$ is defined as same as in \eqref{eq:hat_Sq} and \eqref{eq:macro} in the single-group model. Moreover, with $\psi_g$, \eqref{eq:source_term_iso} can be rewritten as 
\begin{equation}
    \label{eq:re_s}
     S_{g}(t, x, \mu) =  \sum_{g'=1}^{G}\frac{\Sigma_{s, g' \rightarrow g}\left(x \right) + \chi_{g}\nu \Sigma_{f, g'}(x)}{2 \Sigma_{g}(x)} \psi_{g'} + \frac{q_{g}(t, x, \mu)}{\Sigma_{g}(x)}.
\end{equation}
UGKWP for the multi-group model is similar to the single-group model, where the related parameters are substituted with the $g$-th group parameters. The only difference exists in the handling of the source term $S_g(t, x, \mu)$. In particular, the reconstruction $\hat{S}_{g,i+1/2}^n$ in \eqref{eq:rec_hs} changes into 
\begin{equation}
    \label{eq:rec_hs_mul}
    \hS_{g,i+1/2}^n =\left\{
    \begin{array}{ll}
       \hS_{g,i+1/2}^{n,L} = \sum\limits_{g'=1}^G \bar{\Sigma}_{g,g',i+1/2} \psi_{g',i+1/2}^{n,L}, & \text{in cell}~i, \\
       \hS_{g,i+1/2}^{n,R} = \sum\limits_{g'=1}^G \bar{\Sigma}_{g,g' i+1/2} \psi_{g', i+1/2}^{n,R},& \text{in cell}~i+1,
    \end{array}
    \right.
  \end{equation}
with 
\begin{equation}
    \bar{\Sigma}_{g,g',i+1/2} = \frac{\Sigma_{s, g'\rightarrow g}(x_{i+1/2}) + \chi_g \nu \Sigma_{f,g'}(x_{i+1/2})}{2\Sigma_g(x_{i+1/2})}.
\end{equation}
Meanwhile, the advance for $\phi_g$ in \eqref{eq:update_UGKWP_macro} is changed into 
\begin{equation}
    \label{eq:update_mg_UGKWP_macro}
    \begin{aligned}
    \psi_{i,g}^{n+1}&=\left(\psi_{g,i}^{n}  - \frac{\Delta t}{\Delta x}\left(\mH_{g,i+1/2}^{n,  {\rm ma}}-\mH_{g,i-1/2}^{n, {\rm ma}}\right)- \frac{\Delta t}{\Delta x}\left(\mH_{g,i+1/2}^{n,  q}-\mH_{g,i-1/2}^{n, q}\right) \right. \\
    &\qquad  - \frac{\Delta t}{\Delta x}\left(\mH_{g,i+1/2}^{n,  {\rm mi, free}}-\mH_{g,i-1/2}^{n, {\rm mi, free}}\right) - \frac{\Delta t}{\Delta x}\left(\mH_{g,i+1/2}^{n,  {\rm ma, col}}-\mH_{g,i-1/2}^{n, {\rm ma, col}}\right) 
         \\
        & \qquad \left. + \frac{\Delta t}{\tau_{g, i}}\left(
        \sum_{g'=1}^G \bar{\Sigma}_{g,g',i} \phi_{g',i}^{n+1}
                \right) - \frac{\Delta t\psi_{i,g}^{n+1}}{\tau_{g,i}} + \frac{\Delta t\hat{Q}_{i,g}^{n+1}}{\tau_{g,i}} \right),       
    \end{aligned}
\end{equation}
with 
\begin{equation}
    \bar{\Sigma}_{g,g',i} =  \frac{\Sigma_{s, g'\rightarrow g}(x_{i}) + \chi_g \nu \Sigma_{f,g'}(x_{i})}{2\Sigma_g(x_{i})}. 
\end{equation}
Here, the method to obtain the related numerical flux $\mH_{g, i+1/2}^{\ast}$ is the same as that in the single-group model. When updating $\phi_{i,g}^{n+1}$, a linear equation system for $\phi_{i,g}^{n+1}, g = 1, \cdots G$ is solved, the coefficient matrix of which is diagonal dominant for most cases and easy to solve.

\section{Numerical experiments}
\label{sec:num}
In this section, several numerical experiments are tested to validate UGKWP for the neutron transport equations. We begin with the spatial 1D benchmark problems for the single-group model, including the Heaviside source problem, the isotropic scattering problem, and the unsteady problem. Then, two spatially 3D problems, the block problem and the tophat problem, are studied. Finally, two spatially 3D problems for the multi-group model, the KUCA core two-group problem and the four-group axially heterogeneous FBR problem, are conducted. For all the numerical examples, the time step length is set as 
\begin{equation}
    \label{eq:time_step}
    \Delta t = \mathrm{CFL} \, \min\!\left(\frac{\Delta}{v}, \, \frac{3\Delta^2 \Sigma}{2v}\right),
\end{equation}
where $\Delta = \min\{\Delta x, \Delta y, \Delta z\}$, $v$ denotes the neutron velocity, which is commonly set to 1 after non-dimensionalization, and $\mathrm{CFL} \in (0,1)$. For the multi-group neutron transport problem, the neutron velocity $v_g$ varies with energy. It is all set as $v_g = 1$ for all groups after non-dimensionalization here, since only the steady-state problems are considered for the 3D multi-group model, and this will not influence the steady-state solution. This setting follows the strategy proposed in \cite{shuang2019parallel}.  

\subsection{1D single-group problems} 
\label{sec:1D}
In this section, three spatially 1-dimensional examples are studied, including the Heaviside source problem, the steady isotropic scattering problem, and the unsteady problem with an analytical solution, validating the accuracy of UGKWP. 
\subsubsection{1D Heaviside source problem}
\label{sec:1D_heaviside}
The 1D Heaviside source problem for the single-group model is tested. This problem is also studied in \cite{bindra2012radiative, shuang2019parallel}. The computation parameters are set the same as \cite{bindra2012radiative, shuang2019parallel} 
\begin{gather}
    \Sigma (x) = \hat{\Sigma} = 1,\qquad \Sigma_{s}(x) = \Sigma_{f}(x) = 0,\qquad q(x) = \mathrm{H}(x - 0.5),\\
    \phi(0,x,\xi) = 0, \qquad x \in[0,1], \qquad \xi \in [-1, 1],
\end{gather}
where $\mathrm{H}(x)$ is the Heaviside function. The inlet boundary condition \cite{nance1997role,nance1998role} is utilized as 
\begin{equation}
\label{eq:ex_bou}
    \phi(x_{w},\xi) = \left\{
    \begin{array}{ll}
        0, & x_w = 0,~\xi > 0, \\
       \frac{1}{\hat{\Sigma}},  &  x_w = 1,~\xi < 0.
    \end{array}
\right.
\end{equation}
There exists an analytical solution of the neutron angular flux $\phi(x, \xi)$ when the system reaches the steady state \cite{bindra2012radiative}, which has the form below
\begin{equation}
\label{eq:ex1_sol}
    \phi(x,\xi) = \left\{
    \begin{array}{ll}
      \frac{e^{-\hat{\Sigma}(x-0.5)/\xi}}{\hat{\Sigma}},             &x\leqslant 0.5,~\xi <0,  \\[2mm]
      \frac{1}{\hat{\Sigma}},   & x>0.5,~\xi < 0, \\[2mm]
      0,  & x\leqslant 0.5,~\xi > 0,\\
      \frac{1 - e^{-\hat{\Sigma}(x-0.5)/\xi}}{\hat{\Sigma}}, & x>0.5,~\xi > 0.
    \end{array}
   \right.
\end{equation}
In the simulation, the final time is set as $t = 5$, and the spatial mesh is $N_x = 50$ with $\mathrm{CFL}=0.4$. The macroscopic scalar flux $\psi(x)$ is approximated as 
\begin{equation}
\label{eq:ex1_psi}
    \psi(x) \approx \sum_{i=1}^{N} \omega_{i}\phi(x,\xi_{i}).
\end{equation}
where $[\xi_i, \omega_i], i = 1, \cdots, N$ are the Gauss–Legendre quadrature points and corresponding quadrature weight \cite{abramowitz1948handbook}.

\begin{figure}[!hptb]
	\centering
	\subfloat[macroscopic scalar flux $\psi$]{
		\includegraphics[width = 0.4\textwidth]{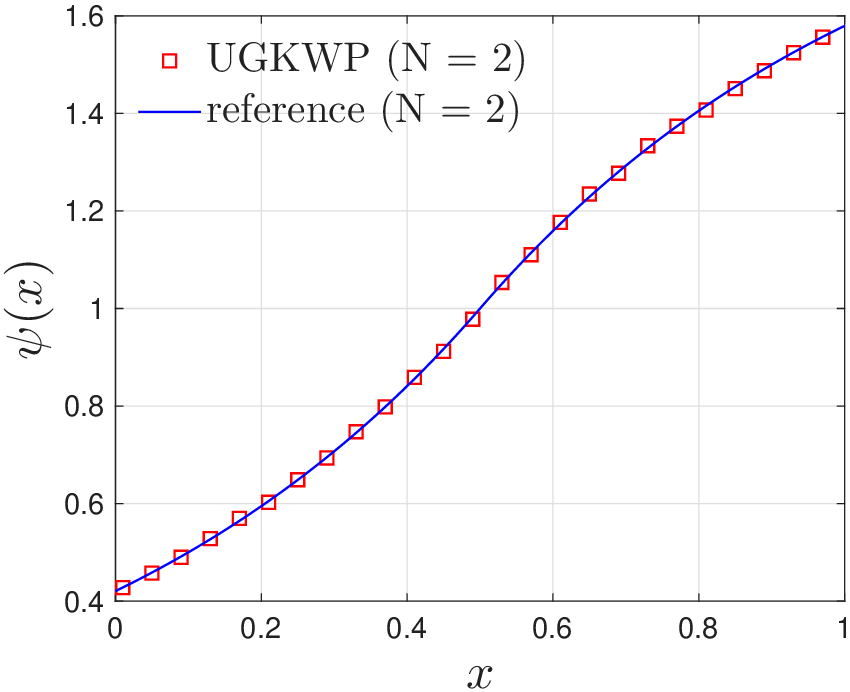}
        \label{fig:1D_heaviside_S2_1}
	}\qquad
        \subfloat[angular flux $\phi$]{
        \includegraphics[width = 0.4\textwidth]{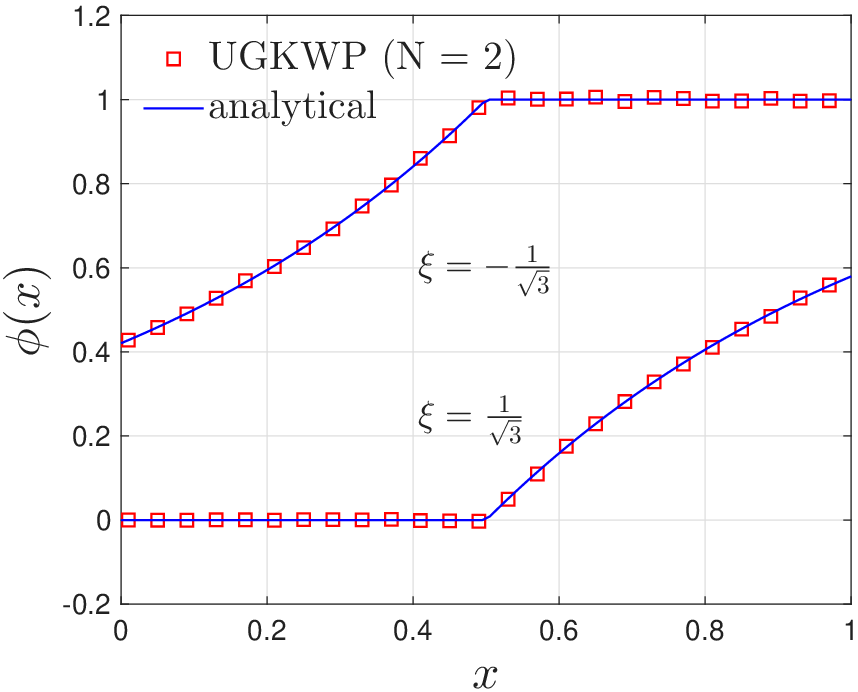}
        \label{fig:1D_heaviside_S2_2}
	}     
	\caption{(1D Heaviside source problem in Sec. \ref{sec:1D_heaviside}) The comparison of the numerical results by UGKWP and the reference solution for $N = 2$. (a) comparison of $\psi$, where the red square is the numerical solution obtained by UGKWP with the distribution of the microscopic velocity two-point distribution, and the blue line is the reference solution obtained by \eqref{eq:ex1_psi} with $N = 2$;   
    (b) the numerical solution of $\phi$ at $\xi = \pm \frac{1}{\sqrt{3}}$, where the red square is the numerical solution obtained by UGKWP with the distribution of the microscopic velocity two-point distribution, and the blue line is the analytical solution \eqref{eq:ex1_sol}.     
    }  
\label{fig:1D_heaviside_S2}
\end{figure}

First, to compare the numerical results with that in \cite{shuang2019parallel}, we first set $N = 2$ the same as that in \cite{shuang2019parallel}, and the detailed quadrature points and weights are 
\begin{equation}
    \label{eq:s2}
    \xi_1 = -\frac{1}{\sqrt{3}}, \qquad \xi_2 = \frac{1}{\sqrt{3}}, \qquad \omega_1 = \omega_2 = 1. 
\end{equation}
In this case, the distribution of the neutron transport velocity $\xi$ in UGKWP is treated as the two-point distribution detailed in \eqref{eq:s2}. Fig. \ref{fig:1D_heaviside_S2} presents the numerical results of the macroscopic scale flux $\psi(x)$, and the neutron angular flux $\phi(x,\xi)$, respectively. Precisely, Fig. \ref{fig:1D_heaviside_S2_1} shows that the macroscopic variable $\psi$ matches well with the reference solution obtained by \eqref{eq:ex1_psi} with $N = 2$. Fig. \ref{fig:1D_heaviside_S2_2} indicates that the microscopic variable $\phi$ also fits well with the analytical solution \eqref{eq:ex1_sol} at the microscopic velocity $\xi = \pm  \frac{1}{\sqrt{3}}$. 
 
 % time step be $\Delta t=\mathrm{CFL}\min\{\frac{\Delta x}{v}, \frac{3\Sigma\Delta x}{2v}\},\: \mathrm{CFL} = 0.4$. 
\begin{figure}[!hptb]
	\centering
    \subfloat[macroscopic scalar flux $\psi$]{
		\includegraphics[width = 0.4\textwidth]{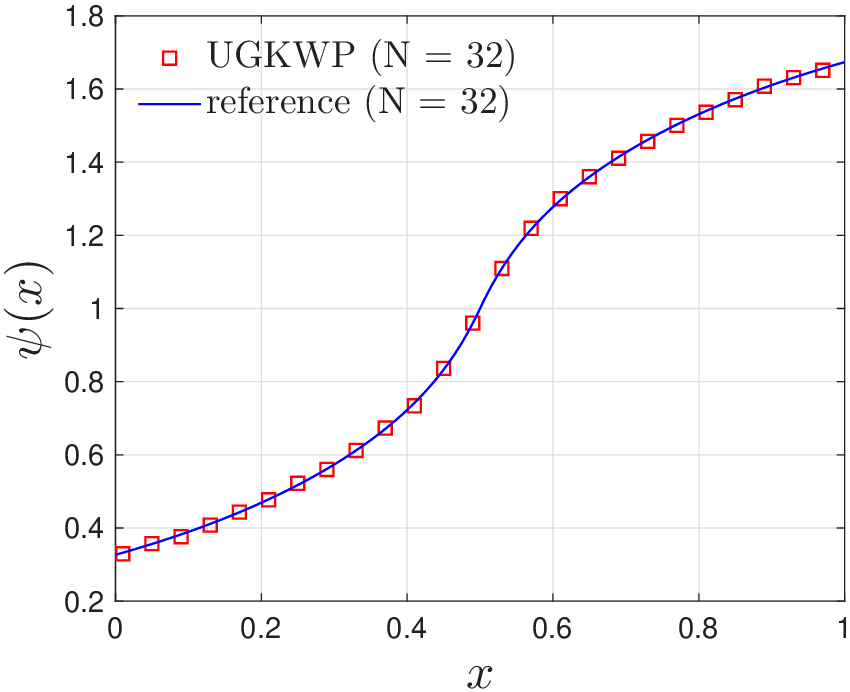}
        \label{fig:1D_heaviside_S32_1}
	}   \qquad 
        \subfloat[angular flux $\phi$]{               
		\includegraphics[width = 0.4\textwidth]{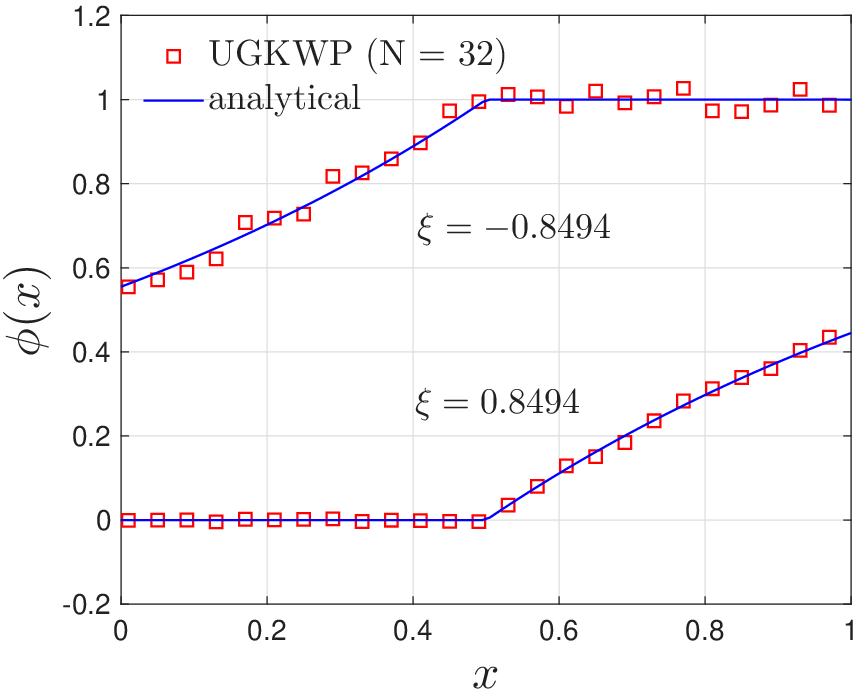}
        \label{fig:1D_heaviside_S32_2}
	}\\
	\subfloat[angular flux $\phi$]{               
		\includegraphics[width = 0.4\textwidth]{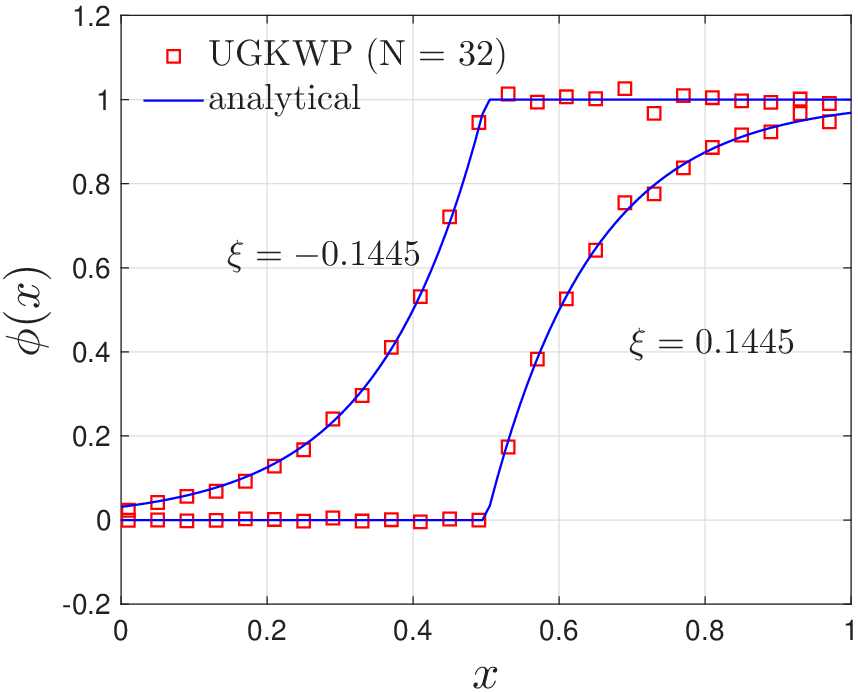}
        \label{fig:1D_heaviside_S32_3}
	} \qquad 
    \subfloat[macroscopic scalar flux $\psi$]{
		\includegraphics[width = 0.4\textwidth]{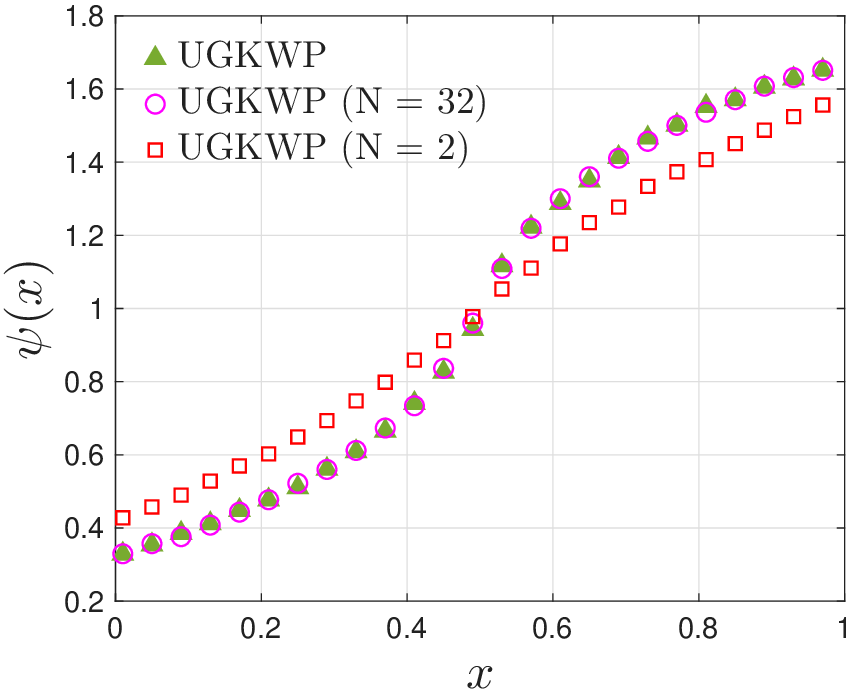}
        \label{fig:1D_heaviside_S32_4}
	} 
    \label{fig:1D_heavisede_S32}
    \caption{(1D Heaviside source problem in Sec. \ref{sec:1D_heaviside}) The comparison of the numerical results by UGKWP and the reference solution by $N = 32$.  (a) comparison of $\psi$, where the red square is the numerical solution obtained by UGKWP with $N = 32$, and the blue line is the reference solution obtained by \eqref{eq:ex1_psi} with $N = 32$;   
    (b) the numerical solution of $\phi$ at $\xi = \pm 0.8494$, where the red square is the numerical solution obtained by UGKWP with $N = 32$, and the blue line is the analytical solution \eqref{eq:ex1_sol}; 
    (c) the numerical solution of $\phi$ at $\xi = \pm 0.1445$, where the red square is the numerical solution obtained by UGKWP with $N = 32$, and the blue line is the analytical solution \eqref{eq:ex1_sol}; 
    (d) the numerical solution of $\psi$ by UGKWP, and UGKWP with $N = 2$ and $32$. 
    }  	
\end{figure}

Then we increase $N$ to $N =32$, and the numerical solution obtained by UGKWP with $N = 32$ is shown in Fig. \ref{fig:1D_heaviside_S32_1}. Fig. \ref{fig:1D_heaviside_S32_1} shows that the macroscopic variable $\psi(x)$ coincides well with the reference solution \eqref{eq:ex1_psi} with $N = 32$. The neuron angular flux $\phi(x, \xi)$ at $\xi = \pm 0.8494$ and $\pm 0.1445$ is plotted in Fig. \ref{fig:1D_heaviside_S32_2} and \ref{fig:1D_heaviside_S32_3}, where the numerical solution also fits well with the analytical solution. Finally, the numerical solution obtained by UGKWP with the distribution function of the microscopic velocity uniform distribution over the interval $[-1,1]$ is illustrated in Fig. \ref{fig:1D_heaviside_S32_4}, where $\psi(x)$ obtained by UGKWP with $N = 2$ and $N = 32$, is also plotted. We can find significant differences when $N = 2$ and $N = 32$. Moreover, the numerical solutions obtained with $N = 32$ and the normal UGKWP match well with each other, indicating the accuracy of the UGKWP method.

\subsubsection{1D isotropic scattering}
\label{sec:1D_isosca}
In this section, the 1D isotropic scattering problem for the single-group model is studied, the setting of which is as below 
\begin{equation}
\label{eq:ex2_para}
    \Sigma (x) = \Sigma_{s}(x) = \tilde{\Sigma}, \qquad \Sigma_{f}(x) = 0,\qquad q(x) = 0.5, \qquad  \phi(0,x,\xi) = 0,\qquad x \in [0, 1],
\end{equation} 
with the vacuum boundary condition utilized. This problem is also studied in \cite{bindra2012radiative, shuang2019parallel}. When $\Sigma \gg 1$, the source term $S$ will lead to strong numerical stiffness for the traditional single-group model, which may bring great numerical challenges. 

\begin{figure}[!hptb]
	\centering
	\subfloat[macroscopic scalar flux $\psi$]{
			\includegraphics[width = 0.28\textwidth]{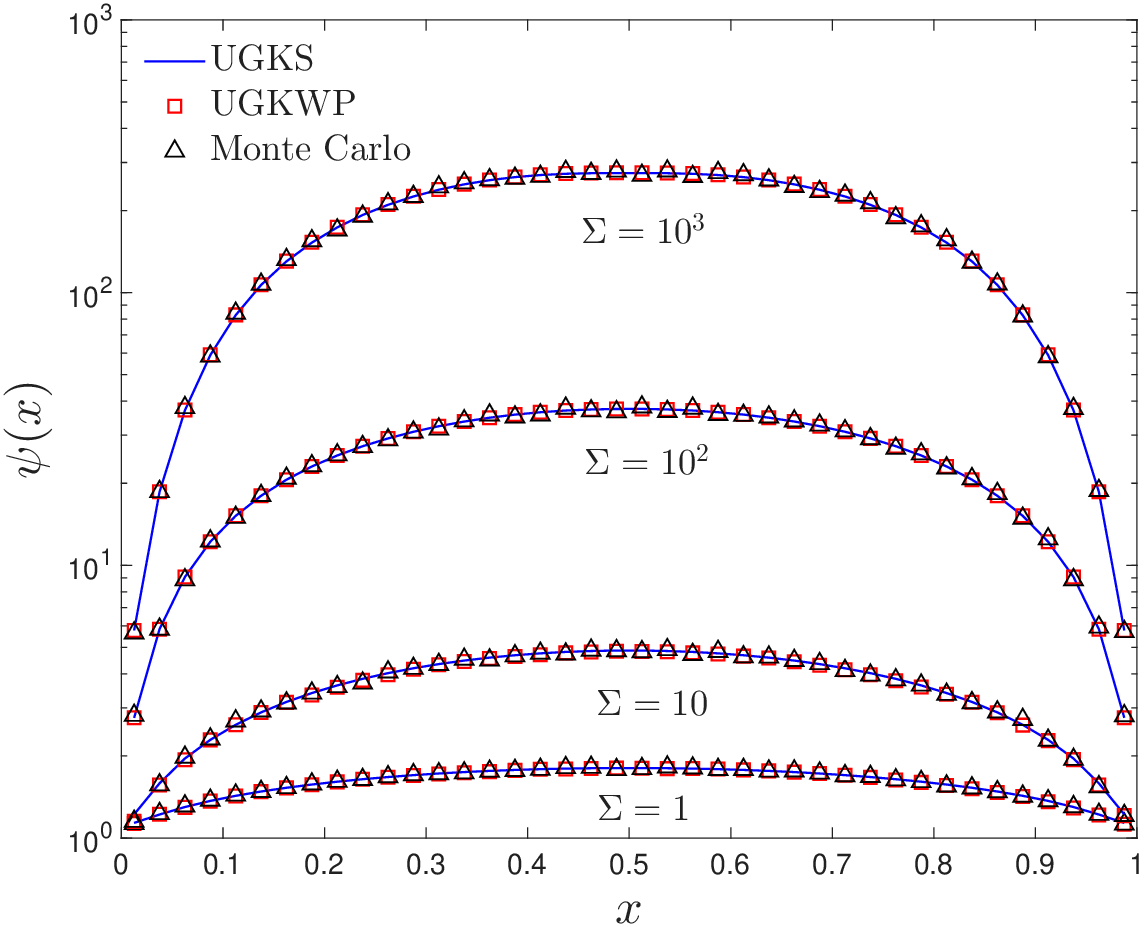}
	        \label{fig:1D_isosca_1}
		}\qquad  
	\subfloat[positive $\psi^{+}$]{
			\includegraphics[width = 0.28\textwidth]{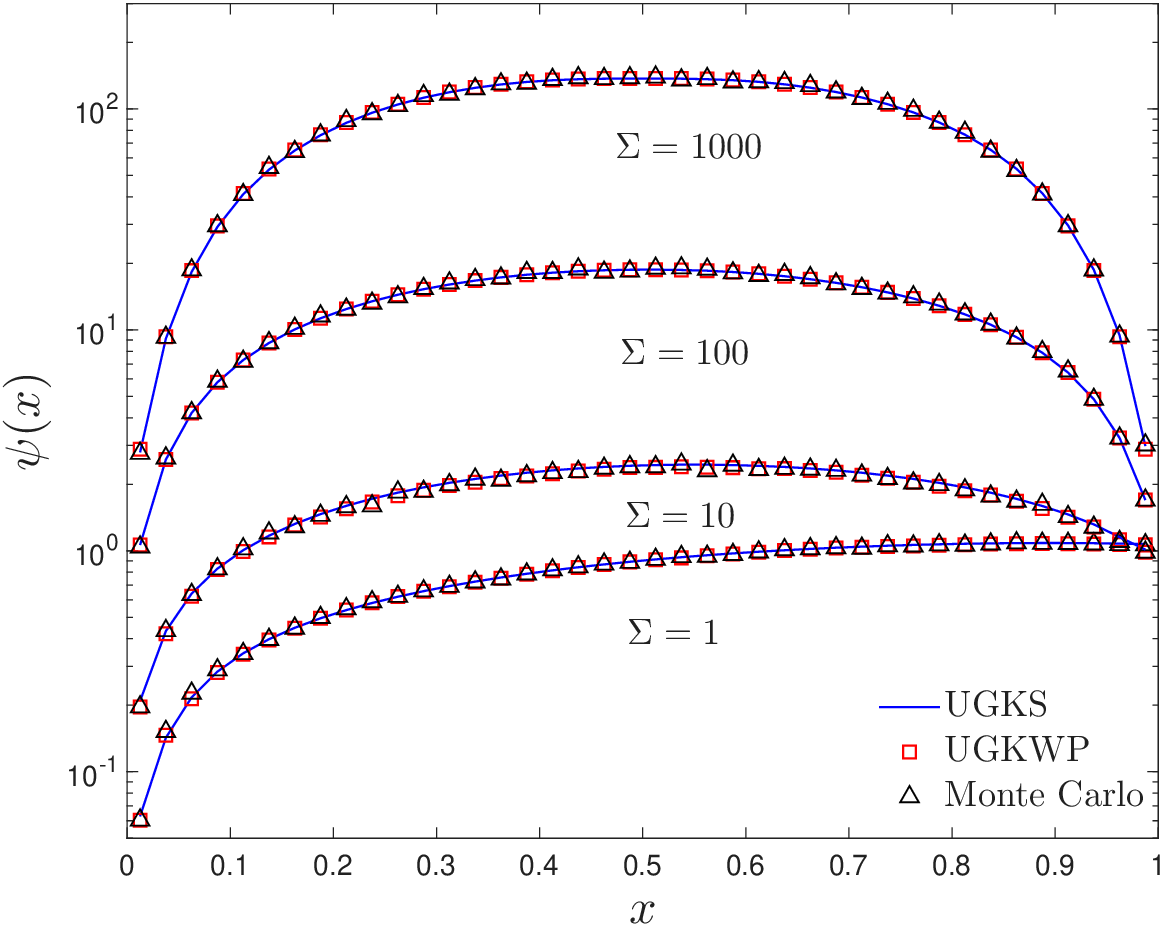}
			\label{fig:1D_isosca_2}
		}\qquad  
	\subfloat[negative $\psi^{-}$]{
			\includegraphics[width = 0.28\textwidth]{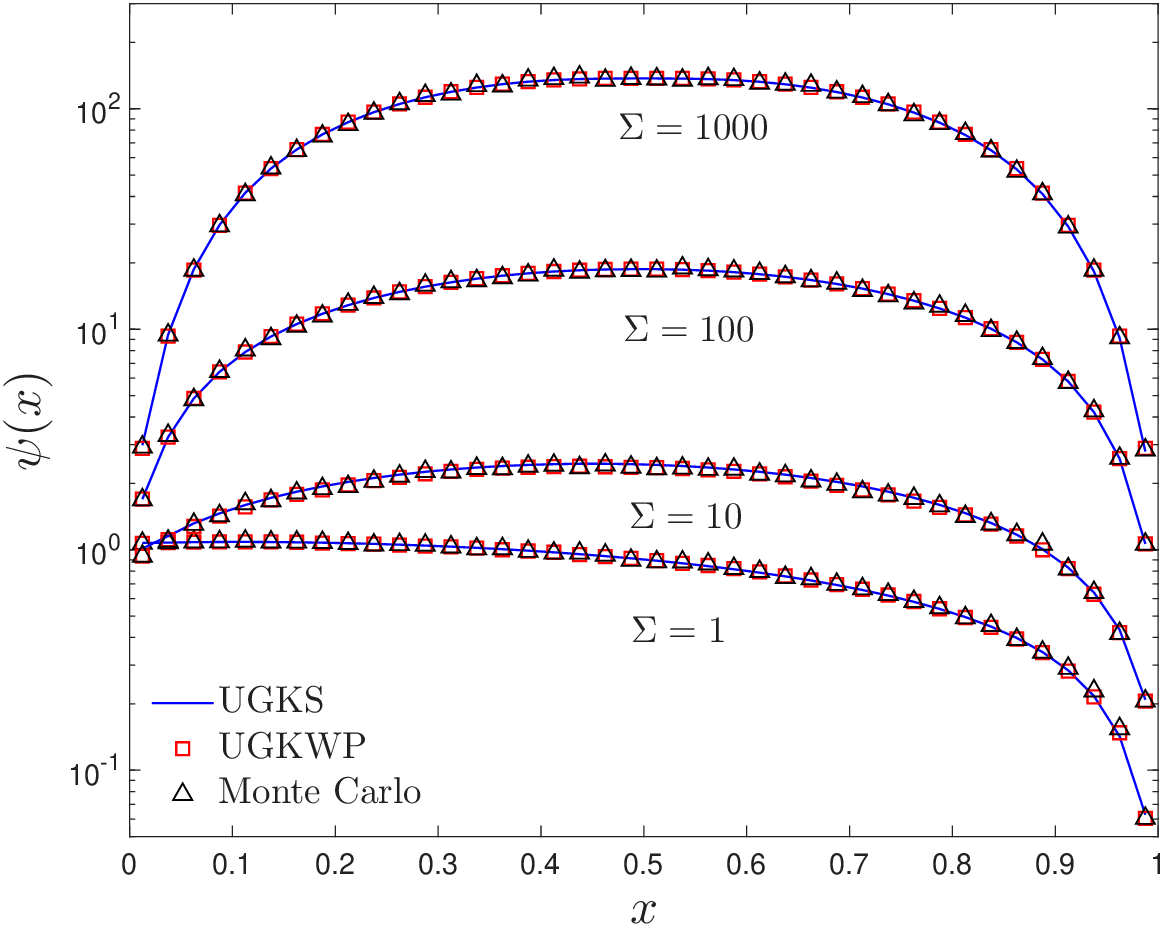}
			\label{fig:1D_isosca_3}
		}\qquad  
	\caption{(1D isotropic scattering problem in Sec. \ref{sec:1D_isosca}) The numerical solution of the macroscopic scalar flux $\psi$ of UGKWP, UGKS and Monte Carlo with different $\Sigma$. 
	    (a) $\psi$; (b) positive $\psi^{+}$; (c) negative $\psi^{-}$.
	    }  
\label{fig:1D_isosca}
\end{figure}

In the simulation, people are interested in the steady state of this problem. Similar as in \cite{shuang2019parallel}, the cases $\Sigma = 1, 10, 10^2$ and $10^3$ are studied. The spatial mesh size is set as $N_x = 40$, and the CFL number is $\mathrm{CFL} = 0.2$. The numerical solution of the macroscopic scalar flux $\psi$ as well as
\begin{equation}
    \label{eq:psi_NP}
    \psi^{+} = \int_{0}^1 \phi \dd \xi, \qquad \psi^{-} = \int_{-1}^0 \phi \dd \xi
\end{equation}
is shown in Fig. \ref{fig:1D_isosca}, where the numerical solution by the UGKS method, and the Monte Carlo method are also plotted. Here, the discrete ordinate method with $N= 32$ is utilized in UGKS. It reveals that the numerical solution is well consistent with the reference solution by the UGKS and the Monte Carlo method for all the $\Sigma$.

Moreover, with the increase of $\Sigma$, the solution of \eqref{eq:1D} will converge to the macroscopic limit. For UGKWP, with the increase of $\Sigma$, it gradually reduces to a discrete approximation of the macroscopic diffusion equation. In this regime, the contribution from particle transport becomes negligible, leading to significantly fewer particles being generated and tracked in each time step. To validate this property, the average computational cost for each time step with different $\Sigma$ is compared, which is shown in Tab. \ref{tab:1D_isosca}. It indicates that with the increase of $\Sigma$, this computational time decreases, indicating improved efficiency of the UGKWP method in the diffusive regime.
The numerical experiments of Tab.~\ref{tab:1D_isosca} and Tab.~\ref{tab:1D_unsteady} are performed on a personal computer equipped with an Intel Core i7-13700K CPU (3.40 GHz), 16 GB RAM.
% we compare the computational cost of the UGKWP method by advancing the same number of time steps under different values of $\Sigma$. The results in Tab. \ref{tab:1D_isosca} demonstrate that as $\Sigma$ increases, the system approaches the diffusive limit, where the UGKWP scheme gradually reduces to a discrete approximation of the macroscopic diffusion equation. In this regime, the contribution from particle transport becomes negligible, leading to significantly fewer particles being generated and tracked in each time step. Consequently, the overall computational cost per time step decreases, indicating improved efficiency of the UGKWP method in the diffusive regime.

\begin{table}[!htbp]
	\centering
	\def\arraystretch{1.5}
	{\footnotesize
		\begin{tabular}{c|c c c c}
			$\Sigma$ & 1 & 10 & 100 & 1000 \\
			\hline
			time & 6.758E-03 & 2.461E-03 & 1.504E-03 & 2.600E-05 \\
		\end{tabular}
		\caption{(1D isotropic scattering problem in Sec. \ref{sec:1D_isosca}) Average computational cost per time step for UGKWP with different total cross-sections $\Sigma$. (unit: \textit{seconds})}
		\label{tab:1D_isosca}
	}
\end{table}

\subsubsection{1D unsteady benchmark}
\label{sec:1D_unsteady}
\begin{figure}[!hptb]
	\centering
		\subfloat[macroscopic scalar flux $\psi$]{
			\includegraphics[width = 0.3\textwidth]{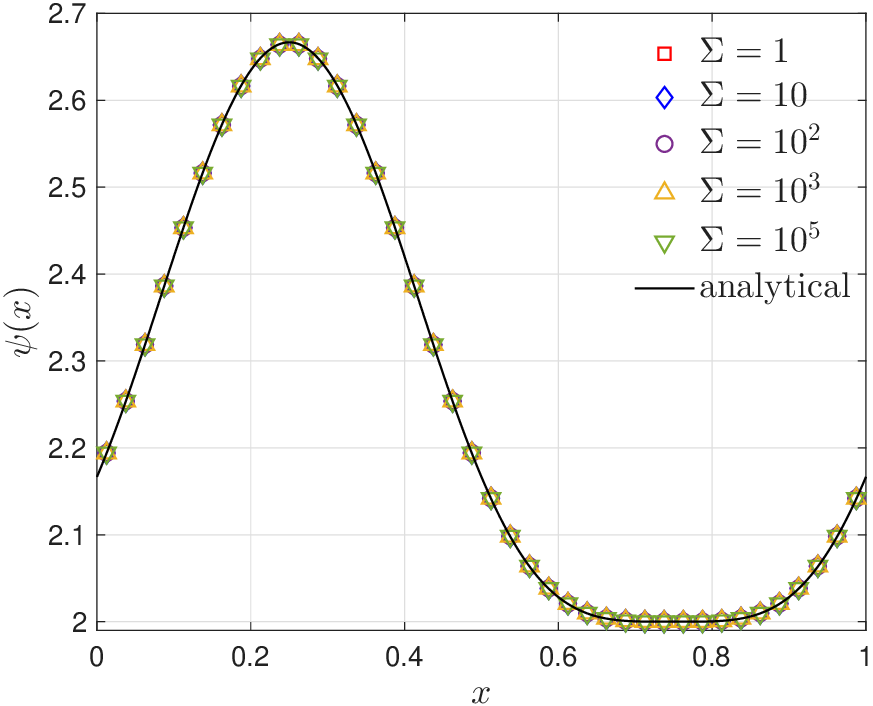}
			\label{fig:1D_unsteady_1}
		}\qquad
		\subfloat[angular flux $\phi$]{
			\includegraphics[width = 0.3\textwidth]{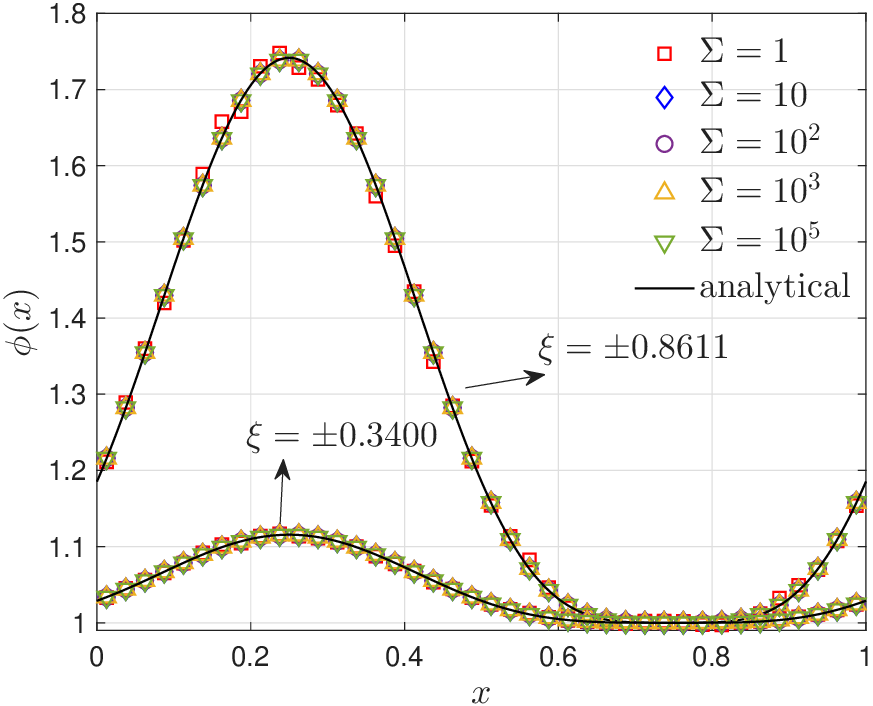}
			\label{fig:1D_unsteady_2}
		} 
		\subfloat[number of particls $N_p$ representing the distribution $\phi$]{
			\includegraphics[width = 0.3\textwidth]{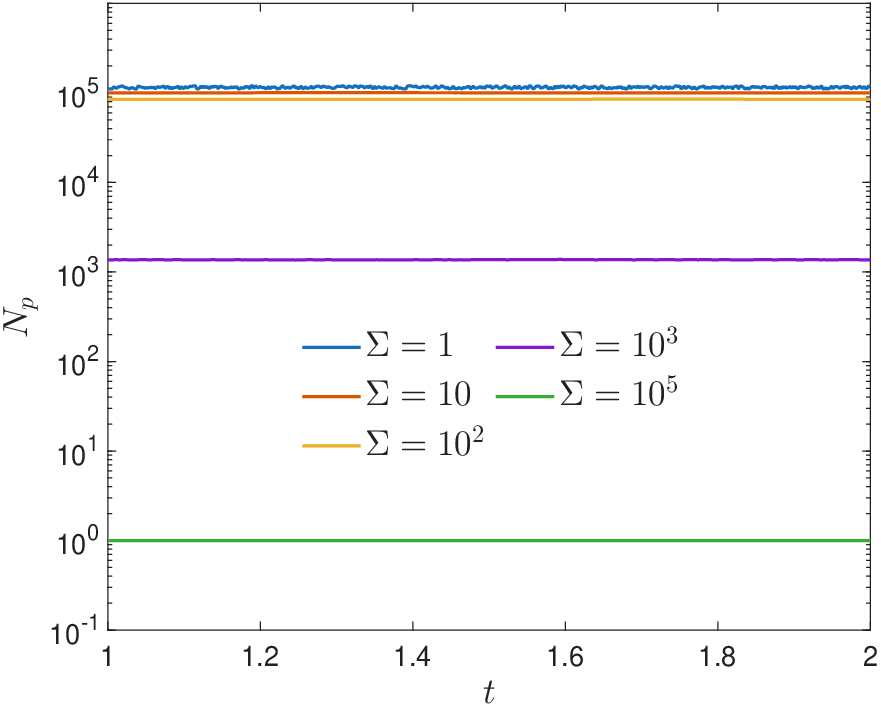}
			\label{fig:num_particle}
		}     
    \caption{(1D unsteady problem in Sec. \ref{sec:1D_unsteady}) The comparison of the numerical solution by UGKWP and the analytical solution. (a) comparison of $\psi$ at $t = 20.75$, where the black line is the analytical solution and other marker styles denote numerical solutions obtained by UGKWP under different total cross-sections $\Sigma$; (b) the numerical solution of $\phi$ at $t = 20.75$ with $\xi = \pm 0.8611$ and $\xi = \pm 0.3400$, where the black line is the analytical solution and other marker styles denote numerical solutions obtained by UGKWP with $N=4$ under different total cross-sections $\Sigma$; (c) number of particles representing the distribution $\phi$ in the system within one evolution period under different total cross-sections.}
    \label{fig:test1_3-varphi}
\end{figure}

In this section, the 1D unsteady benchmark problem for the single-group model is studied, the setting of which is as below
\begin{gather}
	\Sigma(x) = \Sigma_{s}(x) = \{1, 10, 10^{2}, 10^{3},10^{5}\},\quad \Sigma_f(x) = 0,\quad \phi(0, x, \xi) = \xi^{2}\mathrm{cos}^{4}(\pi x)+1, \qquad x\in [0, 1], \\
	 q(t, x,\xi) = -4\pi\xi^{2}(1+\xi)\cos^{3}\left(\pi(x+t)\right)\sin \left(\pi(x+t)\right)+ \left(\xi^{2}\Sigma - \frac{\Sigma_{s}}{3}\right)\cos^{4}\left(\pi(x+t)\right)+(\Sigma - \Sigma_s),
\end{gather}
with the periodic boundary condition utilized. 
There exists an exact solution to this problem as 
\begin{equation}
	\label{eq:test1_3_phi}
	\phi(x,\xi, t) = \xi^{2}\mathrm{cos}^{4}\left(\pi(x+t)\right) + 1.
\end{equation}
In the simulation, the mesh size is set as $N_x = 40$, and the CFL number is $\mathrm{CFL} = 0.3$. To show the stability of UGKWP, the long-time behavior of the macroscopic scalar flux $\psi$ at $t = 20.75$ is plotted in Fig. \ref{fig:1D_unsteady_1}, and it shows that for each $\Sigma$, the numerical solution matches well with the analytical solution. Moreover, the numerical solution of the angular flux $\phi$ with $\xi = \pm 0.8611$ and $\xi = \pm0.3400$ at $t = 20.75$ is presented in Fig. \ref{fig:1D_unsteady_2}, which shows that the numerical solution also fits well with the analytical solution even for $\Sigma = 10^5$, indicating the effectiveness of UGKWP for large $\Sigma$. Furthermore, with the increase of $\Sigma$, the solution to \eqref{eq:1D} asymptotically approaches the macroscopic diffusion limit. In this regime, the UGKWP method effectively reduces to a discrete diffusion solver. As a result, the number of particles representing the distribution $\phi$ within the system decreases substantially, as illustrated in Fig. \ref{fig:num_particle}, leading to a marked reduction in computational time, as reported in Tab. \ref{tab:1D_unsteady}.

\begin{table}[!htbp]
	\centering
	\def\arraystretch{1.5}
	{\footnotesize
		\begin{tabular}{c|c c c c c}
			$\Sigma$ & $1$ & $10$ & $10^{2}$ & $10^{3}$ & $ 10^{5}$ \\
			\hline
			time & 45.28 & 28.33 & 15.73 & 1.76 & 0.21 \\
		\end{tabular}
		\caption{(1D unsteady benchmark problem in Sec~\ref{sec:1D_unsteady}) The total computation time required by UGKWP for $t = 1$ under different total cross-sections $\Sigma$ (unit: \textit{seconds}).}
		\label{tab:1D_unsteady}
	}
\end{table}

\subsection{3D single-group problems}
\label{sec:3D_single_group}
In this section, the two 3D single-group problems are tested, including the three-cubic problem and the top-hat problem, which are also studied in \cite{kobayashi20013d, shuang2019parallel}. In each issue, the computational region consists of three materials with different cross-sections, and two distinct situations are explored. The first is a system of the pure absorption problem with a void region, for which the analytical solution can be obtained through numerical integration. In the second situation, the pure absorption is replaced by a material whose scattering cross-section is $50\%$ of the total cross-section. The related parameters for these two situations are listed in Tab. \ref{tab:3D_sg_cross_section}. 

% \begin{table}[!htbp]
% 	\centering
% 	\def\arraystretch{1.5}
% 	{\footnotesize
% 		\begin{tabular}{c|cccc}
% 			region & $q$ & $\Sigma$ & $\Sigma_{s}$ (Abs.) & $\Sigma_{s}$ (Sca.) \\
% 			\hline
% 			1 & 1 & 0.1 & 0 & 0.05 \\
% 			2 & 0 & $10^{-4}$ & 0 & $0.5 \times 10^{-4}$ \\
% 			3 & 0 & 0.1 & 0 & 0.05 \\
% 		\end{tabular}
% 		\caption{(3D single-group benchmark problem in Sec.~\ref{sec:3D_single_group})
% 			Cross-section parameters and source $q$ for the 3D single-group problem.
% 			Here, Abs. denotes the pure absorption case, and Sca. denotes the scattering case.}
% 		\label{tab:3D_sg_cross_section}
% 	}
% \end{table}
    \begin{table}[!htbp]
      \centering
      \def\arraystretch{1.5}
     {\footnotesize
      \begin{tabular}{c|c|c|c|c}
         \hline \text { Region } & $q$ & $
             \Sigma $& $
            \Sigma_{s} \text { (Abs.) }$& $
            \Sigma_{s} \text { (Sca.) }$ \\
            \hline 1 & 1 & 0.1 & 0 & 0.05 \\
            2 & 0 & $10^{-4} $& 0 & $0.5 \times 10^{-4}$ \\
            3 & 0 & 0.1 & 0 & 0.05 \\
            \hline
      \end{tabular}
      }
      \caption{(3D single-group benchmark problem in Sec. \ref{sec:3D_single_group})
      The cross-section parameters and source $q$ for the 3D single-group problems.
      Here, Abs. is short for the pure absorption case, and Sca. is short for the scattering case. 
      }
      \label{tab:3D_sg_cross_section}
    \end{table}

\subsubsection{The three-cubic problem}
\label{sec:3D_single_group_1}
\begin{figure}[!hptb]
	\centering
	\subfloat[3D view]{
    \label{fig:3D_three_cubic_1}
		\includegraphics[width = 0.3\textwidth]{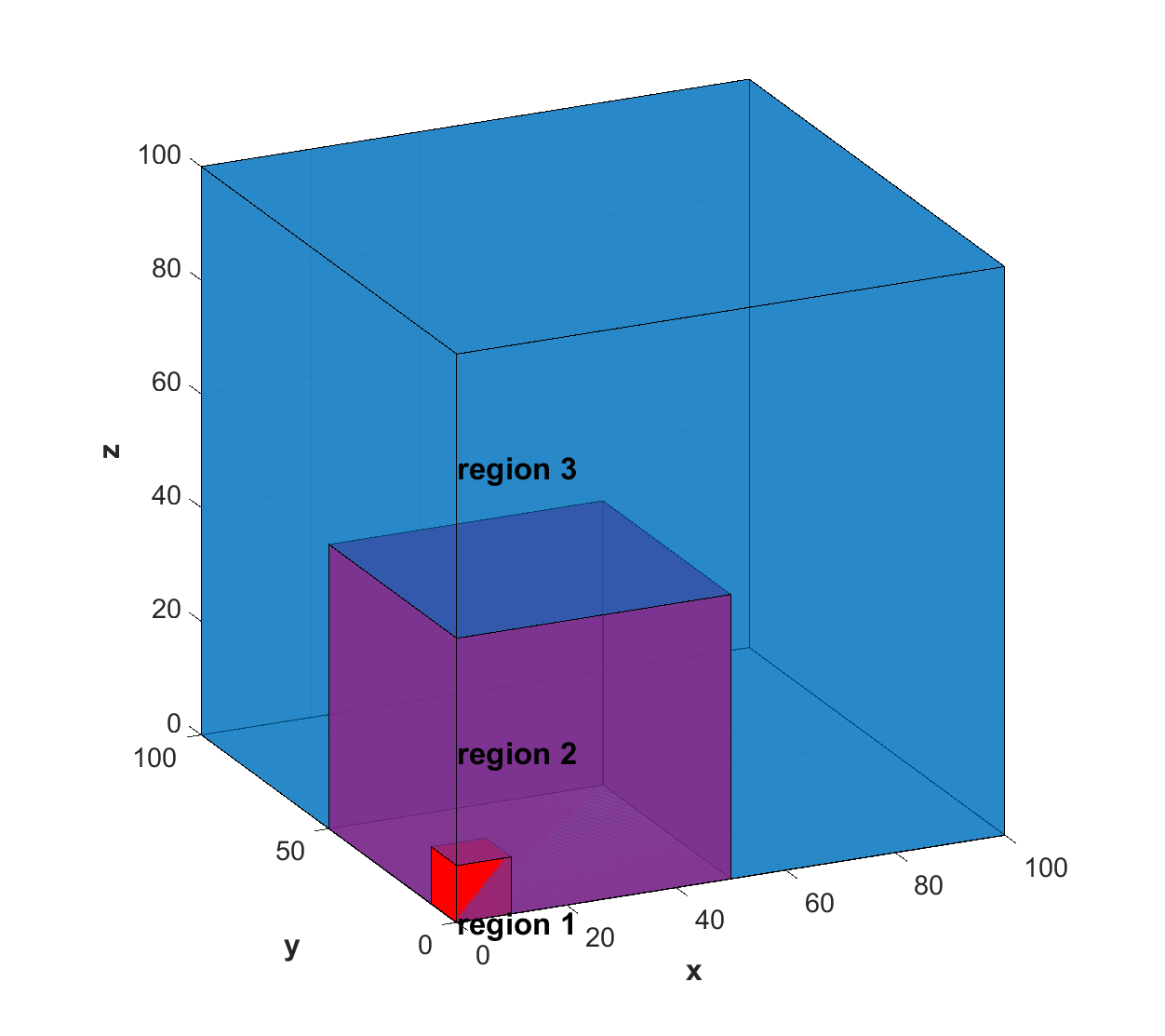}
	}
	\subfloat[boundary condition]{      
        \label{fig:3D_three_cubic_2}
		\includegraphics[width = 0.28\textwidth]{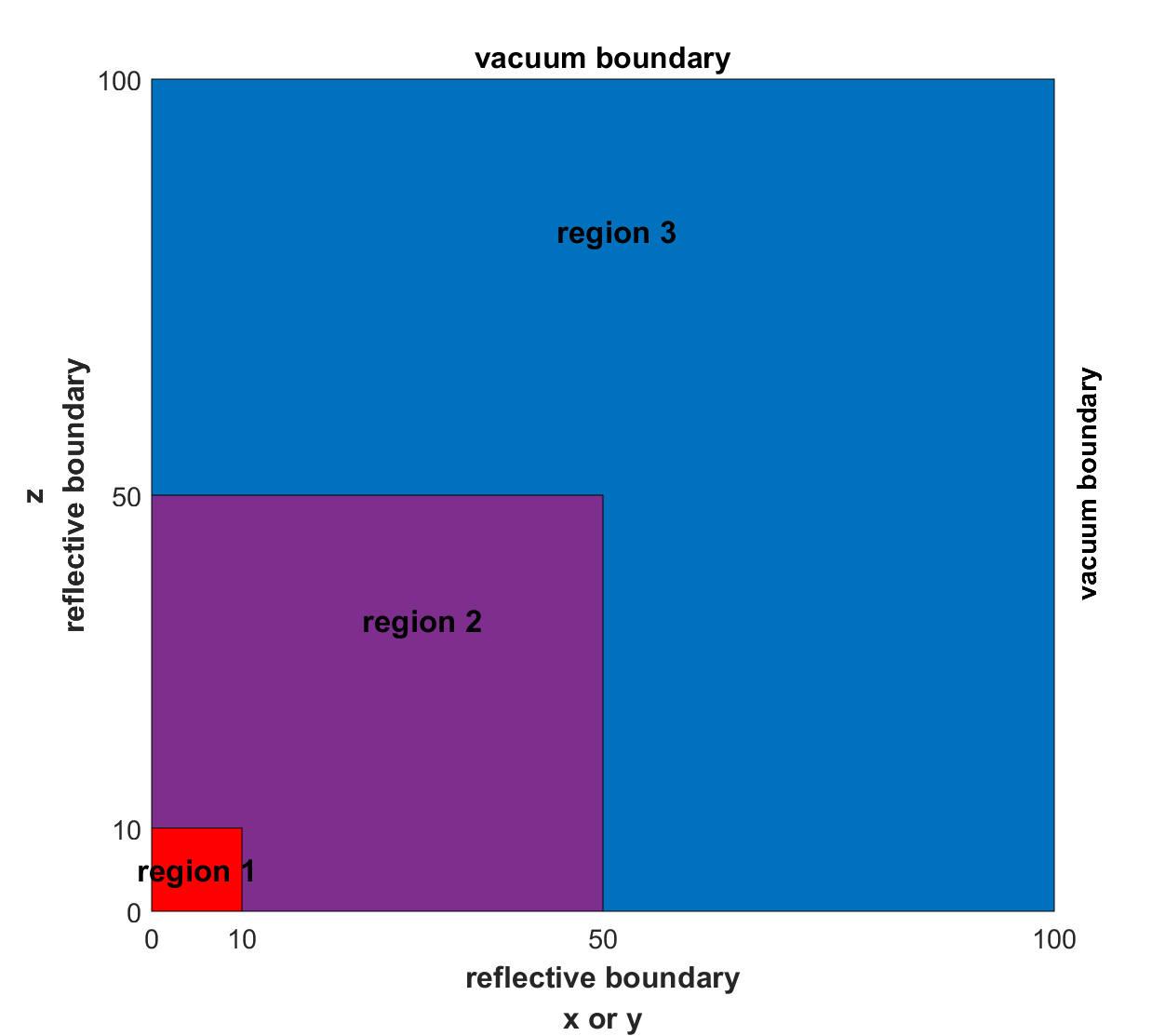}
	}       
	\caption{(The three cubic problem in Sec. \ref{sec:3D_single_group_1}) The geometry profile of the three-cubic problem. (a) the 3D view; (b) the boundary condition.}
    	\label{fig:3D_three_cubic}
\end{figure}

In this section, the three-cubic problem is studied, where the computational region consists of three cuboids with different materials. The geometry profile is a cube of edge length $100$ as shown in Fig. \ref{fig:3D_three_cubic_1}. The boundary condition is shown in Fig. \ref{fig:3D_three_cubic_2}, where the reflective boundary is utilized in the $\rm{XOY}, \rm{XOZ}$ and $\rm{YOZ}$ plane, with a vacuum boundary adopted in other planes. 

\begin{figure}[!hptb]
	\centering
	\subfloat[$\psi (x = z = 5)$]{
        \label{fig:3D_1_psi_abs_1}
		\includegraphics[width = 0.3\textwidth]{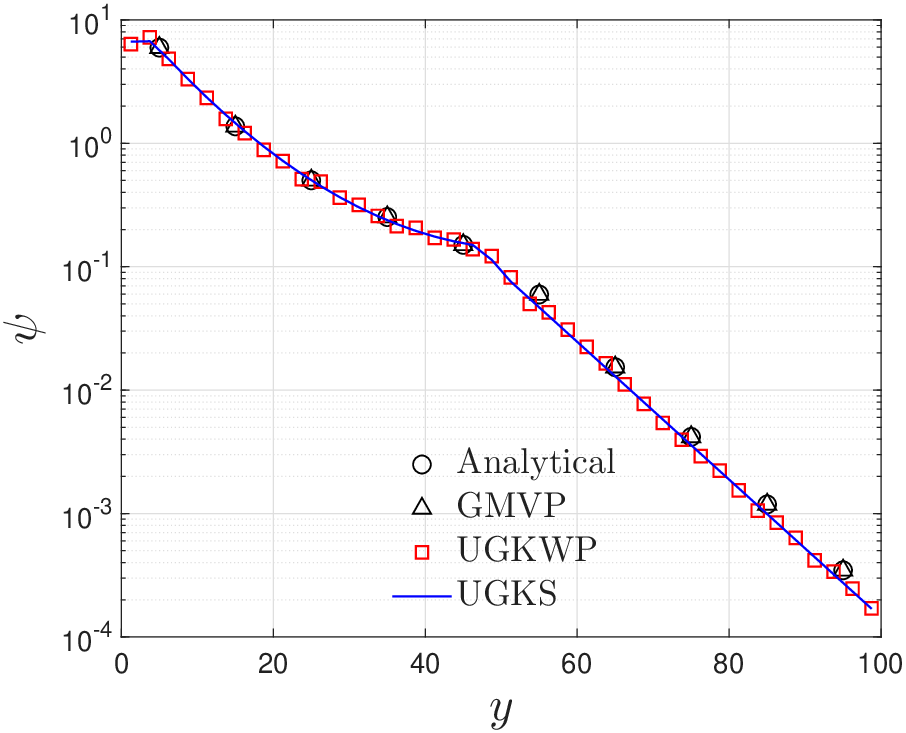}
	} \hfill
	\subfloat[$\psi (x = y = z)$]{  
        \label{fig:3D_1_psi_abs_2}
		\includegraphics[width = 0.3\textwidth]{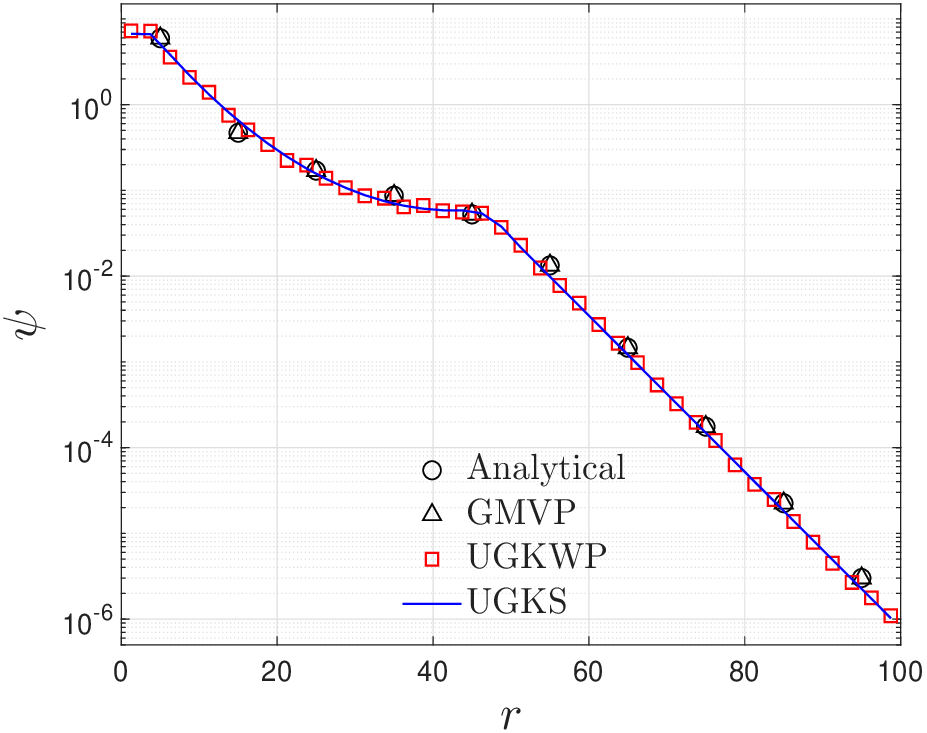}
	}      \hfill
    \subfloat[$\psi (y = 55, z = 5)$]{   
        \label{fig:3D_1_psi_abs_3}
		\includegraphics[width = 0.3\textwidth]{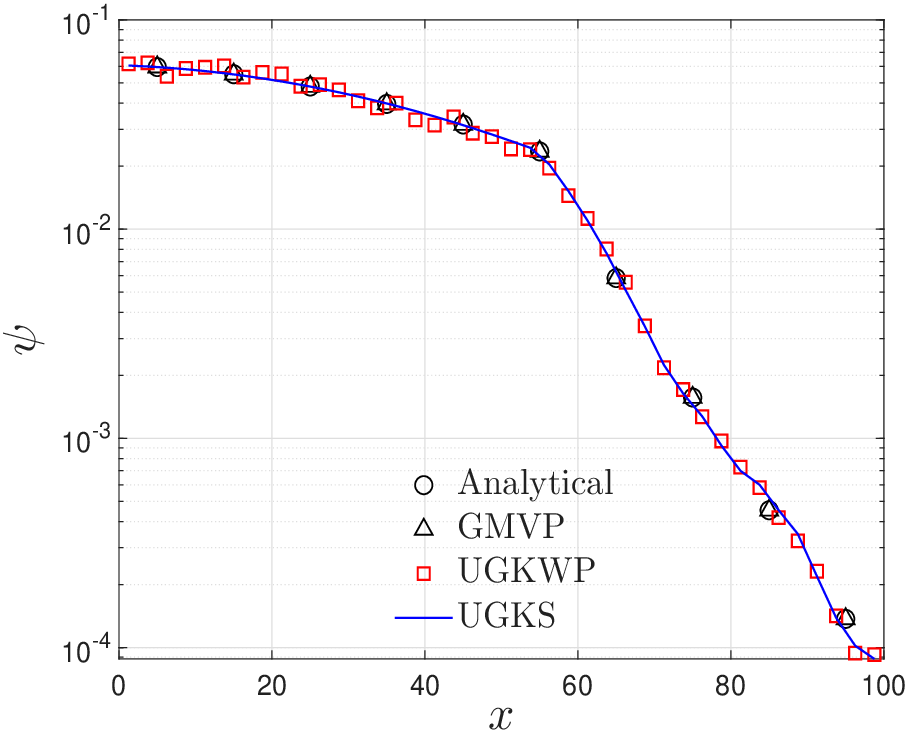}
	}     \hfill
	\caption{(The three-cubic problem in Sec. \ref{sec:3D_single_group_1}) The macroscopic scalar flux $\psi$ of the three-cubic problem for the pure absorption case. Here, the square is the analytical solution, the triangle is the reference solution by GMVP, the blue line is the reference solution by UGKS, and the star line is the numerical solution by UGKWP.  (a) $x = z =5$; (b) $x = y = z$; (c) $y = 55, z= 5$.}
    \label{fig:3D_1_psi_abs}
\end{figure}

In the simulation, people are interested in the steady-state solution. The mesh size is set as $N_x = N_y = N_z = 40$ and the CFL number is $\mathrm{CFL} = 0.2$. The macroscopic flux 
\begin{equation}
    \label{eq:3D_psi}
    \psi(x, y, z) = \frac{1}{4\pi}\int_{4\pi}\phi(x, y, z, \Omega) \dd \Omega
\end{equation}
is studied. The numerical solution of $\psi$ at $x = z = 5$, $x = y = z$ and $y = 55, z = 5$ for the pure absorption problem is plotted in Fig. \ref{fig:3D_1_psi_abs}. For the line $x = z = 5$, which is close to the source region when $y$ is small, the macroscopic flux $\psi$ is large at first and becomes smaller and smaller with the increase of $y$. For the line $x = y = z$, which is a space diagonal of the cube, the macroscopic flux $\psi$ is also decreasing with the increase of $r$, where $r$ is the distance from the origin. For the line $y = 55, z = 5$, which is far from the source region, the macroscopic flux $\psi$ is smaller compared to the other two lines, and still becomes smaller and smaller with the increase of $x$. The steady-state of $\psi$ at these three lines is illustrated in Fig. \ref{fig:3D_1_psi_abs_1}, \ref{fig:3D_1_psi_abs_2} and \ref{fig:3D_1_psi_abs_3}, respectively, where the analytical solution, the numerical solution obtained by GMVP in \cite{kobayashi20013d, mori1994mvp}, and the numerical solution by UGKS (G2-S16) in \cite{shuang2019parallel} is also presented as the reference solution. Fig. \ref{fig:3D_1_psi_abs} shows that for the three lines, the numerical solutions match well with the reference solution for all three lines, even for the line $y = 55, z = 5$, which is far from the source region. 

\begin{figure}[!hptb]
	\centering
	\subfloat[$\psi (x = z = 5)$]{
		\includegraphics[width = 0.3\textwidth]{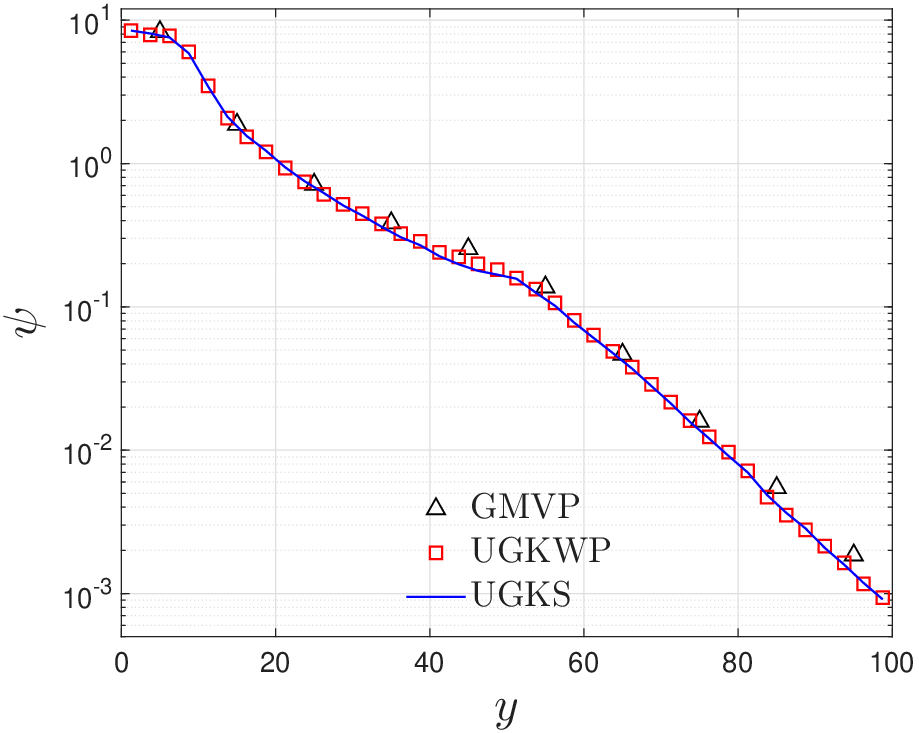}
	} \hfill
	\subfloat[$\psi (x = y = z)$]{               
		\includegraphics[width = 0.3\textwidth]{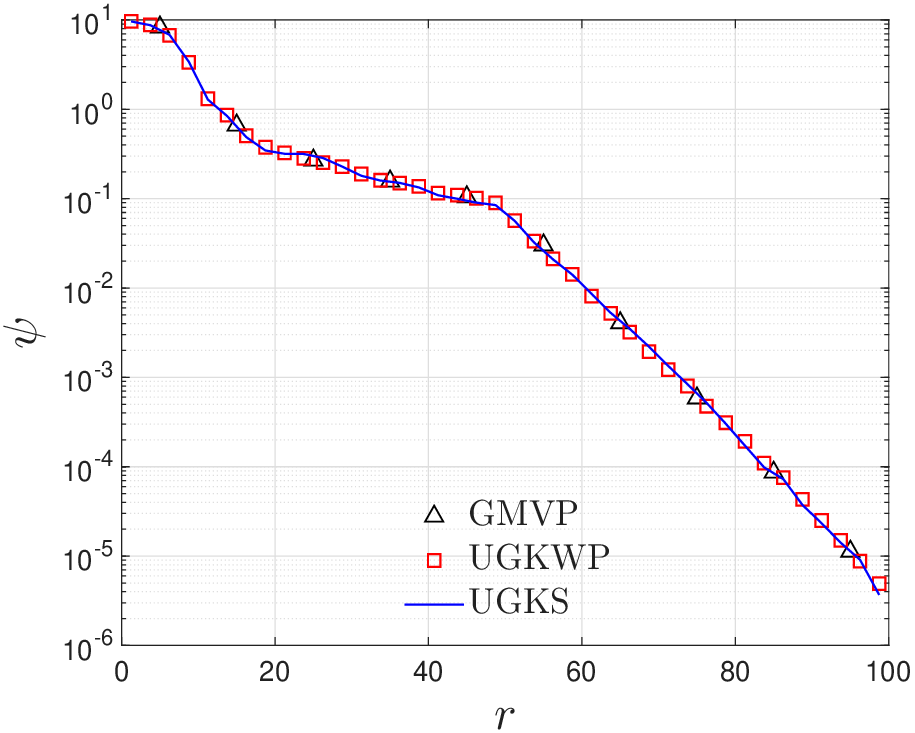}
	}      \hfill
    \subfloat[$\psi (y = 55, z = 5)$]{               
		\includegraphics[width = 0.3\textwidth]{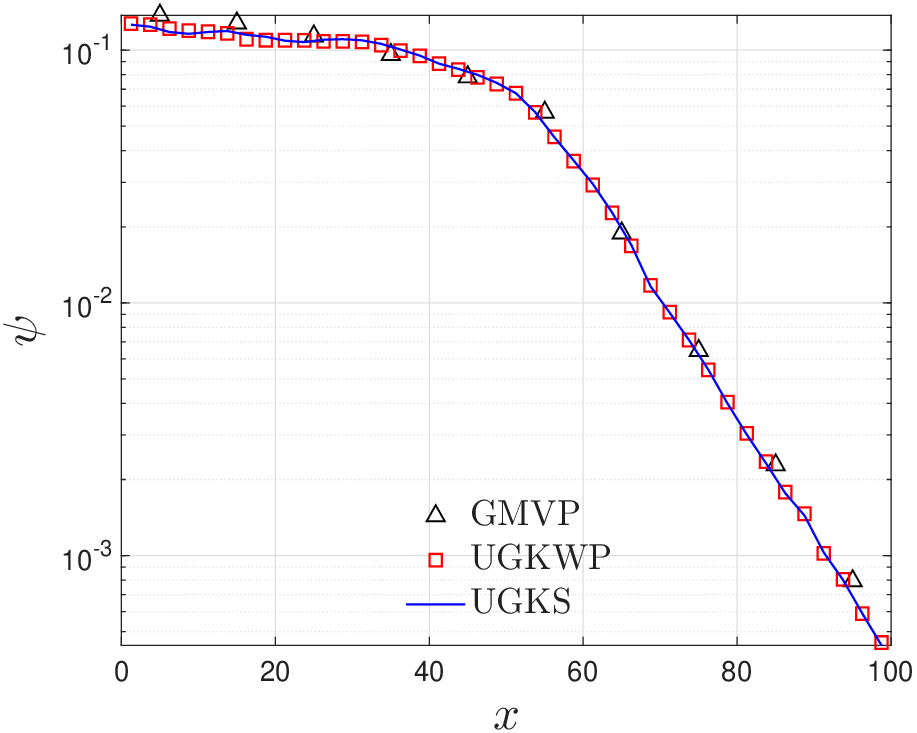}
	}     \hfill
	\caption{(The three-cubic problem in Sec. \ref{sec:3D_single_group_1}) The macroscopic scalar flux $\psi$ of the three-cubic problem for the scattering case. Here, the triangle is the reference solution by GMVP, the blue line is the reference solution by UGKS, and the star line is the numerical solution by UGKWP.  (a) $x = z = 5$; (b) $x = y = z$; (c) $y = 55, z= 5$.}
    \label{fig:3D_1_psi_sca}
\end{figure}

Then the numerical solution $\psi$ for the scattering case along the same lines is plotted in Fig. \ref{fig:3D_1_psi_sca}. Here, the numerical settings are the same as in the pure absorption case. The reference solution obtained by GMVP and UGKS (G2-S16) is also illustrated in Fig. \ref{fig:3D_1_psi_sca}. It shows that the numerical solution agrees well with the reference solution along all three lines.

\subsubsection{The top-hat problem}
\label{sec:3D_single_group_2}
\begin{figure}[!hptb]
	\centering
	\subfloat[3D view]{
		\includegraphics[width = 0.24\textwidth]{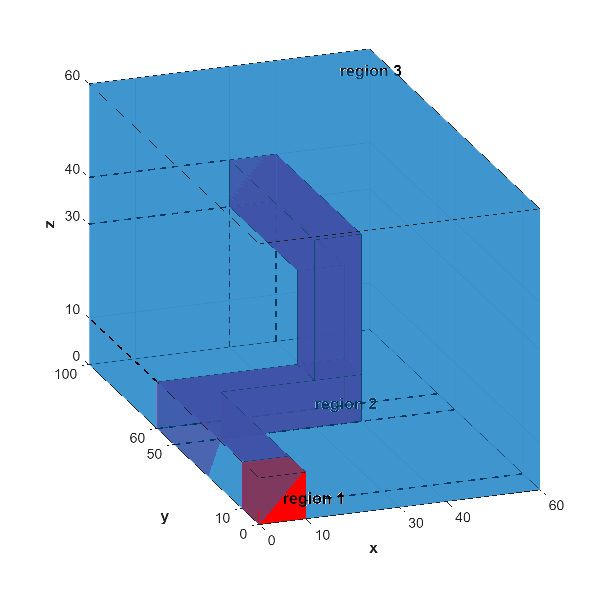}
	} 
	\subfloat[view on plane XOY]{               
		\includegraphics[width = 0.24\textwidth]{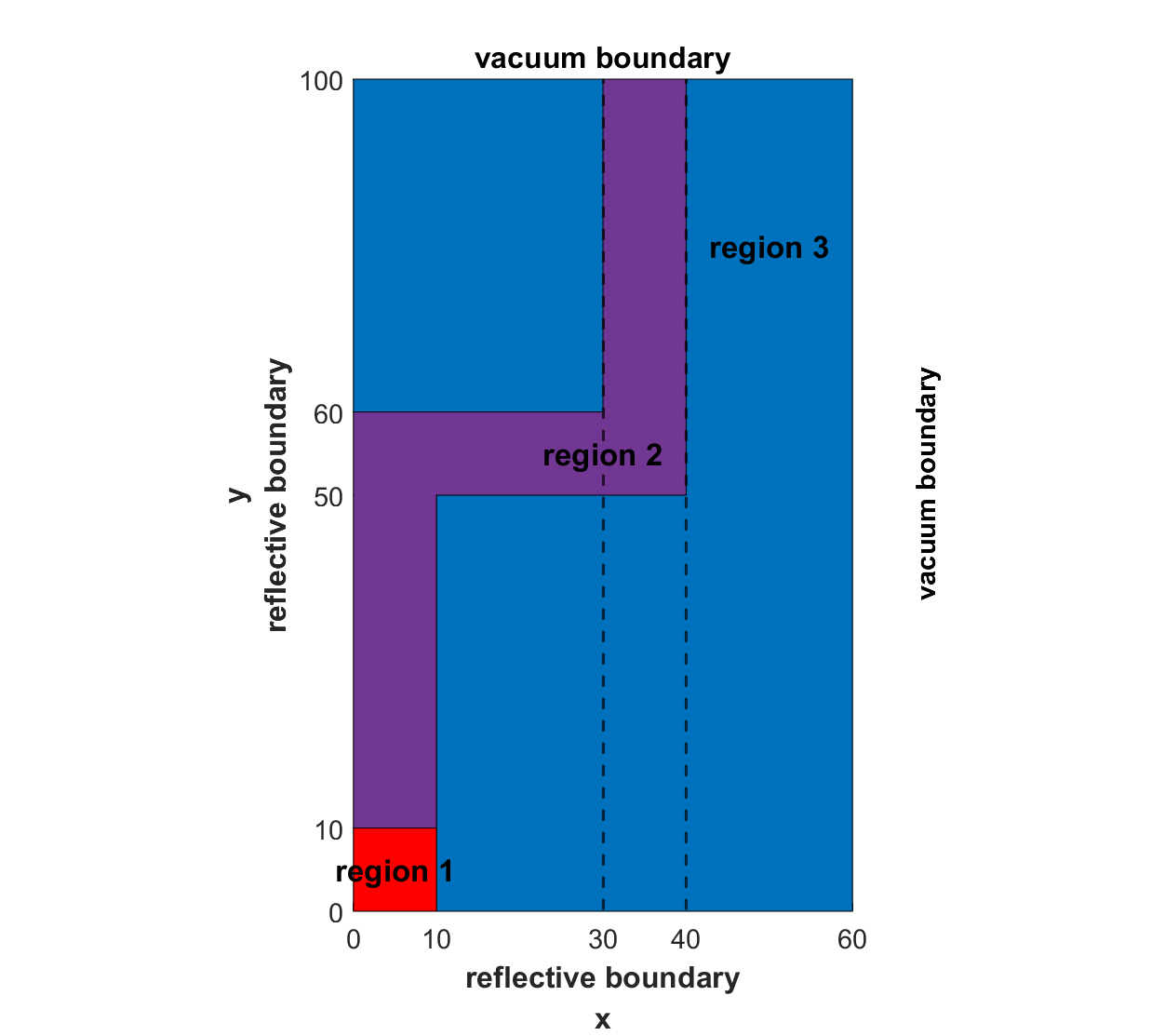}
	} 
        \subfloat[view on plane XOZ]{               
		\includegraphics[width = 0.24\textwidth]{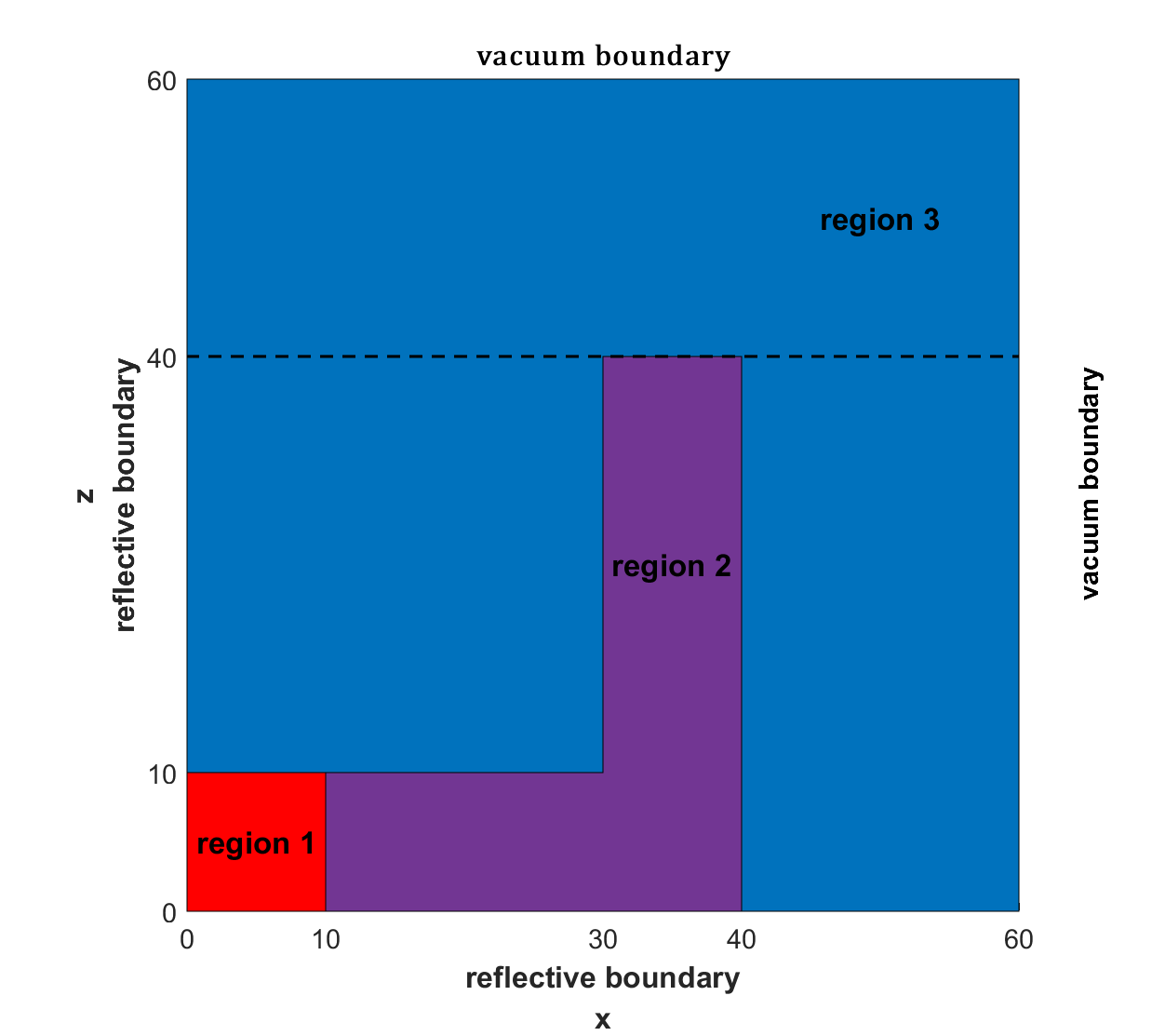}
	} 
        \subfloat[view on plane YOZ]{               
		\includegraphics[width = 0.24\textwidth]{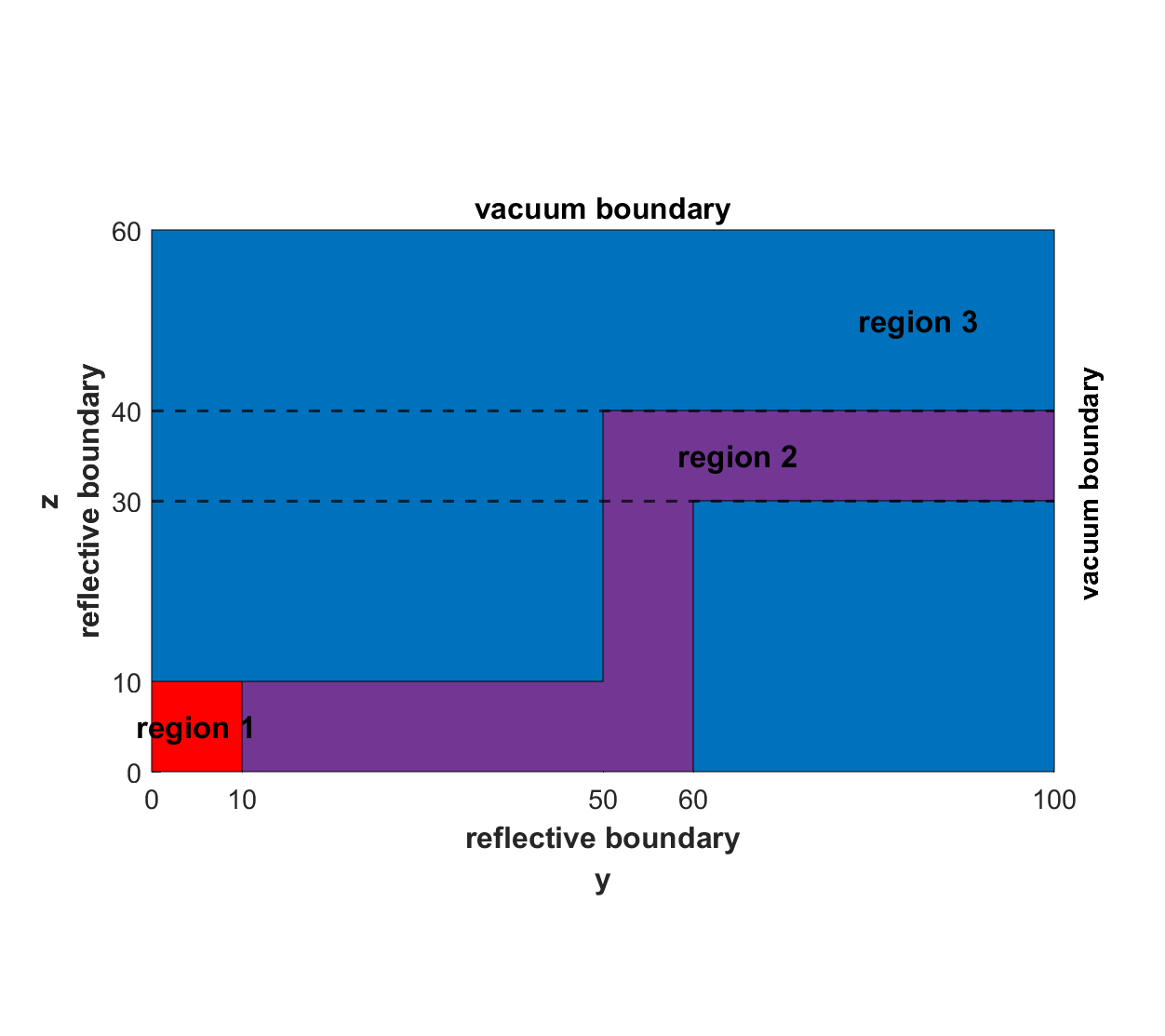}
	}    
	\caption{(The top-hat problem in Sec. \ref{sec:3D_single_group_2}) The geometry profile of the top-hat problem. (a) the 3D view; (b) view on plane XOY; (c) view on plane XOZ; (d) view on plane YOZ.}  
	\label{fig:3D_tophat}
\end{figure}

In this section, the 3D top-hat problem for the single-group model is studied. The geometry profile is a $60\times 100 \times 50$ cubic as illustrated in Fig. \ref{fig:3D_tophat}, where the boundary conditions are also plotted. Here, the reflective boundary condition is applied to the planes XOY, XOZ, and YOZ, while the vacuum boundary condition is used on the other planes. 

\begin{figure}[!hptb]
	\centering
	\subfloat[$\psi (x = z = 5)$]{
		\includegraphics[width = 0.3\textwidth]{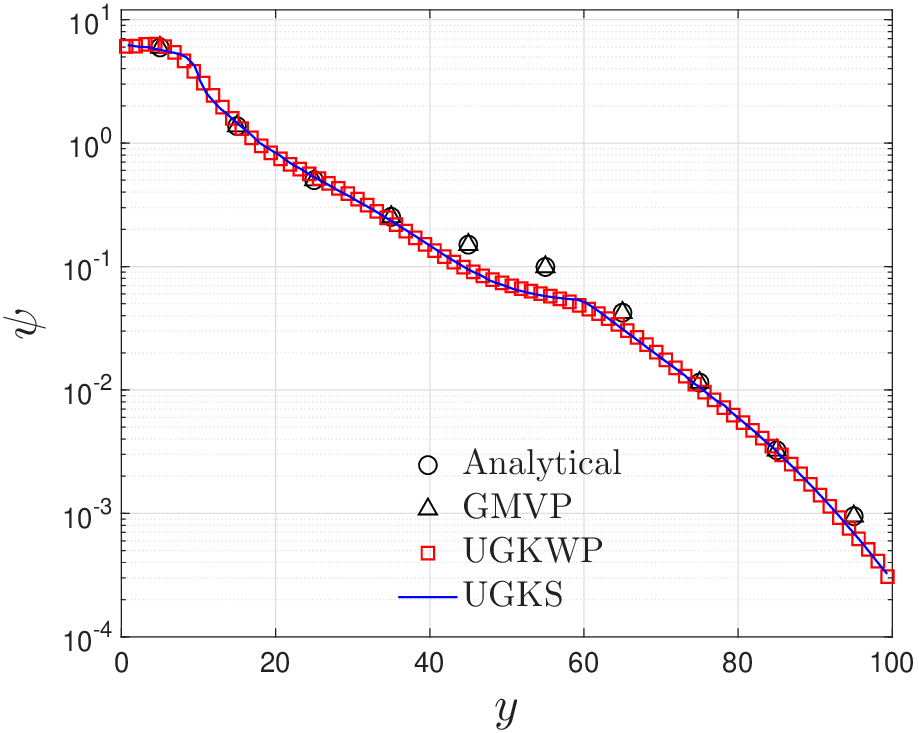}
	}\hfill
	\subfloat[$\psi (y = 55, z = 5)$]{               
		\includegraphics[width = 0.3\textwidth]{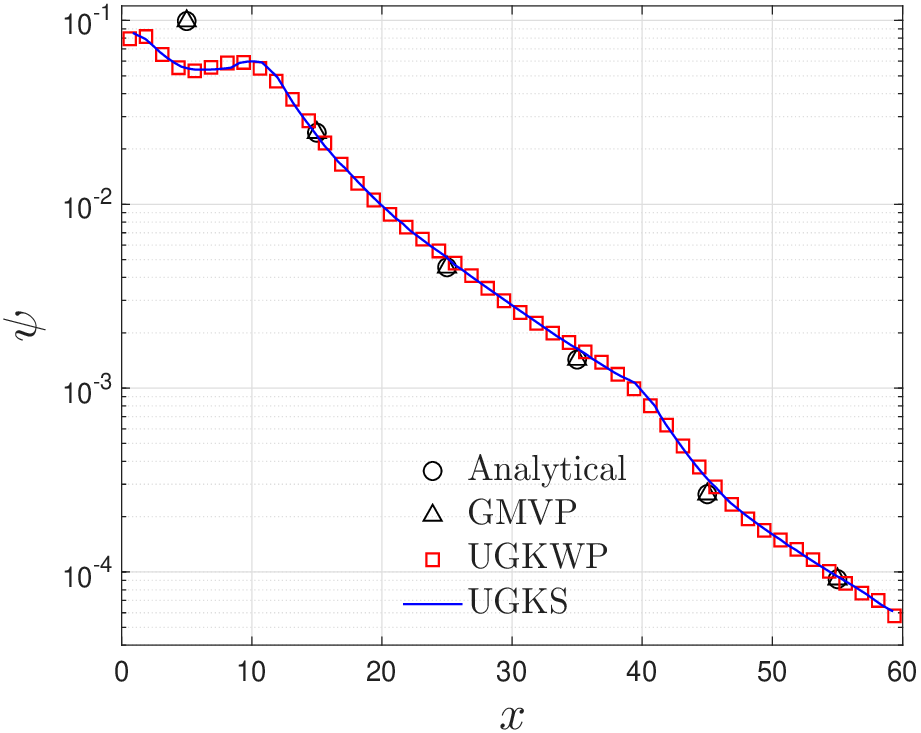}
	} \hfill
    \subfloat[$\psi (y = 95, z = 5)$]{
		\includegraphics[width = 0.3\textwidth]{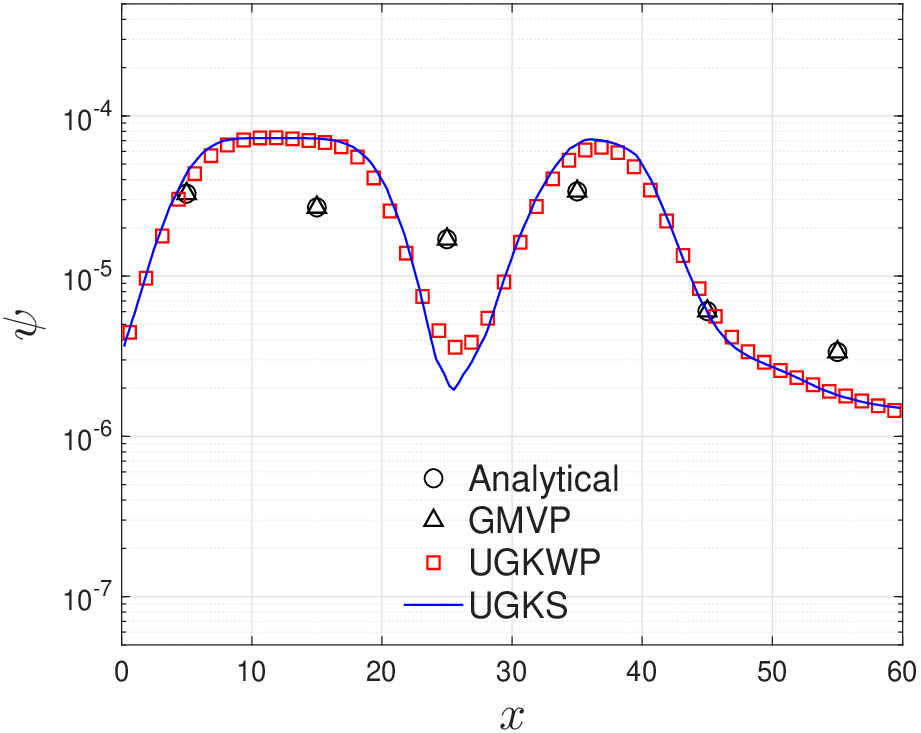}
	}
	\caption{(The top-hat problem in Sec. \ref{sec:3D_single_group_2}) The macroscopic scalar flux $\psi$ of the top-hat problem for the pure absorption case. Here, the triangle is the analytical solution, the square is the reference solution by GMVP, the blue line is the reference solution by UGKS, and the star line is the numerical solution by UGKWP.  (a) $x = z =5$; (b) $y = 55, z = 5$; (c) $y = 95, z= 5$.}  
	\label{fig:3D_tophat_abs}
\end{figure}

For this problem, people are also interested in the steady-state solution. In the simulation, the uniform orthogonal grid $\Delta x = \Delta y = \Delta z = 1.25$ is utilized and the CFL number is $\mathrm{CFL} = 0.2$.. The numerical solution for the pure absorption case is plotted in Fig. \ref{fig:3D_tophat_abs}, where the lines are $x = z = 5$, $y = 55, z = 5$, and $y = 95, z = 5$. The reference solution given by the analytical solution, GMVP method in \cite{kobayashi20013d, mori1994mvp}, UGKS (G2-S16) in \cite{shuang2019parallel} is also presented in Fig. \ref{fig:3D_tophat_abs}. It shows that for the line $x = z = 5$, and $y =  55, z = 5$, the numerical solution by UGKWP matches well with the reference solution by UGKS, while there exists a little difference from the analytical solution and that by GMVP. However, for the line $y = 95, z = 5$, which is far from the source region, the macroscopic flux $\psi$ is quite small and close to zero. In this case, it is quite difficult to capture the behavior of $\psi$ by the numerical method. The difference between the numerical solution and the analytical solution is larger, but the absolute error is still quite small due to the small $\psi$.

\begin{figure}[!hptb]
	\centering
	\subfloat[$\psi (x = z = 5)$]{
		\includegraphics[width = 0.3\textwidth]{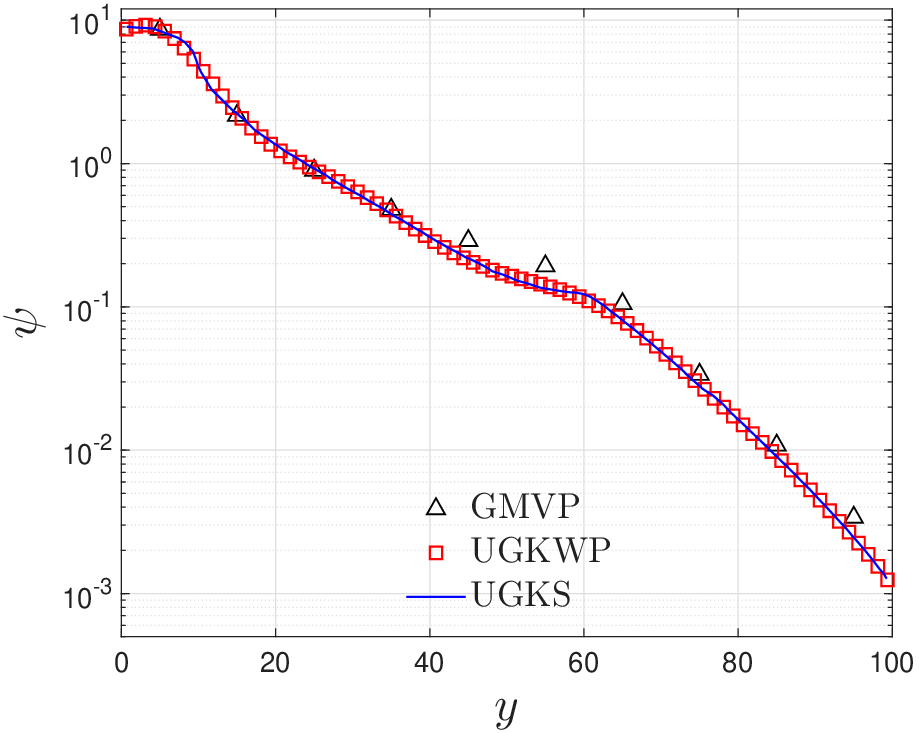}
	}\hfill
	\subfloat[$\psi (y = 55, z= 5)$]{               
		\includegraphics[width = 0.3\textwidth]{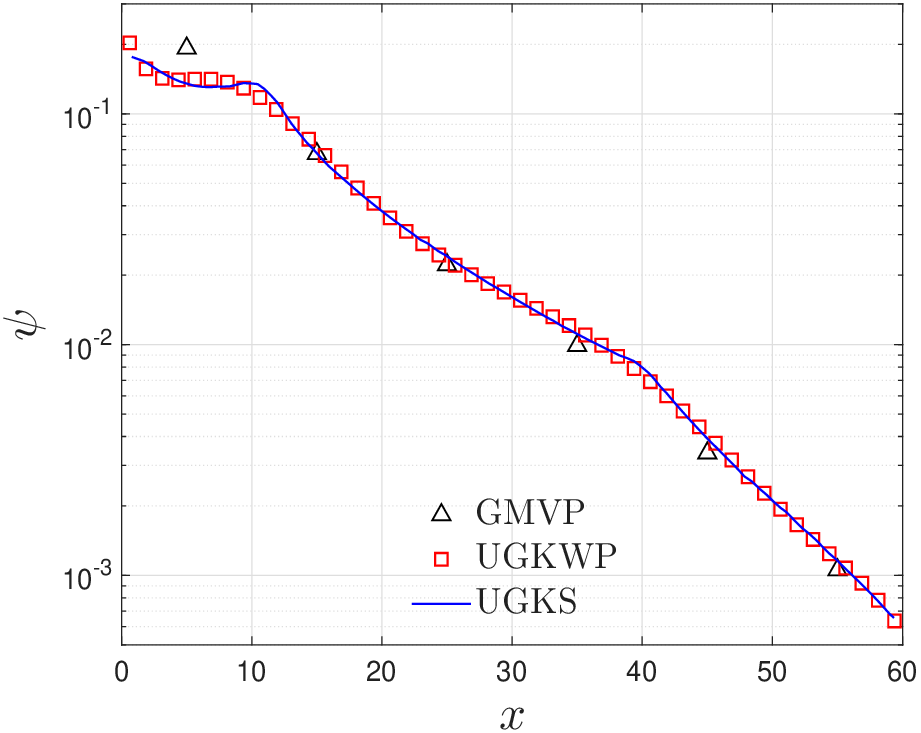}
	} \hfill
    \subfloat[$\psi (y = 95, z = 5)$]{
		\includegraphics[width = 0.3\textwidth]{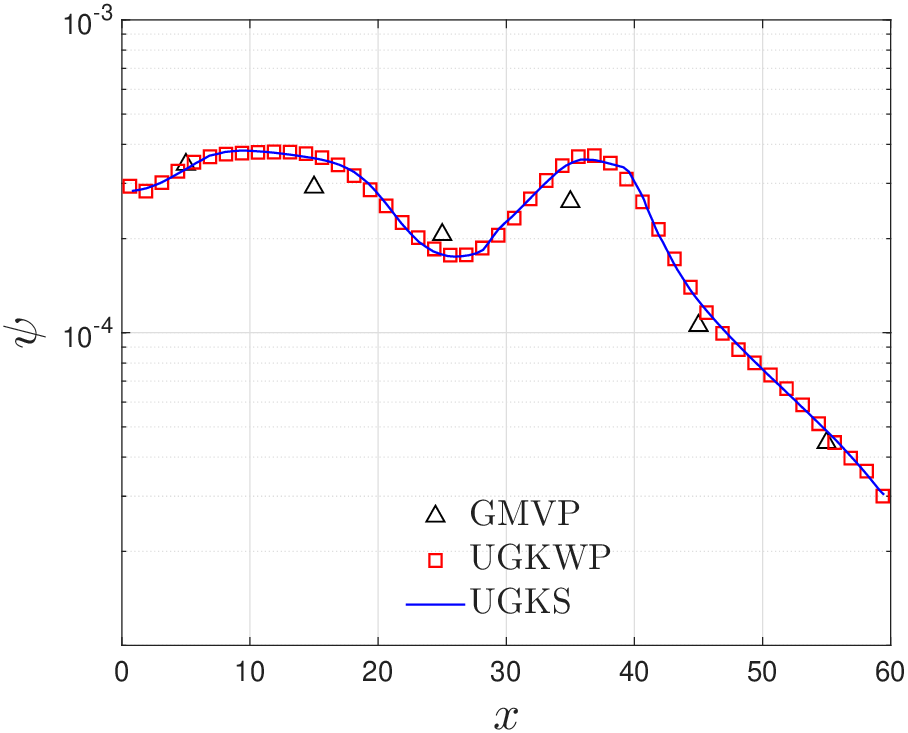}
	}
	\caption{(The top-hat problem in Sec. \ref{sec:3D_single_group_2}) The macroscopic scalar flux $\psi$ of the top-hat problem for the scattering case. Here, the triangle is the reference solution by GMVP, the blue line is the reference solution by UGKS, and the star line is the numerical solution by UGKWP.  (a) $x = z =5$; (b) $y = 55, z = 5$; (c) $y = 95, z= 5$.}  
	\label{fig:3D_tophat_sca}
\end{figure}

Then the numerical solution $\psi$ for the scattering problem at the same lines is plotted in Fig. \ref{fig:3D_tophat_sca} with the same numerical setting as the pure absorption case. Compared to the pure absorption case, the behavior of the numerical solution along $x = z = 5$ and $y = 55, z = 5$ is similar, all of which matches well with the reference solution, while $y = 95, z = 5$, the numerical solution fits much better with the reference solution by GMVP. This may be due to that the numerical solution of $\psi$ is much larger than the pure absorption case, the behavior of which is much easier to capture. 

\subsection{3D multi-group problems}
\label{sec:3D_multi_group}
In this section, the 3D neuron transport problems for the multi-group model are studied. The first problem is the two-group KUCA core problem, and the second is an axially heterogeneous Fast Breeder Reactor (FBR) problem for the four-group model. These two examples are formulated as a $k$-eigenvalue problem \citep[Sec. 1.5]{lewis1984computational}, and the objective is to compute the fundamental neutron multiplication factor $k_{\text{eff}}$ and the corresponding steady-state flux distribution as 
\begin{equation}
    \label{eq:nte_time_independent}
        \mathcal{L}\phi = \displaystyle{\frac{1}{k_\text{eff}}}\mathcal{F}\phi,
\end{equation}
with 
\begin{equation}
    \label{eq:operator_L_F}
        \begin{aligned}
            \mathcal{L}\phi &= \Omega \cdot \nabla \phi(\boldsymbol{x}, \Omega, E)+\Sigma(\boldsymbol{x}, E) \phi(\boldsymbol{x}, \Omega, E)  \\
            & \qquad \qquad -\displaystyle{\frac{1}{4 \pi}} \int_{4\pi}\int_0^{+\infty} \Sigma_{s}\left(\boldsymbol{x}, \Omega' \cdot \Omega, E' \rightarrow E\right) \phi\left(\boldsymbol{x}, \Omega', E'\right) \mathrm{d} E' \mathrm{d} \Omega', \\
            \mathcal{F}\phi &=\displaystyle{\frac{\chi(E)}{4 \pi}} \int_0^{+\infty} \nu \Sigma_{f}\left(\boldsymbol{x}, E'\right) \int_{4\pi} \phi\left(\boldsymbol{x}, \Omega', E'\right) \mathrm{d} \Omega' \mathrm{d} E'. 
        \end{aligned}
  \end{equation}
To show the fundamental neutron multiplication factor $k_{\rm eff}$ and the angular flux $\psi_g(\bmx, \Omega)$ more clearly, the normalization all over the spatial space is applied as 
 \begin{equation}
	\label{eq:normalized_flux}
	\sum_{g=1}^{G}\int_{4\pi}\int_{V} \mathrm{d}\bm{x}\mathrm{d}\Omega \nu \Sigma_{f,g}(\bm{x})\phi_{g}(\bm{x}, \Omega) = 1,
\end{equation}
 where $V$ is the total volume of the spatial space, $G$ is the total number of energy groups, and $\nu\Sigma_{f,g}(\bm{x})$ is the neutron production cross-section for group $g$ at position $\bm{x}$, detailed parameters are shown in App. \ref{APP:3D}.

The main computational challenges of these 3D multi-group problems lie in accurately modeling the strong material heterogeneity between fuel, control, and reflector regions, resolving the multi-group energy structure, and ensuring stable and efficient convergence of the eigenvalue iteration, especially under the high-absorption scattering effects introduced by control rods. 

\subsubsection{Two-group KUCA core problem}
\label{sec:kuca}
\begin{figure}[!hptb]
	\centering
	\subfloat[3D view]{
		\includegraphics[width = 0.24\textwidth]{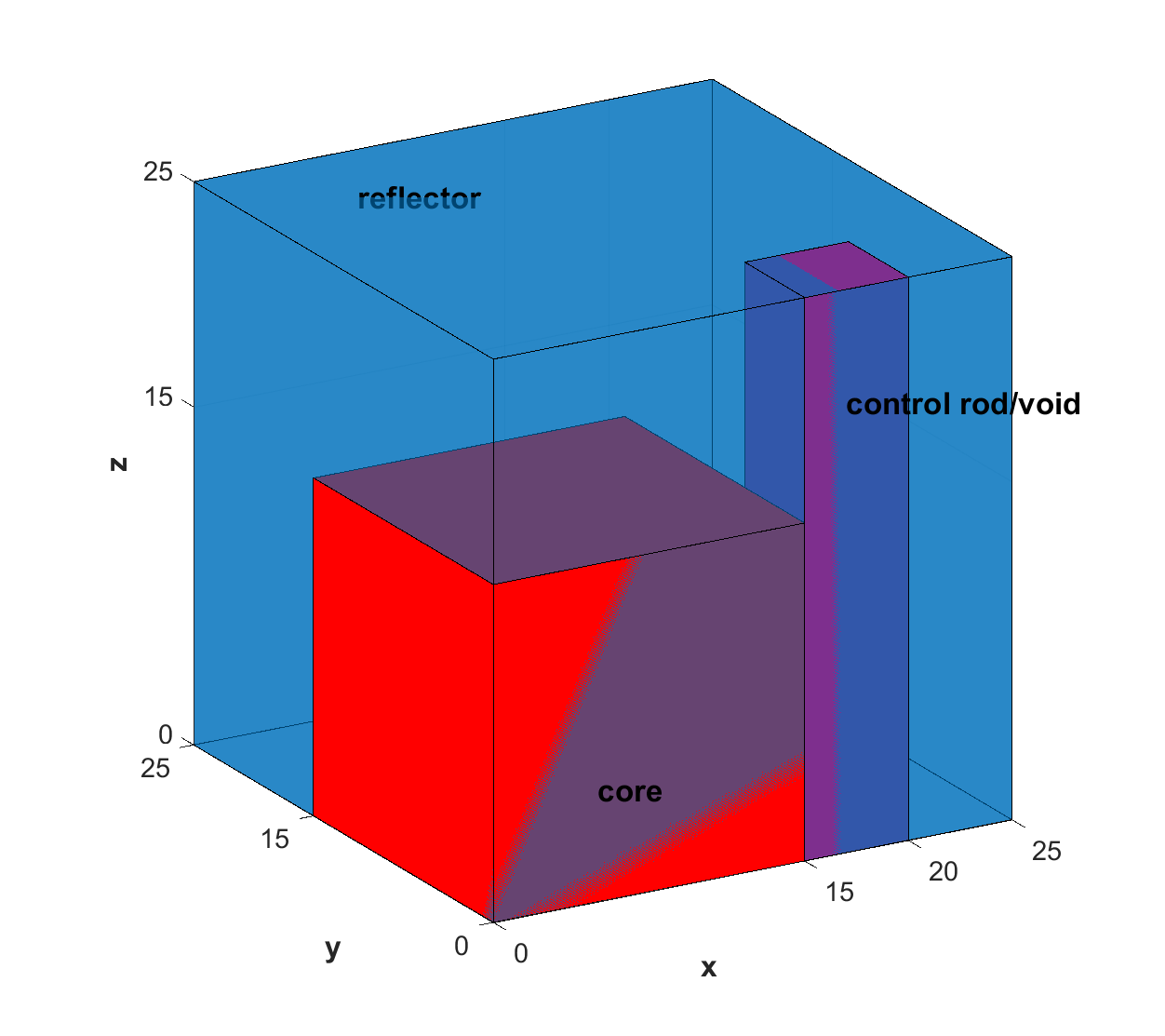}
	}
	\subfloat[boundary condition: XOY]{               
		\includegraphics[width = 0.24\textwidth]{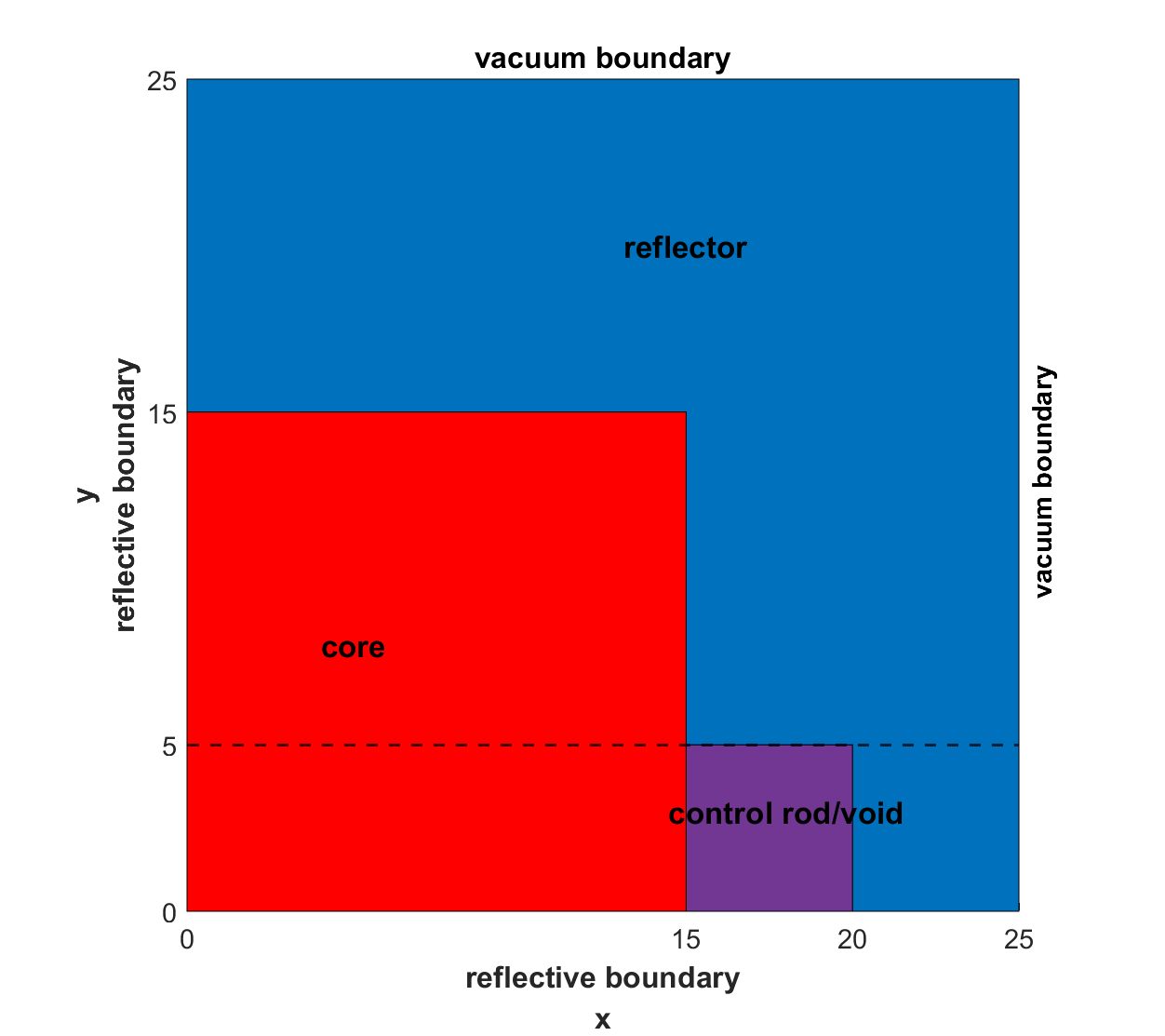}
	}
	\subfloat[boundary condition: XOZ]{               
		\includegraphics[width = 0.24\textwidth]{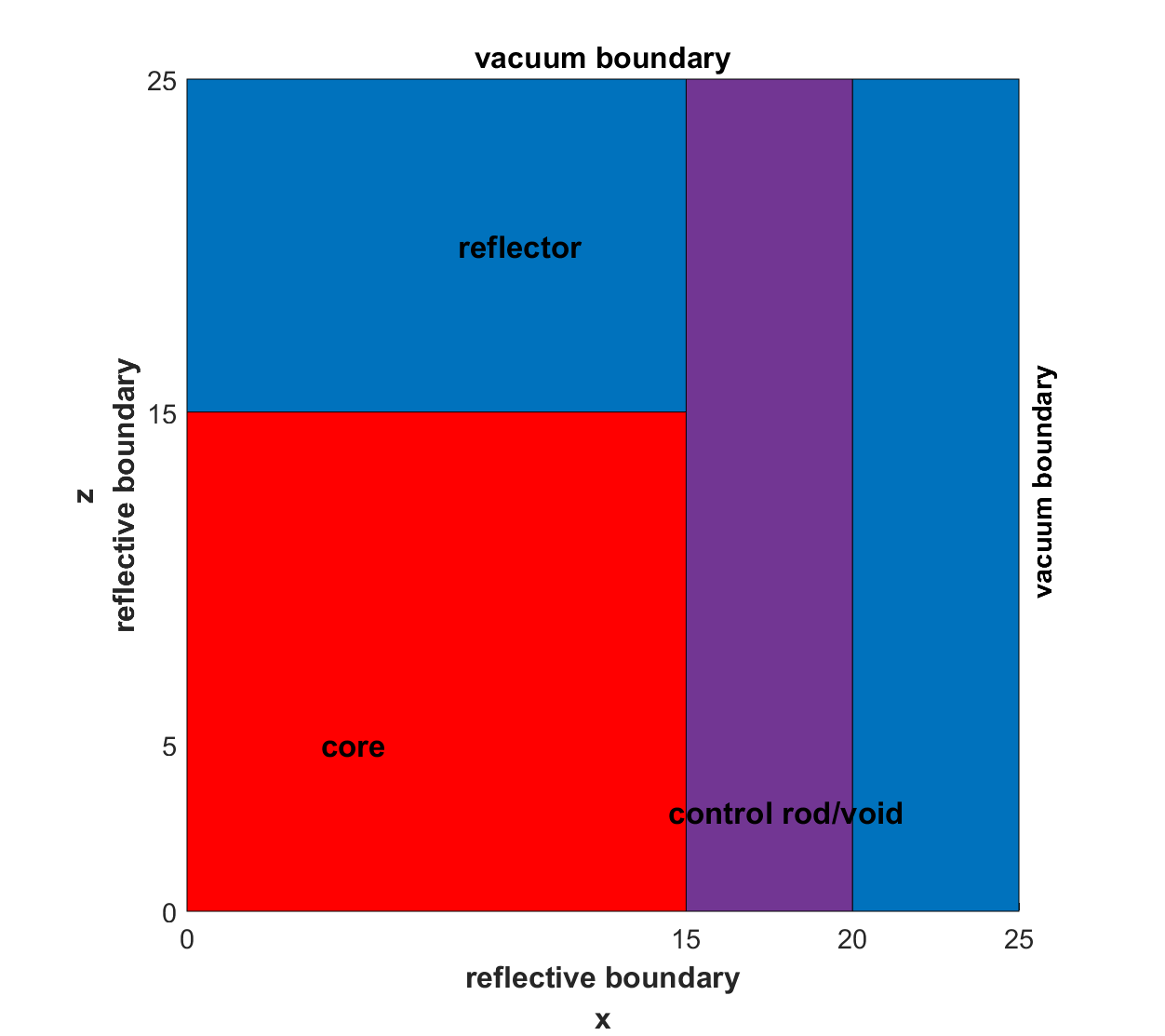}
	}
	\subfloat[boundary condition: YOZ]{               
		\includegraphics[width = 0.24\textwidth]{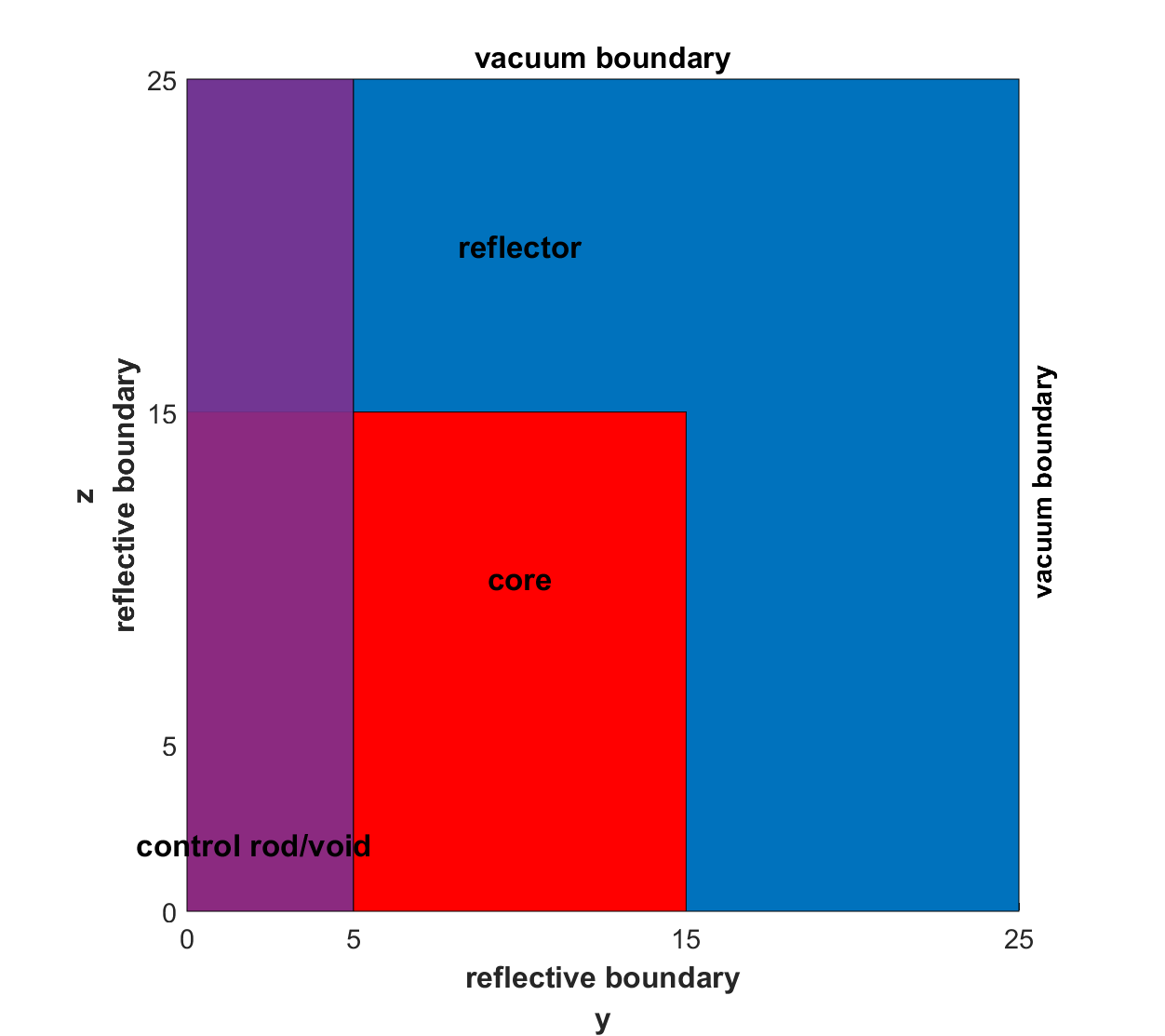}
	}    
	\caption{(The KUCA core two-group benchmark problem in Sec. \ref{sec:kuca}) 
		The computational domain and boundary settings.}  
	\label{fig:kuca_domain}
\end{figure}

In this section, the simplified model of a Light Water Reactor \cite{bunker1983early, duran2022verification}, based on the Kyoto University Critical Assembly (KUCA) \cite{takeda19913}, is studied to validate the UGKWP method for the multi-group neutron transport simulation. Due to the symmetry of the reactor, one-eighth of the total geometry is shown in Fig. \ref{fig:kuca_domain}, which consists of a cube reactor with an edge length of $25$. The core fuel region, shaped as a cube with an edge length of $15$, is at the center, which contains the fissile material and serves as the main site of the neutron-induced fission. A control rod region, shaped as a $5\times 5 \times 25$ rectangular prism, is adjacent to one side of the core fuel region. It is a movable piece of neutron-absorbing material utilized to control the reactor. The rest regions are the reflector, forming a layer along all the boundaries. It reflects the escape of neutrons back into the core, thereby improving the neutron economy. More details on the components of a nuclear reactor can be found in \cite{lamarsh1975introduction}, and the detailed cross section parameters and the energy group structure are listed in Tab. \ref{tab:KUCA_cross} and \ref{tab:KUCA_energy}, respectively.

The specification of the problem provides two different configurations. The first case is called rod-out, where the control rod is completely removed, leaving an empty region. The second case is called rod-in, where the control rod is completely inserted. Moreover, to quantify the impact of control rod insertion, the control rod (CR) worth is calculated \cite{abramowitz1948handbook}. This critical safety parameter is defined as the total change in reactivity upon rod insertion. First, define the reactivity $\rho$ in terms of the effective multiplication factor $k_\text{eff}$ as
\begin{equation}
	\rho_{s} = \frac{k_\text{eff, s} - 1}{k_\text{eff, s}}, \qquad \text{s} = \text{out, in},
\end{equation}
where $k_\text{eff, out}$ and $k_\text{eff, in}$ are the effective multiplication factors for the control rods removed and inserted, respectively. Then, the CR worth, denoted as $\Delta\rho$, is the difference between the reactivity of these two states
\begin{equation}
    \label{eq:cr_worth}
	\Delta\rho = \left( \frac{k_\text{eff, out} - 1}{k_\text{eff, out}} \right) - \left( \frac{k_\text{eff, in} - 1}{k_\text{eff, in}} \right).
\end{equation}

\begin{figure}[!hptb]
	\centering
	\subfloat[rod-out, $\psi_1$]{
		\includegraphics[width = 0.45\textwidth, trim=10 0 0 2, clip]{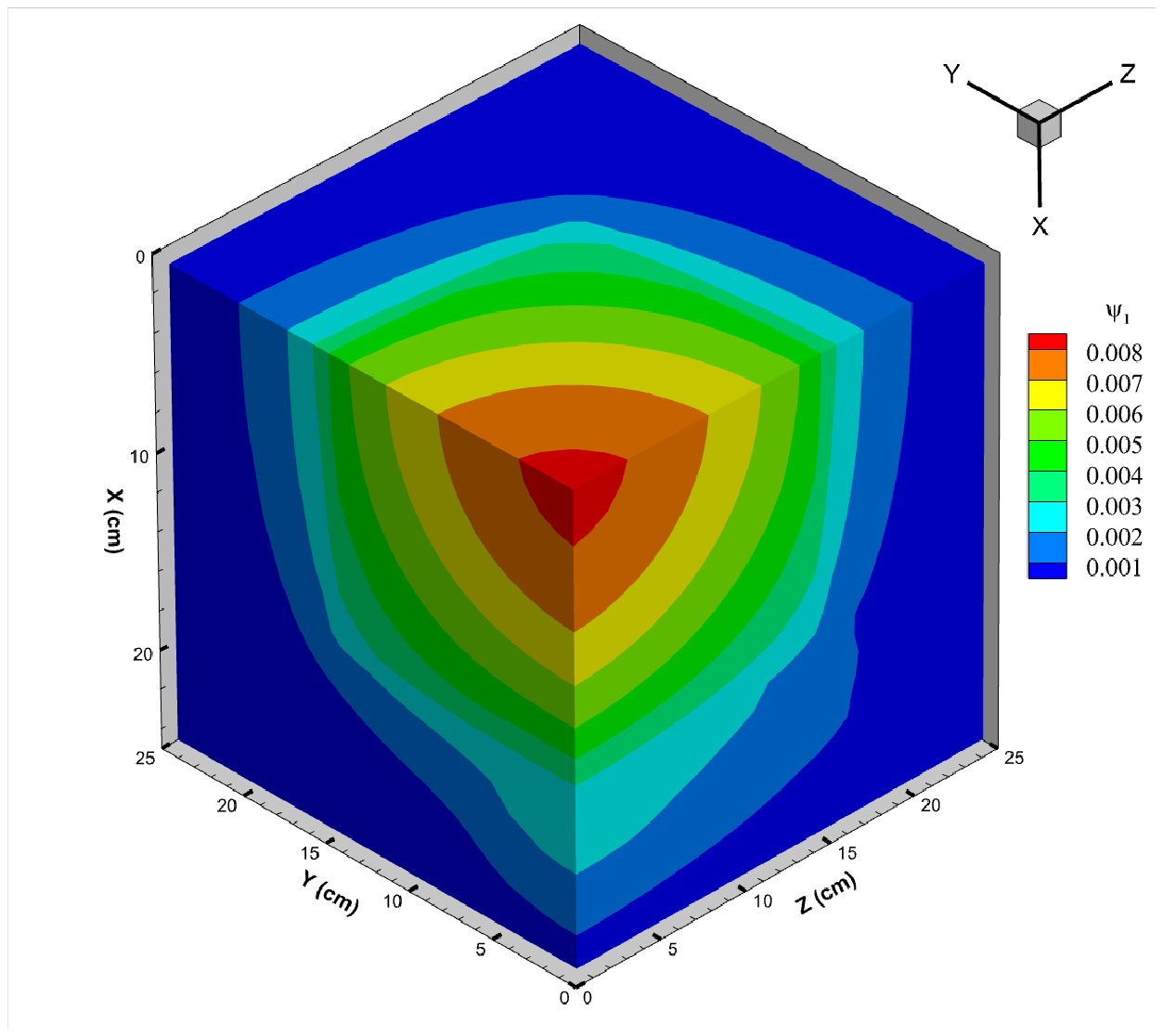}
            \label{fig:kuca_po_1}
	}\hfill
	\subfloat[relative error of $\psi_1$ with UGKS]{               
		\includegraphics[width = 0.45\textwidth, trim=10 0 0 2, clip]{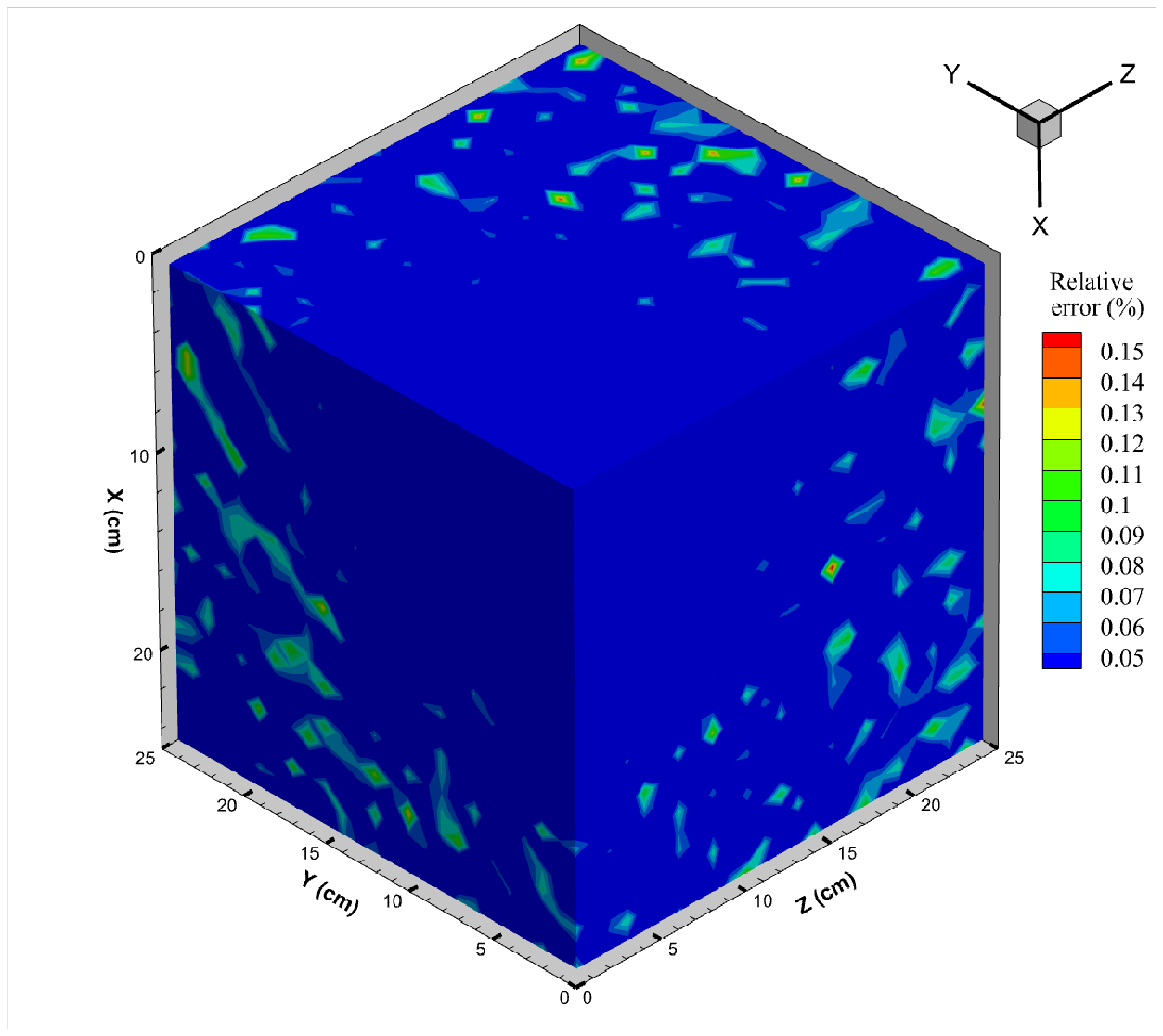}
            \label{fig:err_kuca_po_1}
	}\\
	\subfloat[rod-out, $\psi_2$]{
		\includegraphics[width = 0.45\textwidth, trim=10 0 0 2, clip]{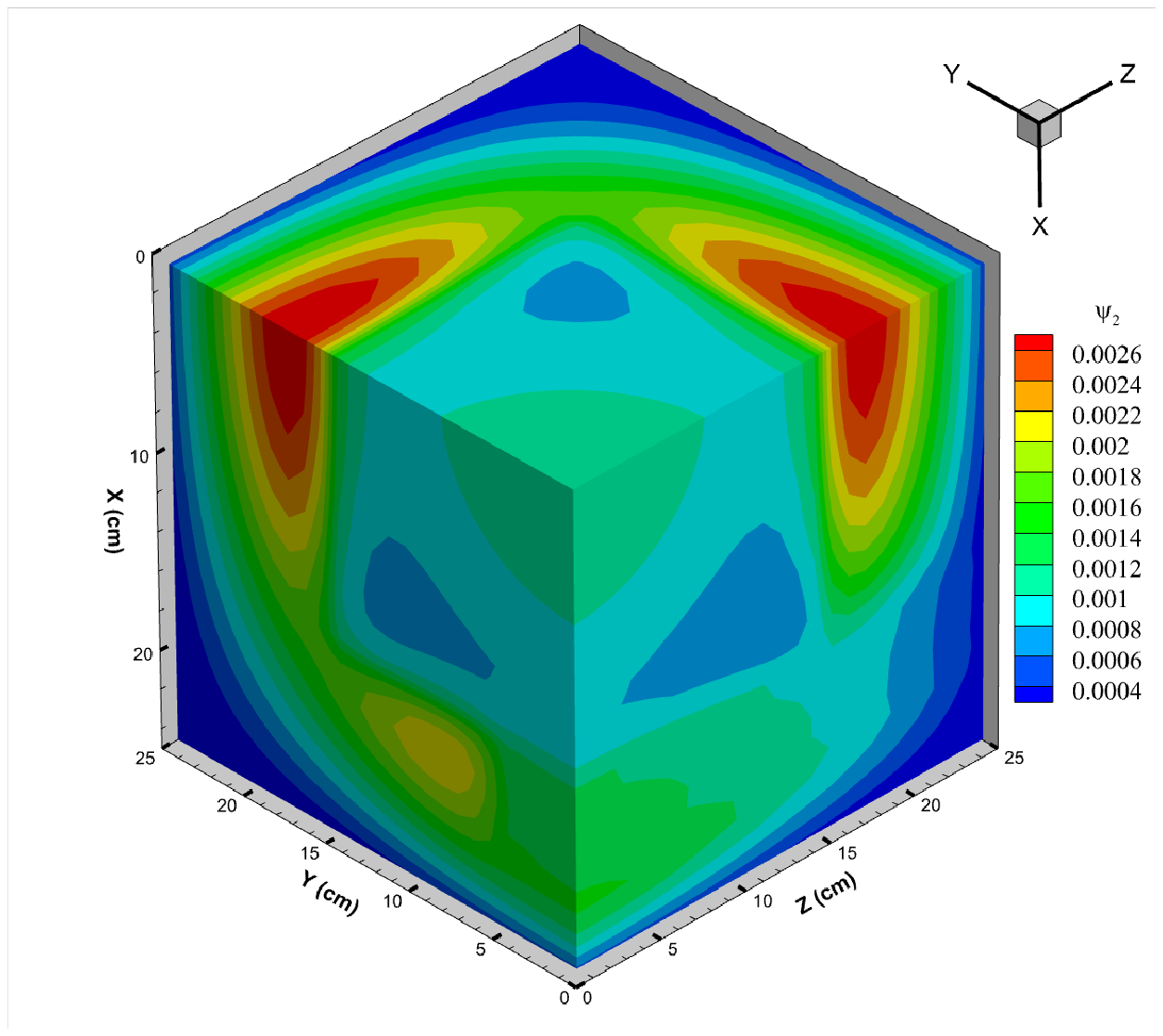}
            \label{fig:kuca_po_2}
	}\hfill
	\subfloat[relative error of $\psi_2$ with UGKS]{
		\includegraphics[width = 0.45\textwidth, trim=10 0 0 2, clip]{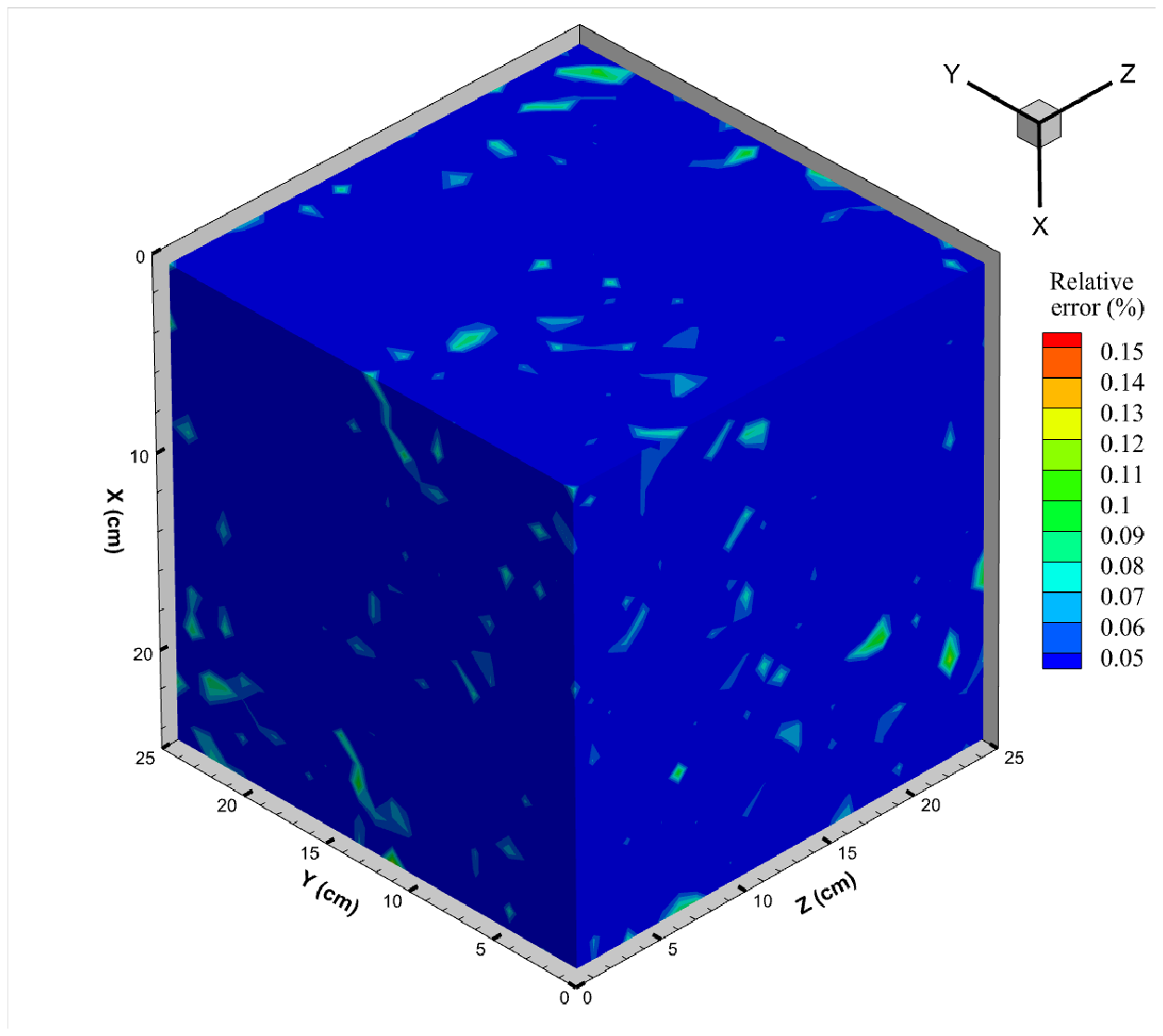}
            \label{fig:err_kuca_po_2}
	}
	\caption{(The KUCA core two-group benchmark problem in Sec. \ref{sec:kuca}) The numerical solution of $\psi_1$ and $\psi_2$ for rod-out, and the relative error with UGKS. Here, the left column is the numerical solution of $\psi_1$ and $\psi_2$, and the right column is the relative error with UGKS. }
	\label{fig:kuca_rod_out}
\end{figure}

In the simulation, the mesh size is set as $N_x = N_y = N_z = 25$ and the CFL number is $\mathrm{CFL} = 0.8$. The parameters of the detailed group-wise cross-section are provided in Appendix \ref{APP:kuca}. The multiplication factors $k_\text{eff}$ and the CR worth are presented in Tab. \ref{tab:3D_kuca_k_eff}, where the related results obtained by Monte Carlo (MC) \cite{takeda19913}, UGKS \cite{shuang2019parallel}, and MOC \cite{liu2011new} are listed. Here, the result by MC is taken as the reference solution. The multiplication factors $k_\text{eff}$ obtained by UGKWP are more consistent with Monte Carlo compared with UGKS and MOC. Moreover, the CR worth $\Delta \rho$ obtained by these three methods is all very close to that of Monte Carlo, and the error is less than $1\%$.

\begin{table}[!htbp]
	\centering
    \def\arraystretch{1.5}
    {\footnotesize
    \begin{tabular}{c|c c|c c |c}
		% \hline Code & Rod out &   &  Rod in &   \\
            code  & rod out  $k_\text{eff}$  & error & rod in $k_\text{eff}$ &error   &CR worth $\Delta \rho$ \\
		\hline
		Monte Carlo & $0.9778 \pm 0.0005$ & - & $0.9624\pm0.0005$ &  - & $(1.64\pm0.07)$E-02\\
     	UGKWP & 0.9774 & 5.11E-04 & 0.9622 &  2.08E-04   & 1.616E-02 \\
           UGKS & 0.9772 & 6.13E-04 & 0.9620 &  4.16E-04   & 1.617E-02\\
		MOC & 0.9764 & 1.43E-03 & 0.9609 &  1.56E-03  & 1.652E-02
		\end{tabular}
	\caption{(The KUCA core two-group benchmark problem in Sec. \ref{sec:kuca}) The $k_\text{eff}$ and the CR worth $\Delta \rho$ obtained by UGKWP for rod-out and rod-in, and the comparison with Monte Carlo, UGKS, and MOC.} 
	\label{tab:3D_kuca_k_eff}
    }
\end{table}

\begin{figure}[!hptb]
	\centering
	\subfloat[$\psi_{1}$ (rod-out, $x = y = z$)]{
		\includegraphics[width = 0.23\textwidth]{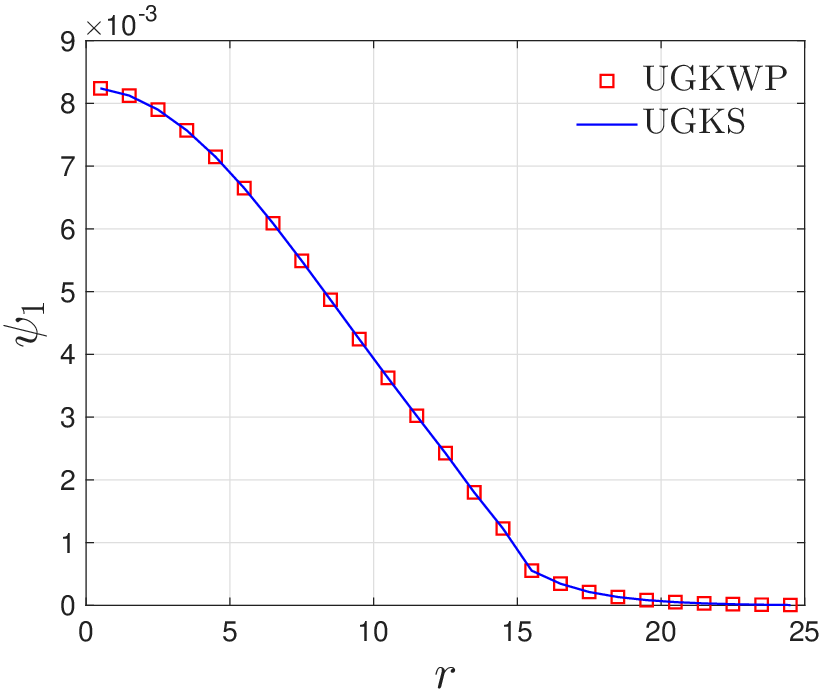}
	}\hfill
	\subfloat[$\psi_{1}$ (rod-out, $x = 16, y = 5$)]{               
		\includegraphics[width = 0.23\textwidth]{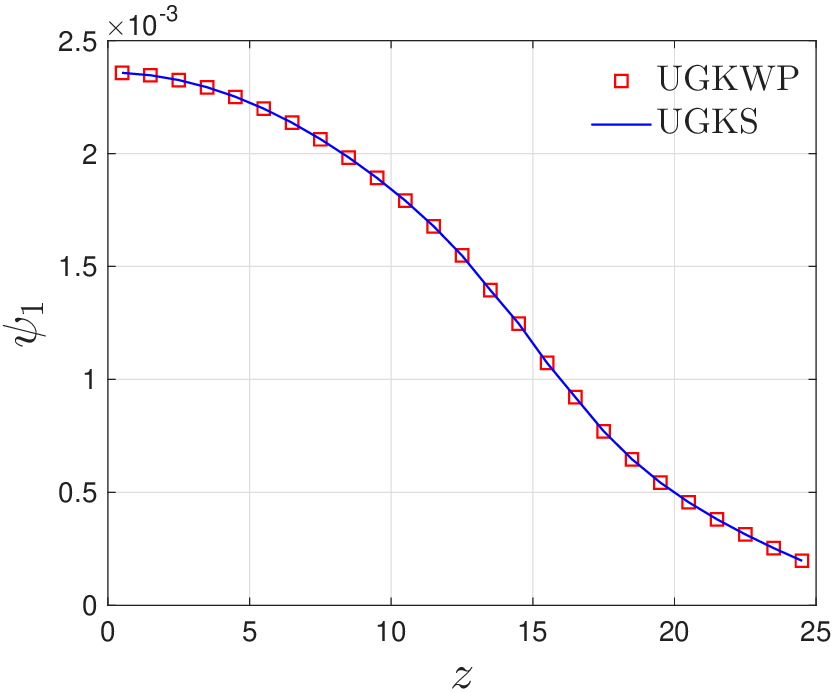}
	}\hfill
	\subfloat[$\psi_{2}$ (rod-out, $x = y = z$)]{
		\includegraphics[width = 0.23\textwidth]{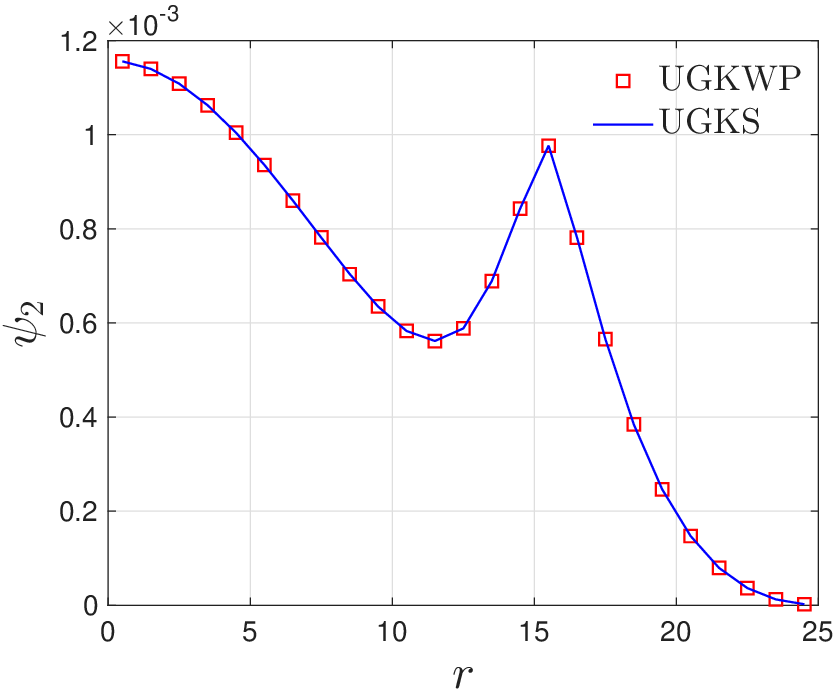}
	}\hfill
	\subfloat[$\psi_{2}$ (rod-out, $x = 16, y = 5$)]{               
		\includegraphics[width = 0.23\textwidth]{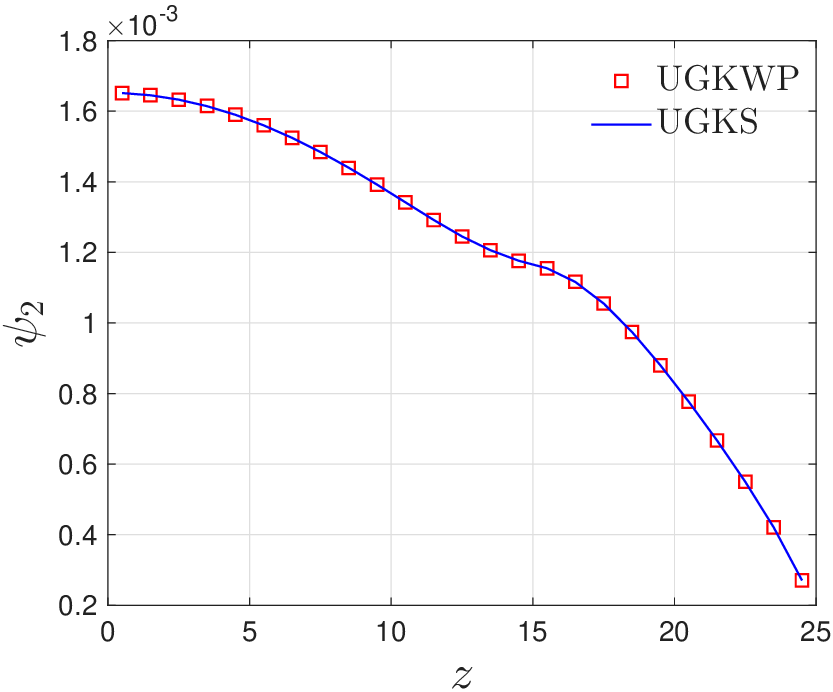}
	}
	\caption{(The KUCA core two-group benchmark problem in Sec. \ref{sec:kuca}) 
    The numerical solution of $\psi_1$ and $\psi_2$ for rod-out along different lines.     (a) $\psi_1$ along $x = y = z$; (b) $\psi_1$ along $x = 16, y = 5$; (c) $\psi_2$ along $x = y = z$; (d) $\psi_2$ along $x = 16, y = 5$.}
	\label{fig:kuca_po_along_line}
\end{figure}

The numerical solution of the macroscopic scalar flux $\psi_1$ and $\psi_2$ for Group 1 and Group 2 for rod-out is shown in Fig. \ref{fig:kuca_po_1} and \ref{fig:kuca_po_2}, respectively, where the behavior of $\psi$ for the two groups is quite different. Moreover, a point-wise relative error of the numerical solution between UGKWP and UGKS \cite{shuang2019parallel} is shown in Fig. \ref{fig:err_kuca_po_1} and \ref{fig:err_kuca_po_2}. For both macroscopic scalar flux, this error is at the order of $10^{-3}$, indicating the accuracy of UGKWP. Moreover, the numerical solution along $x = y = z$ and $x = 16, y = 5$ for the two groups is plotted in Fig. \ref{fig:kuca_po_along_line}. Along the main diagonal, the neutron traverses various material regions, including fuel, clad, and moderator. Even though, the numerical solution by UGKWP and UGKS is visually indistinguishable, demonstrating an excellent agreement of UGKPW and UGKS in capturing the sharp flux gradients at material interfaces and the overall flux shape. A similar level of agreement is observed along the transverse line, which highlights the flux variation across different fuel assemblies. These line-by-line comparisons effectively confirm the high fidelity of the UGKWP solution.

\begin{figure}[!hptb]
	\centering
	\subfloat[rod-in, $\psi_1$]{
		\includegraphics[width = 0.45\textwidth, trim=10 0 0 2, clip]{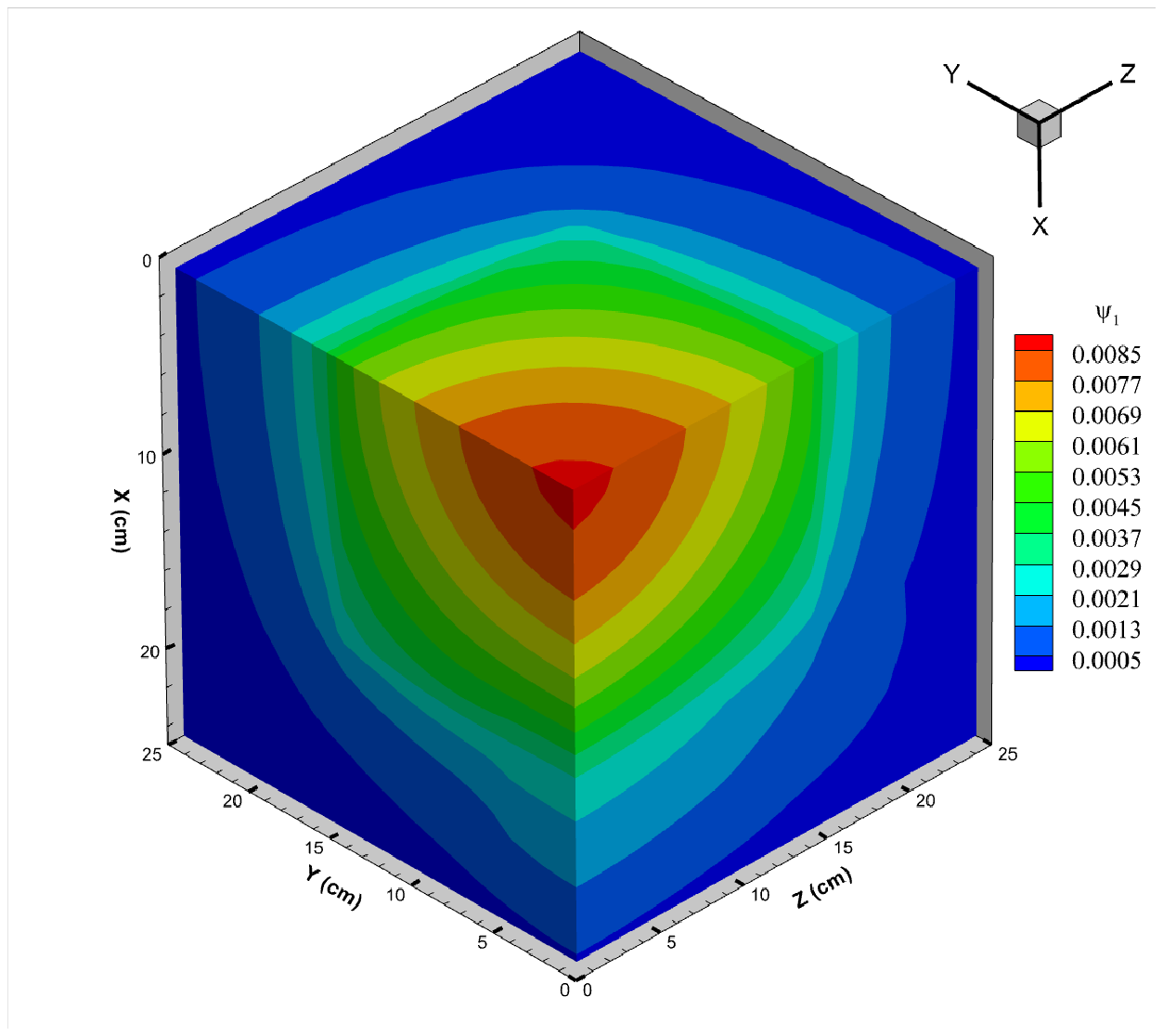}
            \label{fig:kuca_pi_1}
	}\hfill
	\subfloat[relative error of $\psi_1$ with UGKS]{               
		\includegraphics[width = 0.45\textwidth, trim=10 0 0 2, clip]{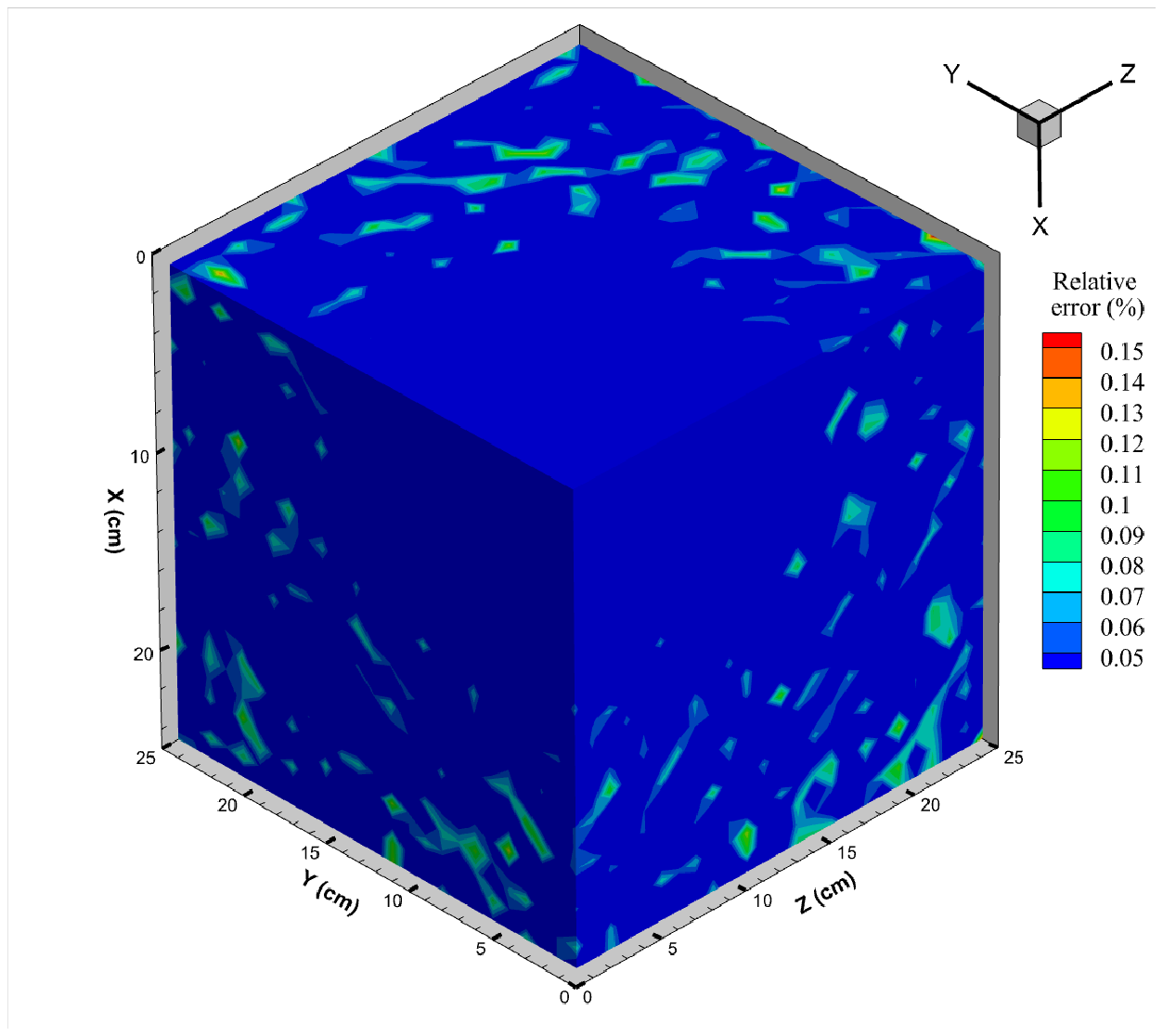}
            \label{fig:err_kuca_pi_1}
	}\\
        \subfloat[rod-in, $\psi_2$]{
		\includegraphics[width = 0.45\textwidth, trim=10 0 0 2, clip]{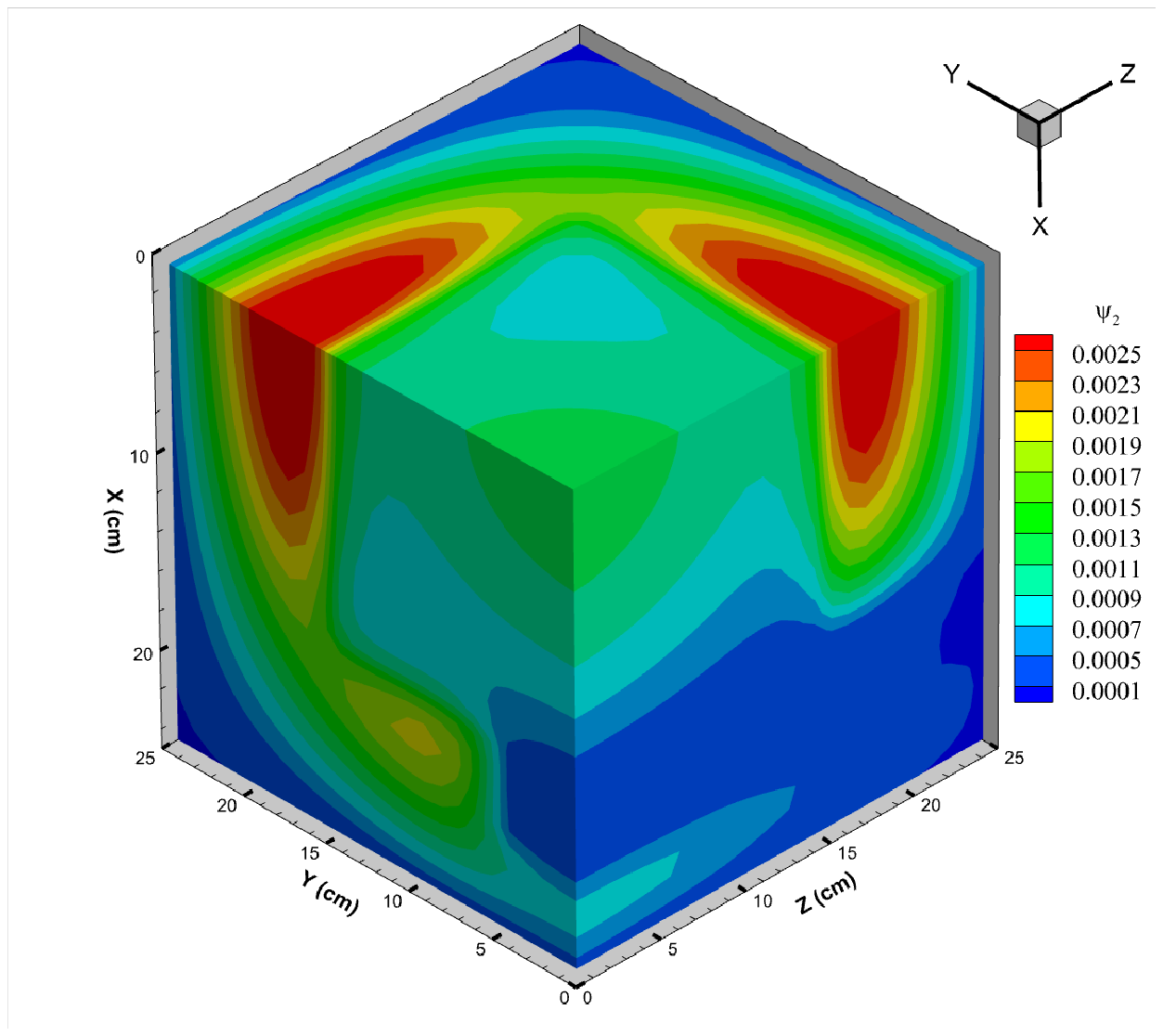}
            \label{fig:kuca_pi_2}
	}\hfill
	\subfloat[relative error of $\psi_2$ with UGKS]{               
		\includegraphics[width = 0.45\textwidth, trim=10 0 0 2, clip]{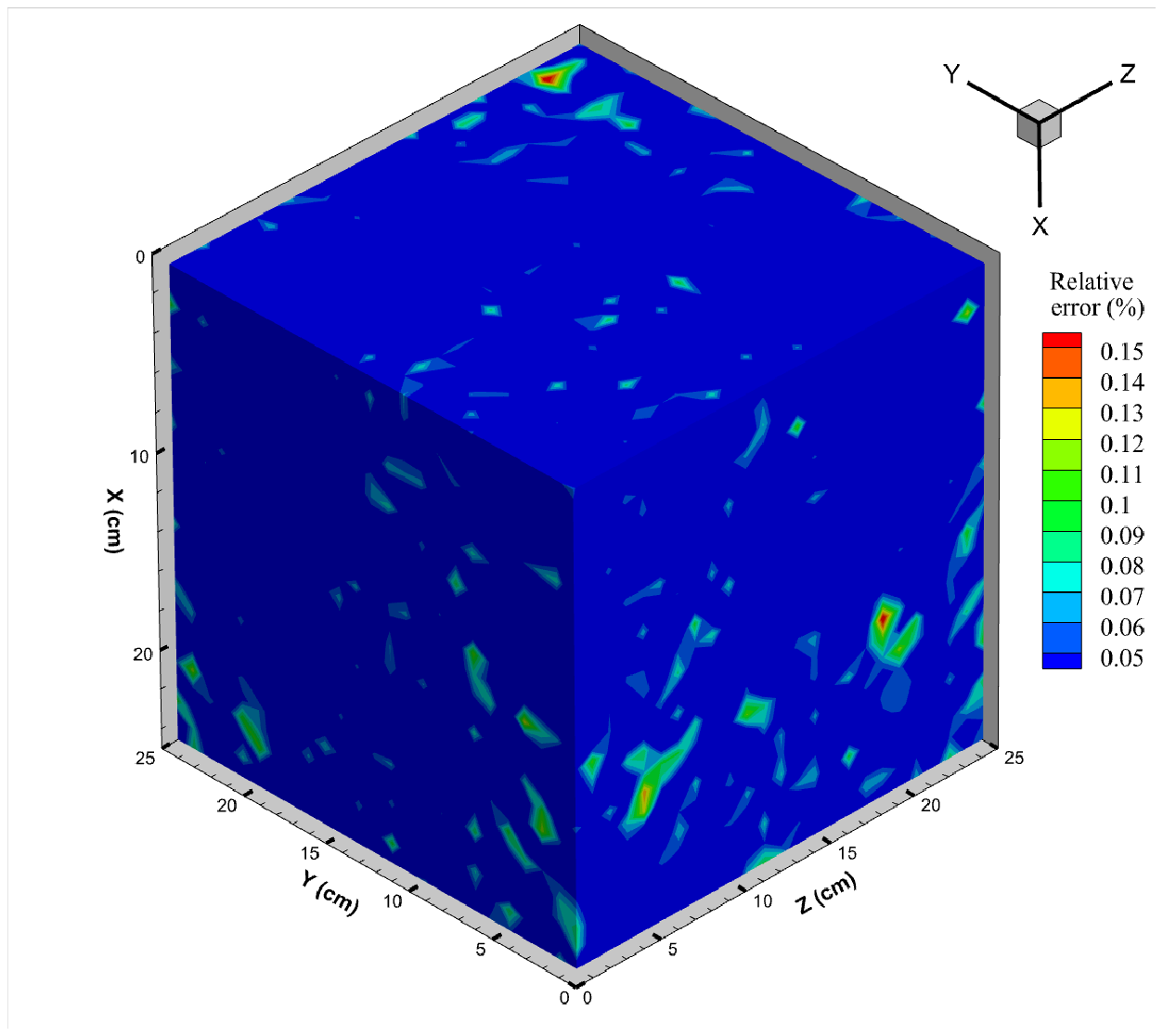}
            \label{fig:err_kuca_pi_2}
	}
	\caption{(The KUCA core two-group benchmark problem in Sec. \ref{sec:kuca}) 
    The numerical solution of $\psi_1$ and $\psi_2$ for rod-in, and the relative error with UGKS. Here, the left column is the numerical solution of $\psi_1$ and $\psi_2$, and the right column is the relative error with UGKS.   }  
\label{fig:kuca_rod_in}
\end{figure}

The numerical solution of the macroscopic scalar flux $\psi_1$ and $\psi_2$ of Group 1 and Group 2 for rod-in is plotted in Fig. \ref{fig:kuca_pi_1} and \ref{fig:kuca_pi_2}, respectively. Compared with Fig. \ref{fig:kuca_rod_out}, the influences of the control rod can be visually shown from the distribution of $\psi_1$ and $\psi_2$. With the control rod in, the values of $\psi_1$ and $\psi_2$ are relatively smaller at the same location, especially for $\psi_2$ near the control rod. A similar relative error between UGKWP and UGKS is plotted in Fig. \ref{fig:err_kuca_pi_1}, and \ref{fig:err_kuca_pi_2}, where this relative error of the macroscopic scalar flux for both groups remains consistently $0.15\%$ throughout the computational domain. The numerical solution along $x = y = z$ and $x = 16, y = 5$ also indicates the consistency of UGKWP and UGKS.

\begin{figure}[!hptb]
	\centering
	\subfloat[$\psi_{1}$ (rod-in, $x=y=z$)]{
		\includegraphics[width = 0.23\textwidth]{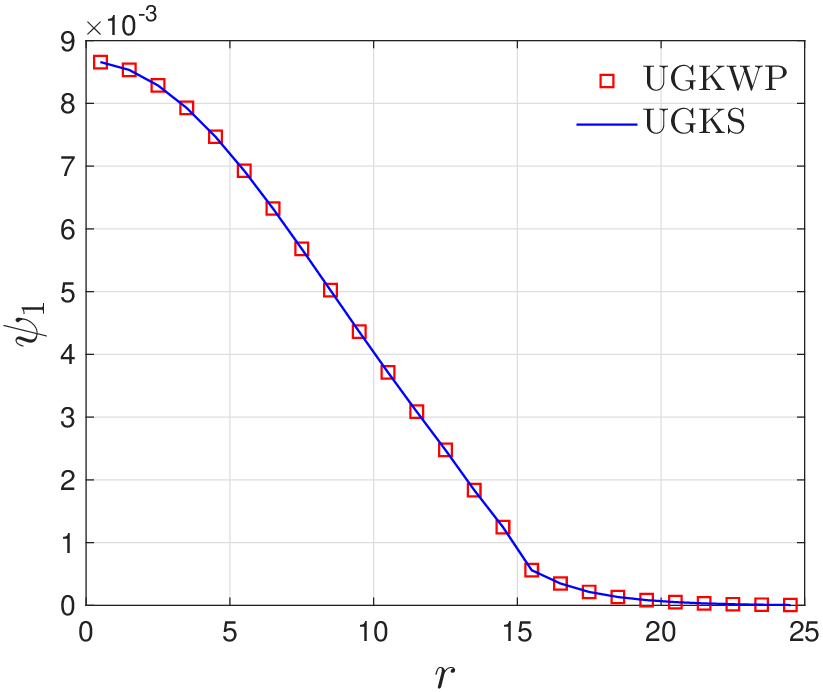}
	}\hfill
	\subfloat[$\psi_{1}$ (rod-in, $x=16, y=5$)]{               
		\includegraphics[width = 0.23\textwidth]{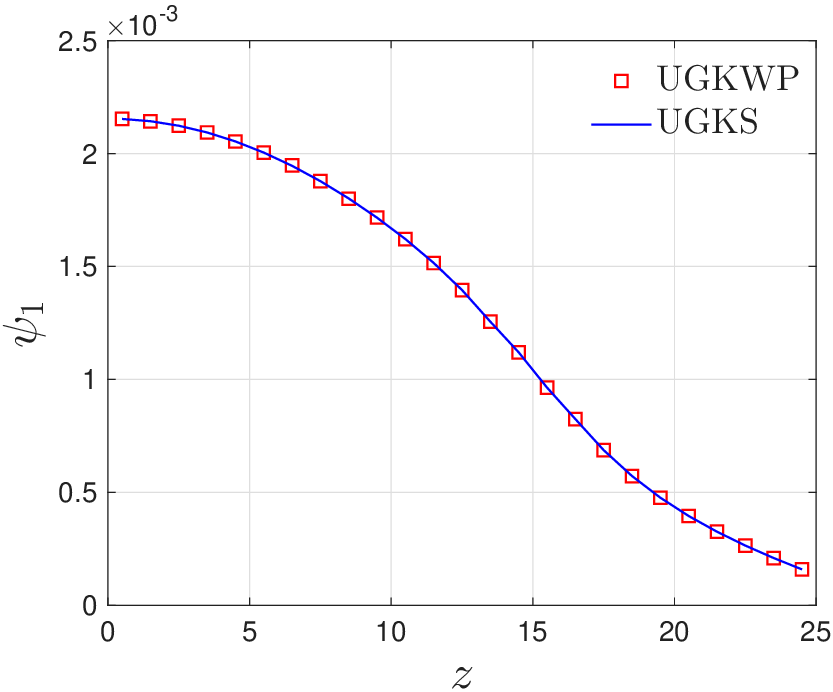}
	}\hfill
	\subfloat[$\psi_{2}$ (rod-in, $x=y=z$)]{
		\includegraphics[width = 0.23\textwidth]{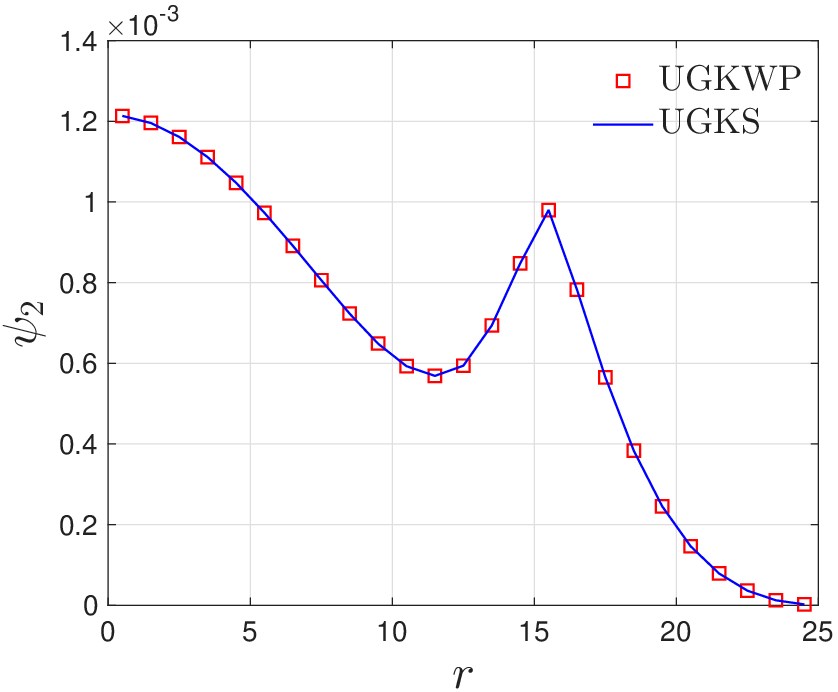}
	}\hfill
	\subfloat[$\psi_{2}$ (rod-in, $x=16, y=5$)]{               
		\includegraphics[width = 0.23\textwidth]{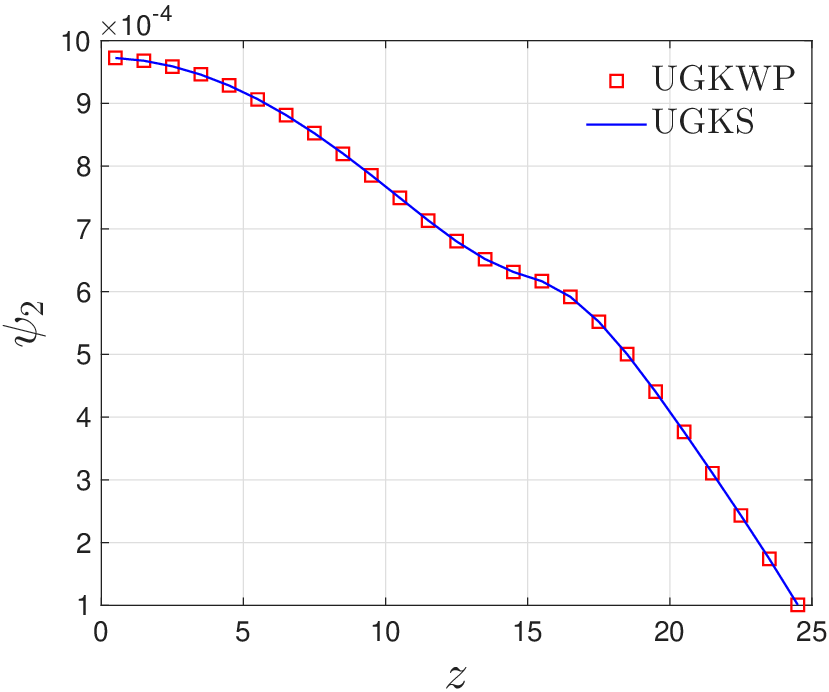}
	}
	\caption{(The KUCA core two-group benchmark problem in Sec. \ref{sec:kuca}) 
    (The KUCA core two-group benchmark problem in Sec. \ref{sec:kuca}) 
    The numerical solution of $\psi_1$ and $\psi_2$ for rod-in along different lines. (a) $\psi_1$ along $x = y = z$; (b) $\psi_1$ along $x = 16, y = 5$; (c) $\psi_2$ along $x = y = z$; (d) $\psi_2$ along $x = 16, y = 5$.    
    }  
	\label{fig:kuca_pi_along_line}
\end{figure}
\subsubsection{Four-group axially heterogeneous FBR}
\label{sec:FBR}
\begin{figure}[!hptb]
	\centering
	\subfloat[$xOy$]{
		\includegraphics[width = 0.45\textwidth]{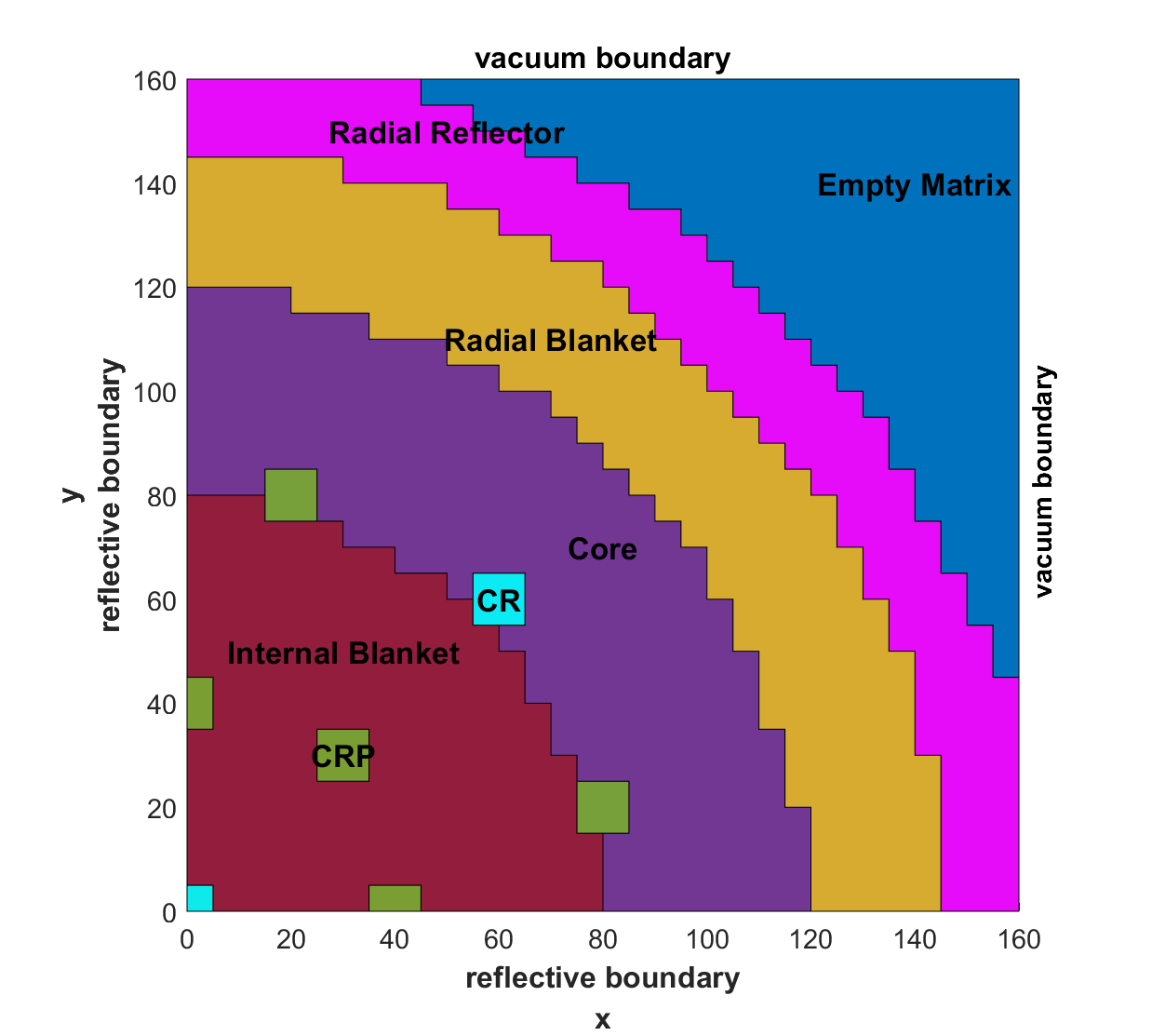}
	}
	\subfloat[3D view]{
		\includegraphics[width = 0.45\textwidth]{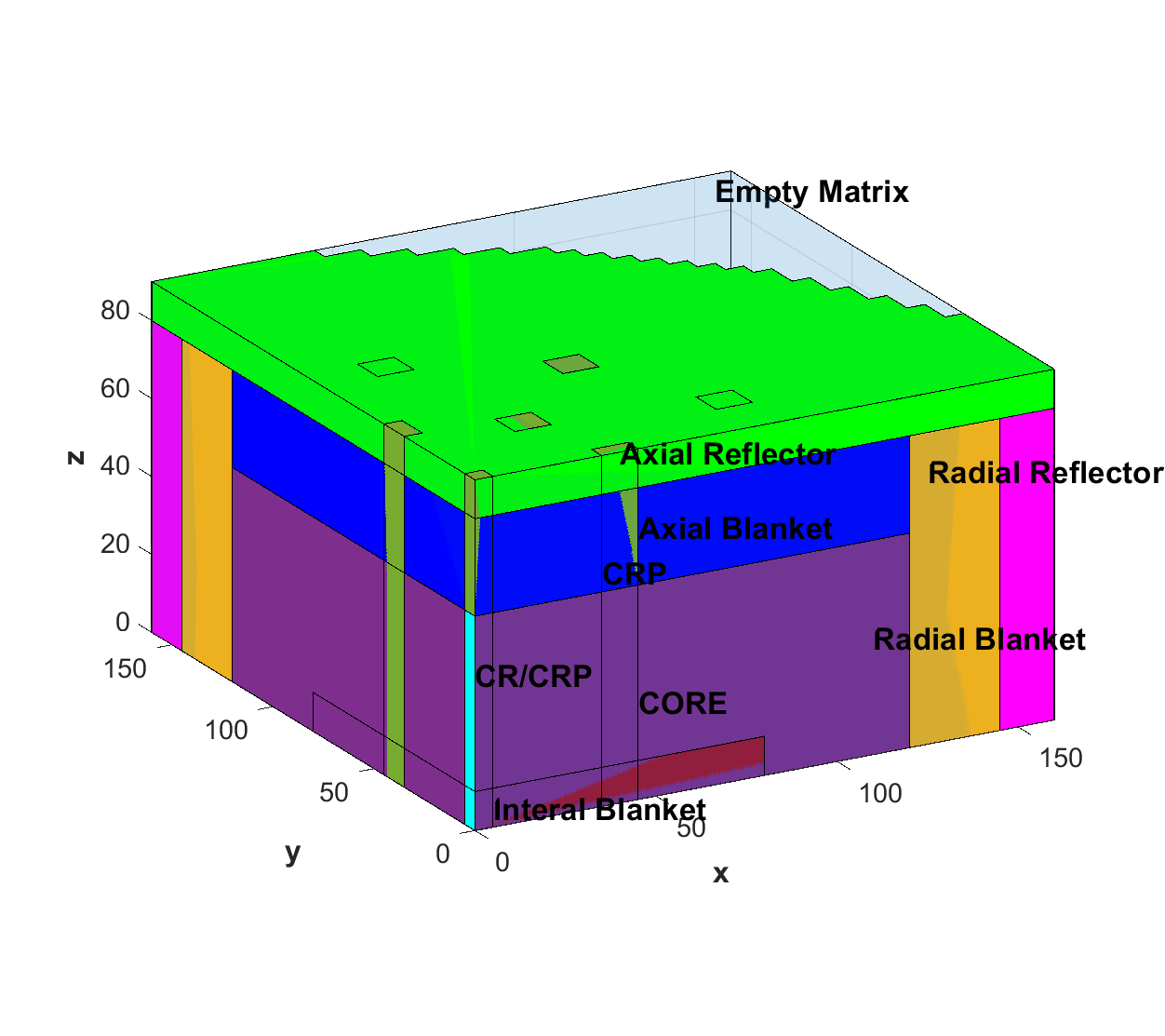}
	}
	\caption{(Four-group axially heterogeneous FBR problem in Sec.~\ref{sec:FBR}) The computational domain and boundary settings in FBR problem.}  
	\label{fig:FBR_domain}
\end{figure}

In this example, a four-group neutron transport simulation is carried out for an axially heterogeneous fast breeder reactor (FBR) core \citep[Ch.~1]{waltar1981fast}, following the benchmark in \cite{takeda19913}. The geometry is shown in Fig. \ref{fig:FBR_domain} and represents a one-eighth symmetric section of the reactor, with overall dimensions of $160~\times 160~ \times 90$. The core consists of multiple regions with distinct materials, including a fissile core region, axial and radial blankets, axial and radial reflectors, and control rod regions, resulting in a highly heterogeneous configuration.

The fissile core region is the main source of fission neutrons. It is surrounded by blanket regions (axial, radial, and internal), which consist of fertile material used for breeding new fissile isotopes \cite{lamarsh1975introduction}. The reflector regions along both axial and radial boundaries help return escaping neutrons back to the core. Control rod regions, inserted at the center and off-center, are used to absorb excess neutrons and control reactivity. In the absence of control rods, these zones are replaced with sodium, which serves as a low-absorption moderator. 
\begin{table}[!htbp]
	\centering
	\def\arraystretch{1.5}
	{\footnotesize
		\begin{tabular}{c|c c|c c|c}
			code & rod out $k_\text{eff}$ & error &rod in $k_\text{eff}$ & error & CR worth \\
			\hline
			Monte-Carlo & $0.9709 \pm 0.0002$ & - & $1.0005 \pm 0.0002$ & - & $(3.05 \pm 0.03)$E-02 \\
			IUGKS       & 0.97077 & 1.34E-04 & 1.00072 & 2.20E-04 & 3.083E-02 \\
			UGKWP       & 0.97080 & 1.03E-04 & 1.00068 & 1.80E-04 & 3.075E-02 \\
		\end{tabular}
		\caption{(Four-group axially heterogeneous FBR problem in Sec. \ref{sec:FBR}) The $k_\text{eff}$ and the CR worth $\Delta \rho$ obtained by UGKWP for rod-out and rod-in, and the comparison with Monte Carlo, and IUGKS.}
		\label{tab:3D_FBR_k_eff}
	}
\end{table}

The problem is also a $k$-eigenvalue problem. Compared to the KUCA benchmark, which involves thermal neutron groups, this problem features a fast neutron spectrum characteristic of fast breeder reactors \cite{waltar1981fast}. Group-wise cross-sections, core composition data, and the energy group structure are provided in Tab. \ref{tab:FBR_cross_1}, \ref{tab:FBR_cross_2}, and \ref{tab:FBR_energy}, respectively. Similarly, the rod-out and rod-in cases are considered, with the multiplication factor $k_\text{eff}$ and the CR worth $\Delta \rho$ studied. 

% In this study, two cases are considered. In the first case, the control rods are fully inserted as in Fig. \ref{fig:FBR_domain}. In the second case, all control rods are removed. The main challenges include strong axial and radial heterogeneity, energy-group coupling, and reactivity control sensitivity under different absorber configurations. 

In the simulation, the mesh size is set as $N_x = N_y = 64, N_z = 36$ and the CFL number is $\mathrm{CFL} = 0.4$. The numerical results of $k_\text{eff}$ and CR worth are listed in Tab. \ref{tab:3D_FBR_k_eff}, where the numerical results by MC \cite{takeda19913} and IUGKS (G2-S8 quadrature) \cite{tan2020time} are all listed, and MC is taken as the reference solution. Tab. \ref{tab:3D_FBR_k_eff} shows that compared to UGKS, UGKWP has a small error for both $k_\text{eff}$ and $\Delta \rho$ in rod-out, and a similar trend is observed for rod-in, demonstrating a slightly better agreement with MC.

\begin{figure}[!hptb]
 \centering
 \subfloat[rod-in, $\psi_1$]{
 	\includegraphics[width = 0.45\textwidth, trim=10 0 0 2, clip]{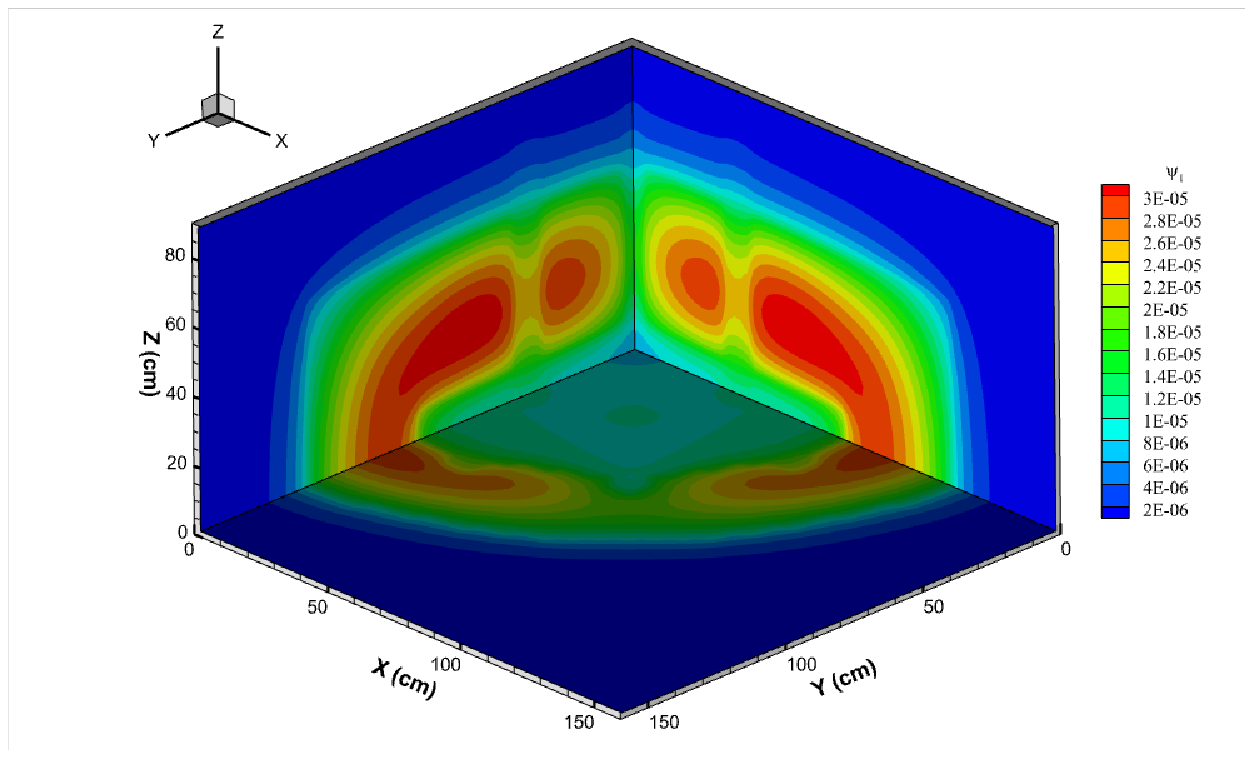}
 }\hfill
 \subfloat[relative error of $\psi_1$ with IUGKS]{               
 	\includegraphics[width = 0.45\textwidth, trim=10 0 0 2, clip]{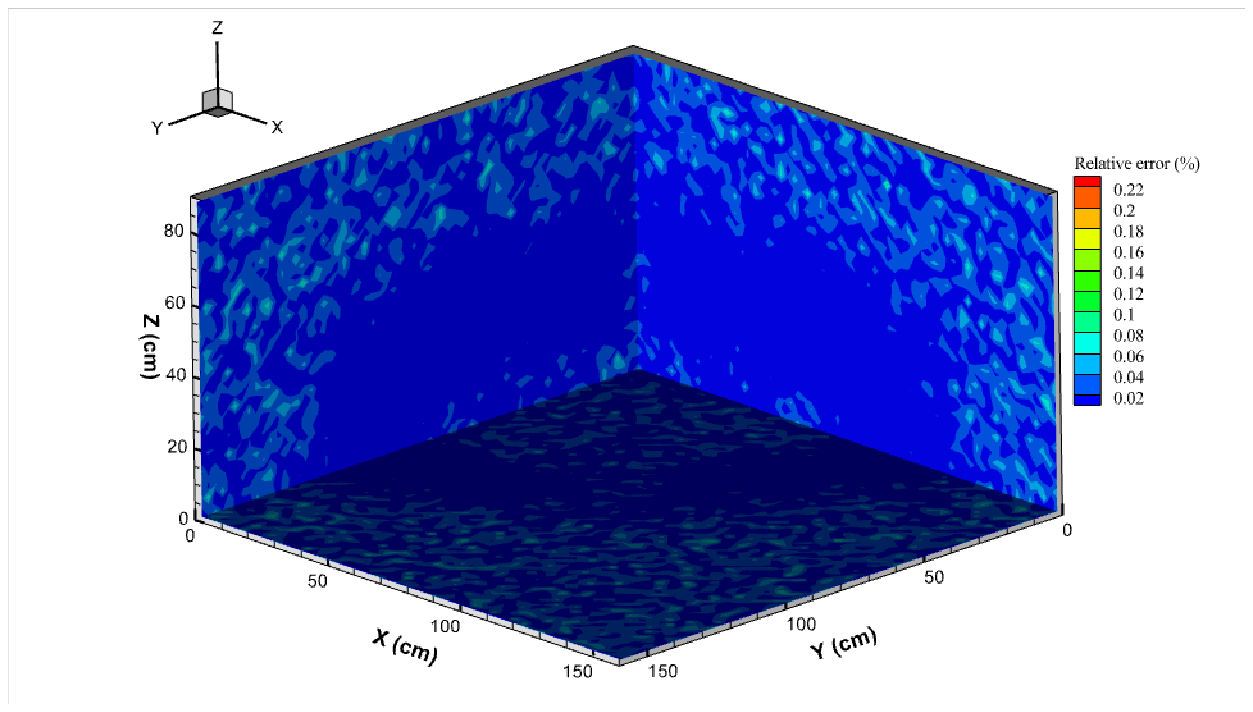}
 }\\
 \subfloat[rod-in, $\psi_2$]{
 	\includegraphics[width = 0.45\textwidth, trim=10 0 0 2, clip]{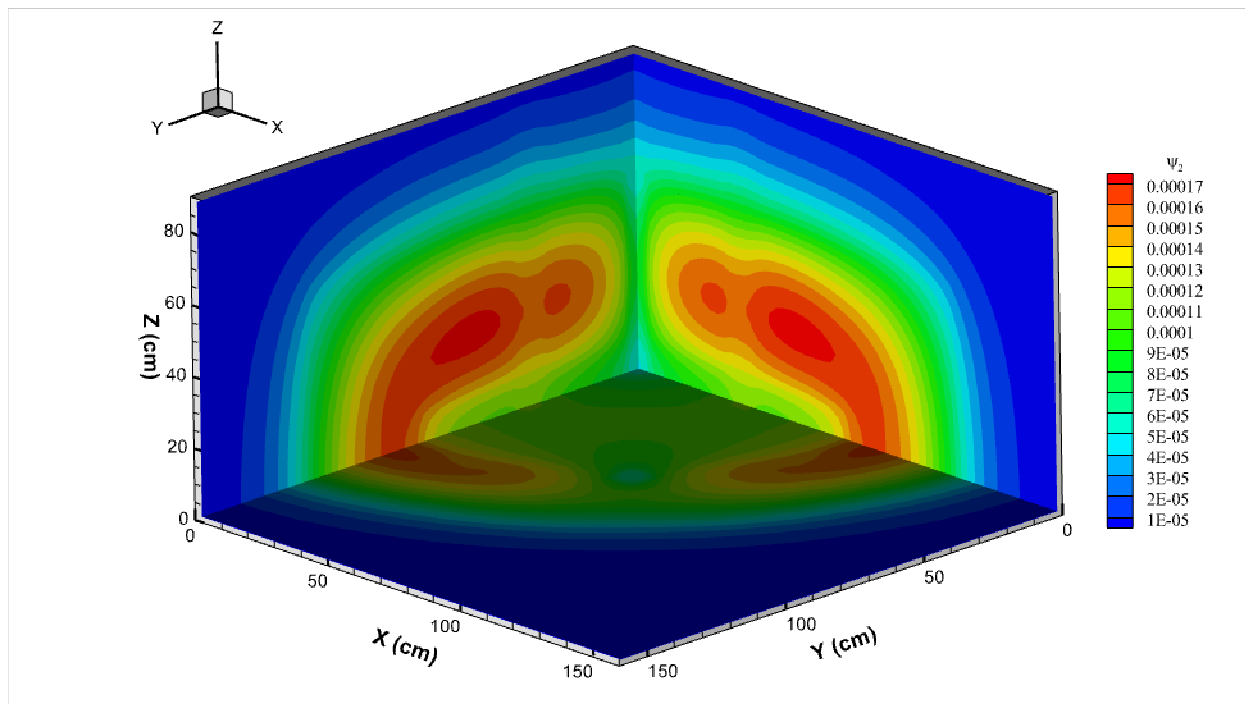}
 }\hfill
 \subfloat[relative error of $\psi_2$ with IUGKS]{               
 	\includegraphics[width = 0.45\textwidth, trim=10 0 0 2, clip]{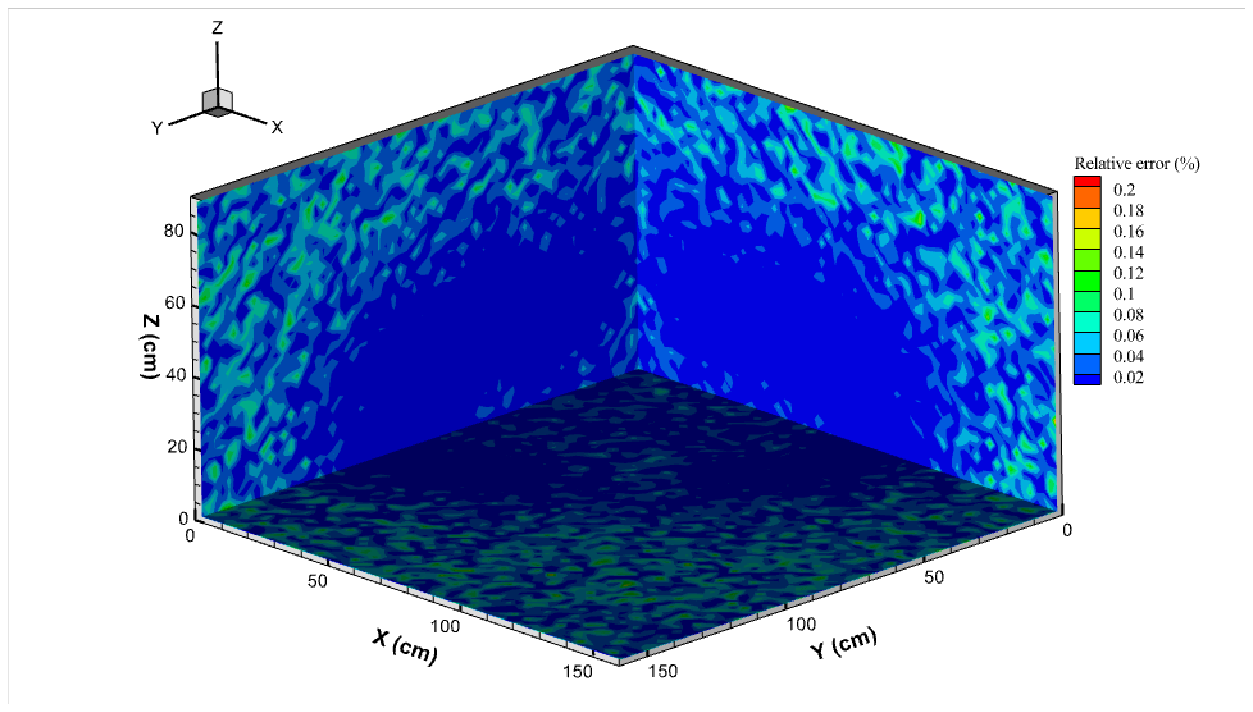}
 }\\
 \subfloat[rod-in, $\psi_3$]{
 	\includegraphics[width = 0.45\textwidth, trim=10 0 0 2, clip]{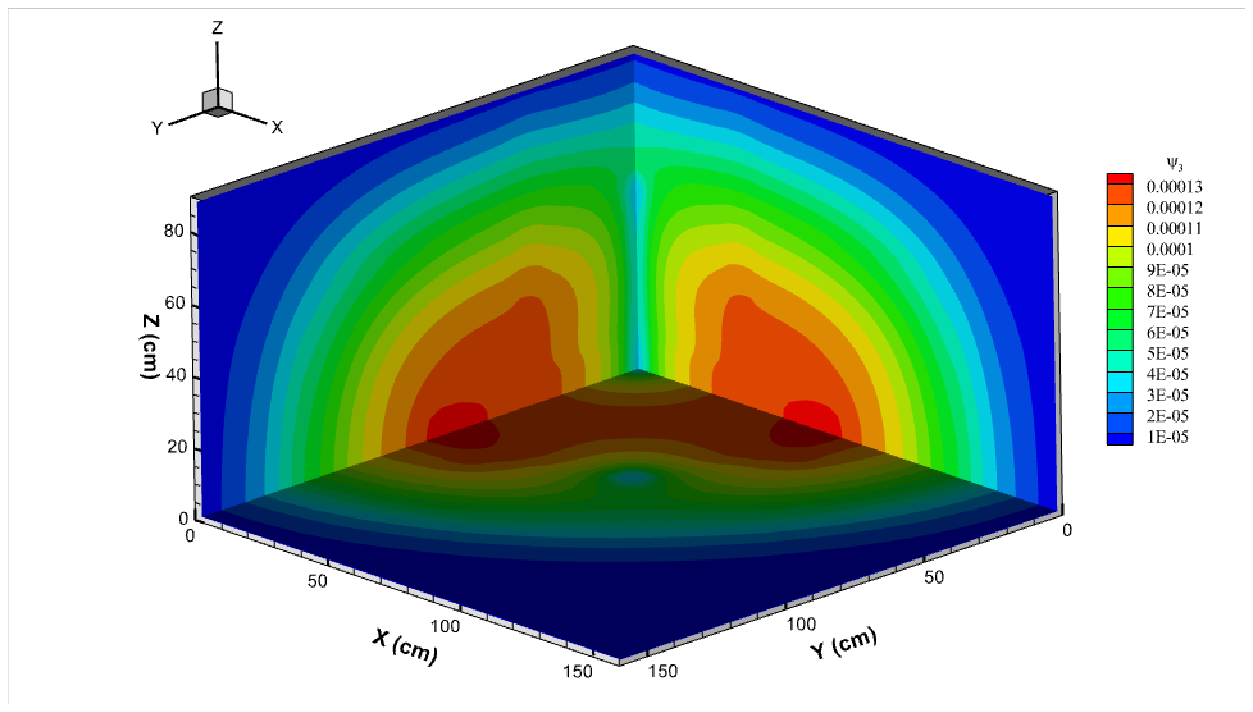}
 }\hfill
 \subfloat[relative error of $\psi_3$ with IUGKS]{               
 	\includegraphics[width = 0.45\textwidth, trim=10 0 0 2, clip]{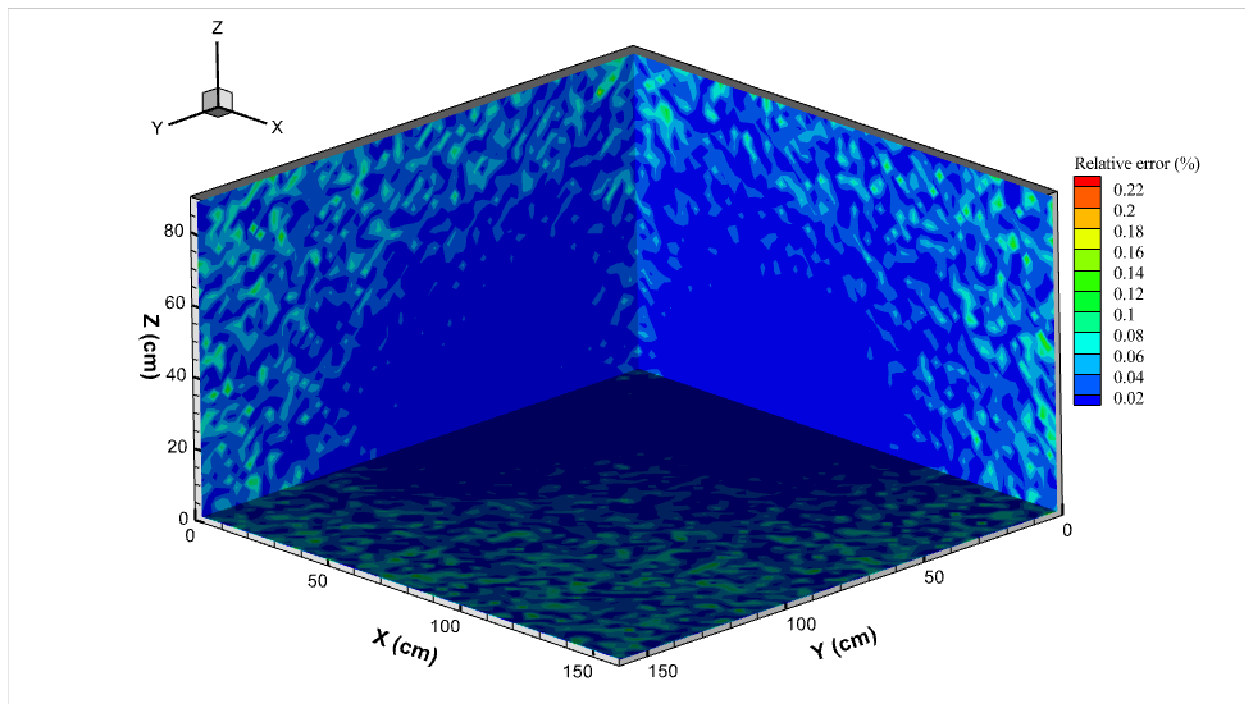}
 }\\
 \subfloat[rod-in, $\psi_4$]{
 	\includegraphics[width = 0.45\textwidth, trim=10 0 0 2, clip]{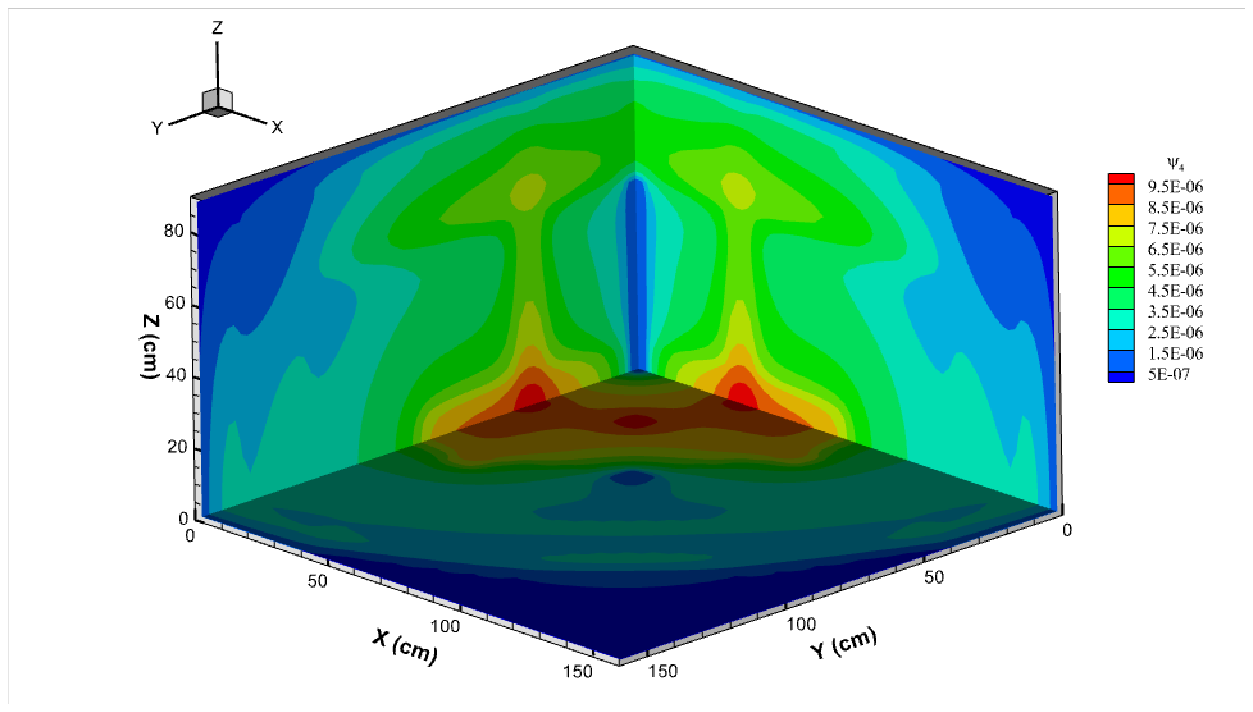}
    \label{fig:FBR_rod_in_4}
 }\hfill
 \subfloat[relative error of $\psi_4$ with IUGKS]{               
 	\includegraphics[width = 0.45\textwidth, trim=10 0 0 2, clip]{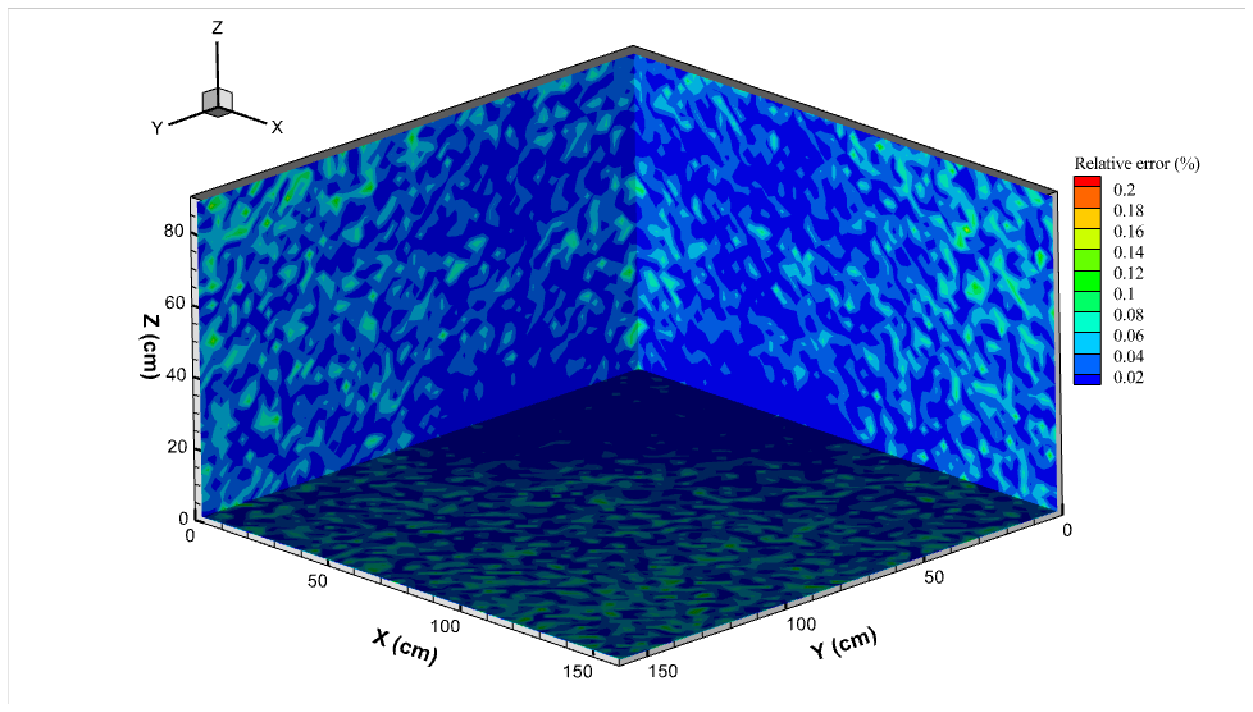}
 }
 \caption{(Four-group axially heterogeneous FBR problem in Sec. \ref{sec:FBR}) 
The numerical solution of $\psi_i, i = 1, \cdots, 4$ for rod-in, and the relative error with IUGKS. Here, the left column is the numerical solution of $\psi_i, i = 1, \cdots, 4$, and the right column is the relative error with IUGKS. }  
 \label{fig:FBR_rod_in}
 \end{figure}

\begin{figure}[!hptb]
	\centering
	\subfloat[$\psi_{1}$ ($y=z=0$)]{
		\includegraphics[width = 0.23\textwidth]{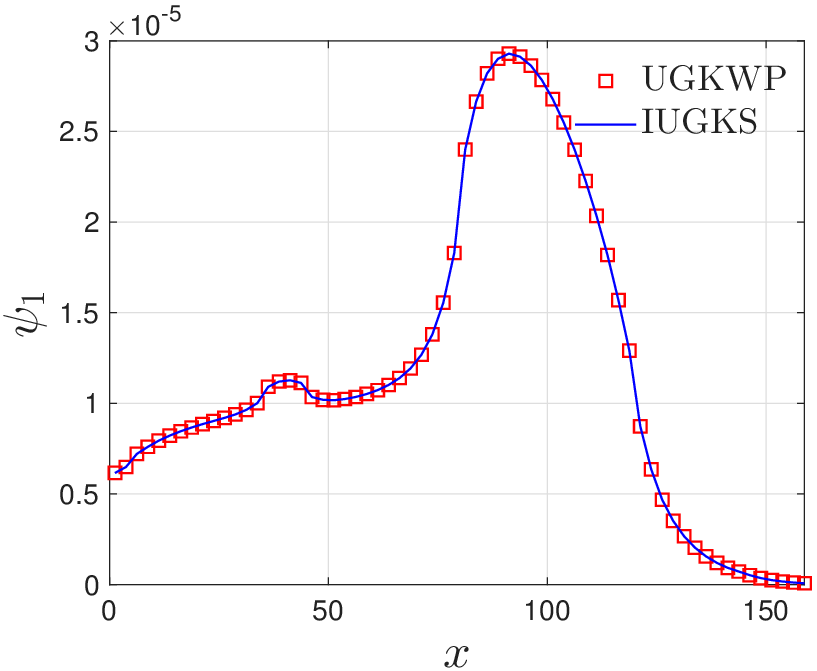}
	} \hfill
    \subfloat[$\psi_{2}$ ($y=z=0$)]{
		\includegraphics[width = 0.23\textwidth]{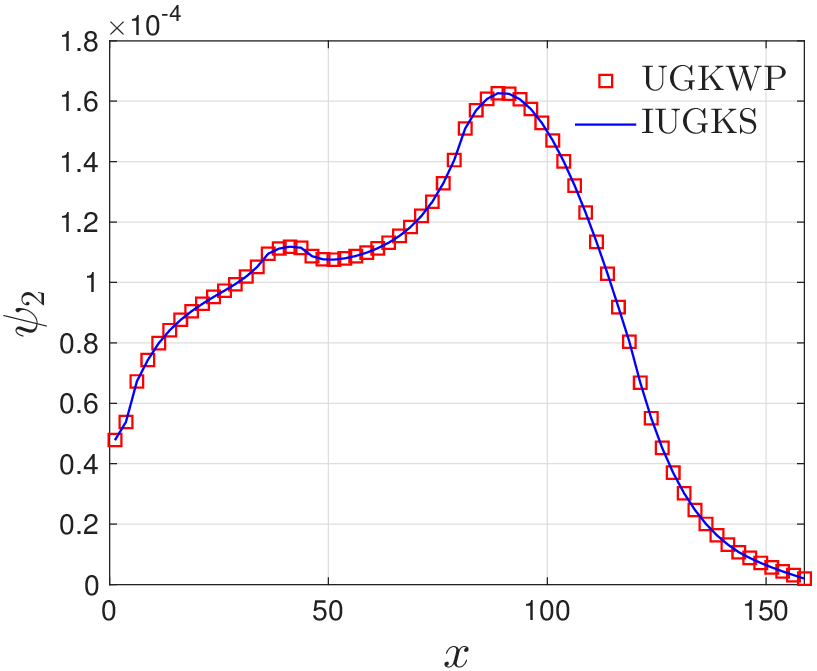}
	} \hfill
    \subfloat[$\psi_{3}$ ($y=z=0$)]{
		\includegraphics[width = 0.23\textwidth]{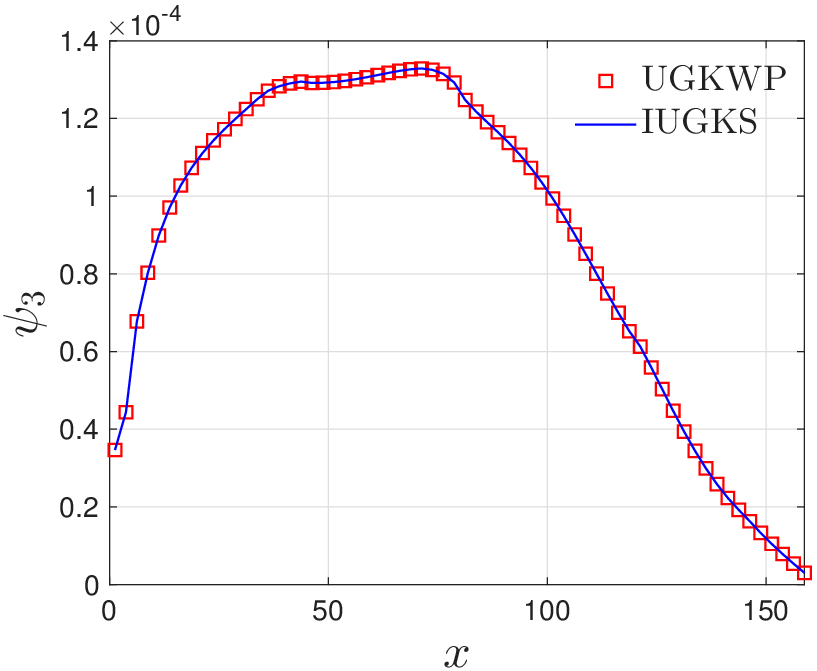}
	} \hfill
    \subfloat[$\psi_{4}$ ($y =z =0$)]{
		\includegraphics[width = 0.23\textwidth]{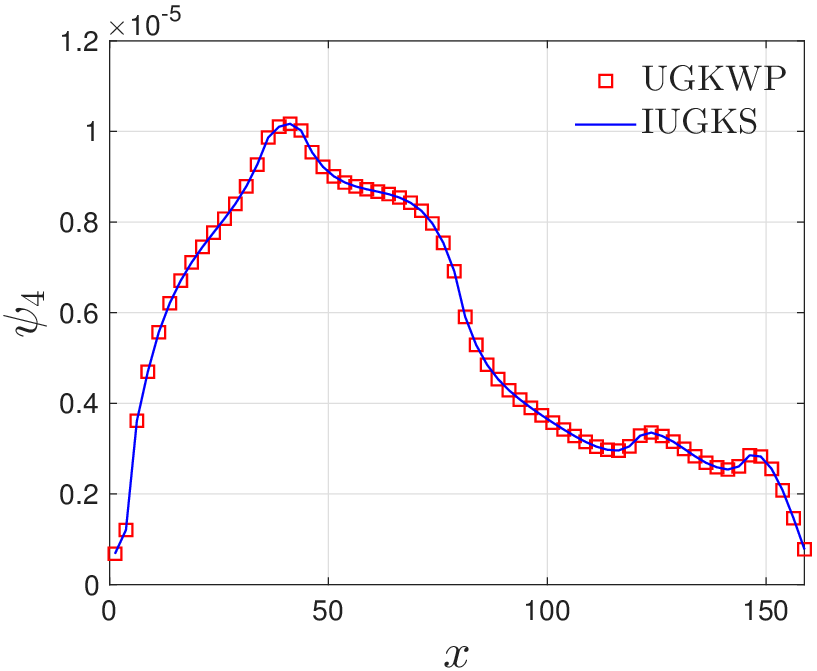}
	} \\
	\subfloat[$\psi_{1}$ ($x=y, z=37.5$)]{               
		\includegraphics[width = 0.23\textwidth]{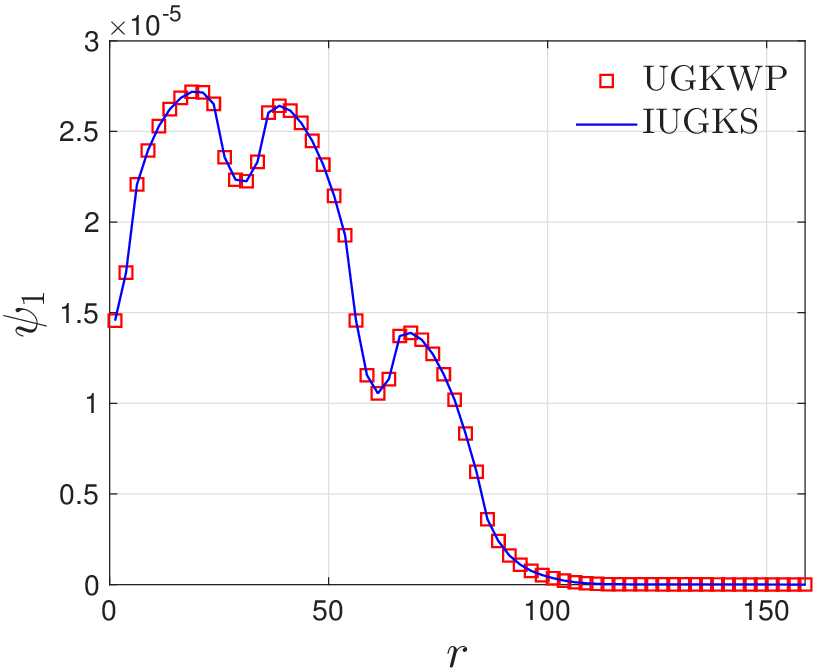}
	}\hfill
	\subfloat[$\psi_{2}$ ($x=y, z=37.5$)]{               
		\includegraphics[width = 0.23\textwidth]{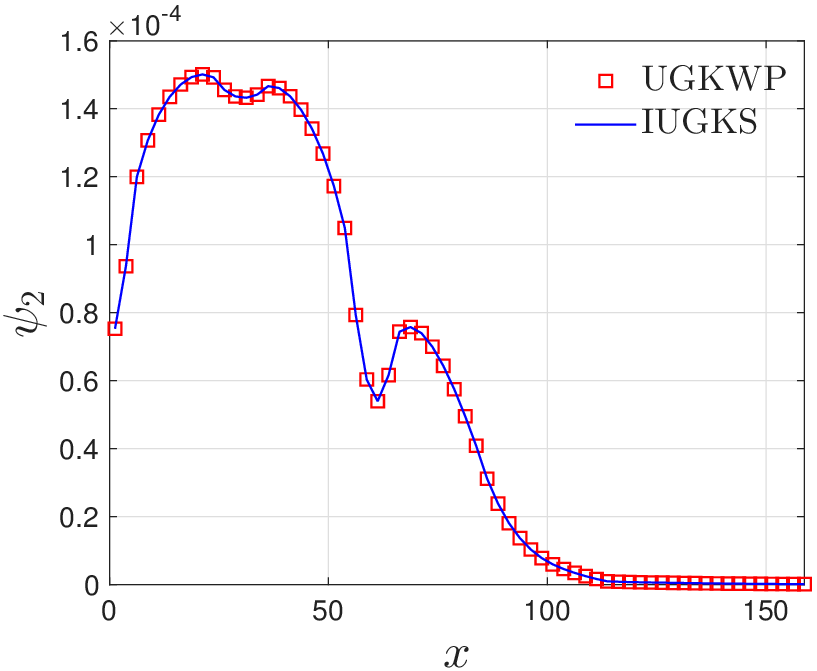}
	} \hfill
	\subfloat[$\psi_{3}$ ($x=y, z=37.5$)]{               
		\includegraphics[width = 0.23\textwidth]{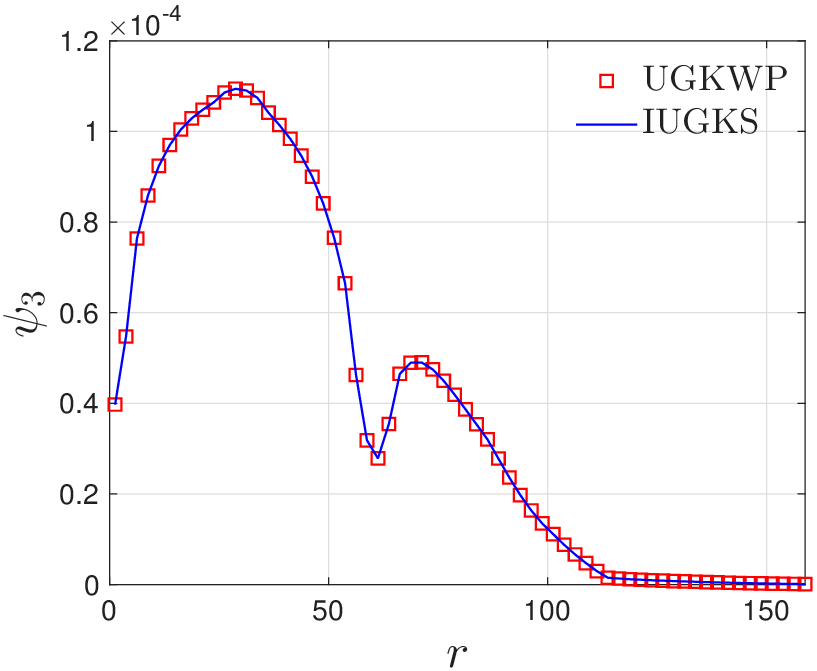}
	}\hfill
	\subfloat[$\psi_{4}$ ($x=y, z=37.5$)]{               
		\includegraphics[width = 0.23\textwidth]{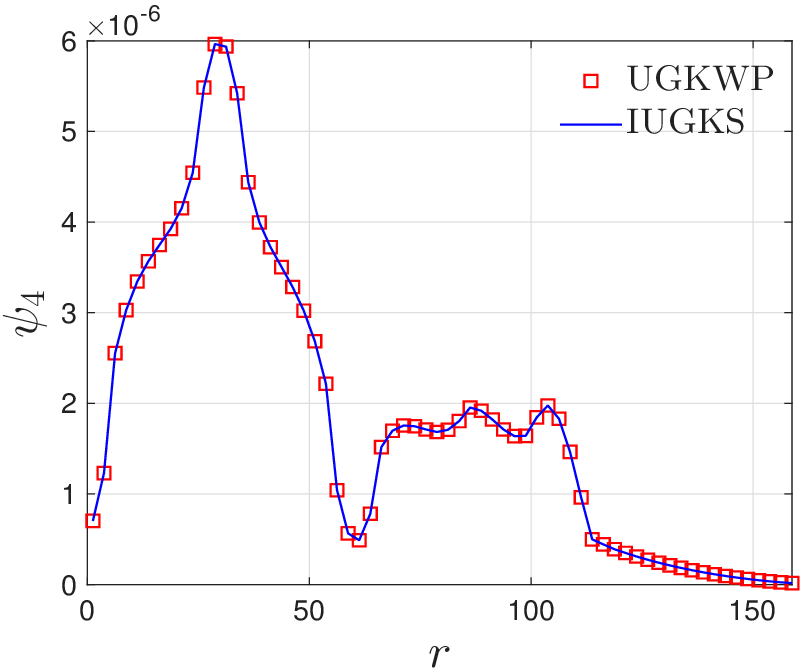}
	}
	\caption{(Four-group axially heterogeneous FBR problem in Sec. \ref{sec:FBR}) The numerical solution of $\psi_i, i = 1, \cdots, 4$ for rod-in along different lines. The first row is $y = z = 0$, and the second row is $x = y, z = 37.5$.  }  
	\label{fig:fbr_pi_34_along_line}
\end{figure}

The numerical solution of the macroscopic scalar flux $\psi_i, i = 1, 2, 3, 4$ and the relative error with UGKS for rod-in is plotted in Fig. \ref{fig:kuca_pi_1} and \ref{fig:kuca_pi_2}, respectively. The behavior of the four macroscopic scalar flux is quite different, while the relative error with UGKS is all quite small, less than $0.2\%$. It also shows that the effect of the control rod position is stronger when the group energy is smaller, as shown in Fig. \ref{fig:FBR_rod_in_4}. Moreover, the macroscopic scalar along $y = z = 0$, and $x = y, z = 37.5$ is shown in Fig. \ref{fig:fbr_pi_34_along_line}. From it, we can clearly see the sharp variation of the macroscopic scalar flux, and UGKWP can well capture it. It also shows that the numerical solution of UGKWP and UGKS is nearly the same, indicating the consistency of these two methods. 

 \begin{table}[!htbp]
 	\centering
 	\def\arraystretch{1.5}
 	{\footnotesize
 		\begin{tabular}{c||c|c c c c c c}
 			code &groups & core & IB & RB & AB & CRP & CR \\
 			%\hline
 			%\multicolumn{7}{c}{Exact Monte-Carlo} \\
 			\hline
 			&1 & 2.0334 & 1.2562 & 0.1589 & 0.2941 & 1.2230 & 1.2150 \\
 			Monte-Carlo&2 & 11.501 & 11.237 & 1.6736 & 3.0115 & 8.6309 & 6.7771 \\
 			&3 & 8.2784 & 11.931 & 2.1521 & 3.5408 & 7.5128 & 3.9699 \\
 			&4 & 0.3409 & 0.7432 & 0.1797 & 0.3752 & 0.5307 & 0.0867 \\
 			%\hline
 			%\multicolumn{7}{c}{IUGKS} \\
 			\hline
 			&1 & 2.0362 & 1.2657 & 0.1611 & 0.2966 & 1.2273 & 1.2338 \\
 			IUGKS&2 & 11.498 & 11.261 & 1.6762 & 3.0242 & 8.6445 & 6.8419 \\
 			&3 & 8.2689 & 11.946 & 2.1535 & 3.5573 & 7.5092 & 4.0205 \\
 			&4 & 0.3402 & 0.7421 & 0.1802 & 0.3769 & 0.5312 & 0.0875 \\
 			%\hline
 			%\multicolumn{7}{c}{UGKWP} \\
 			\hline
 			&1 & 2.0346 & 1.2599 & 0.1598 & 0.2950 & 1.2254 & 1.2247 \\
 			UGKWP&2 & 11.500 & 11.252 & 1.6757 & 3.0212 & 8.6391 & 6.8193 \\
 			&3 & 8.2703 & 11.942 & 2.1541 & 3.5525 & 7.5099 & 4.0044 \\
 			&4 & 0.3403 & 0.7424 & 0.1800 & 0.3762 & 0.5312 & 0.0872 \\
 		\end{tabular}
 		\caption{(Four-group axially heterogeneous FBR problem in Sec. \ref{sec:FBR}) The region-averaged macroscopic scalar flux for different materials (with the units of 1E-05).
                }
 		\label{tab:3D_FBR_region_flux}
 	}
 \end{table}

To show the behavior of different materials, the region-averaged macroscopic scalar flux is listed in Tab. \ref{tab:3D_FBR_region_flux}. For different materials, the behavior of $\psi$ is changed significantly. Generally speaking, for the core region, which is the main source of the fission neutrons, $\psi$ is the largest. Moreover, compared with different energy groups, Group 4, which has the smallest energy, exhibits the smallest macroscopic scalar flux for all materials. The relative error for the region-averaged macroscopic scalar flux of IUGKS and UGKWP with MC is listed in Tab. \ref{tab:relative_error_region_flux}, showing the accuracy of UGKWP from another perspective.

  \begin{table}[!htbp]
 	\centering
 	\def\arraystretch{1.5}
 	{\footnotesize
 	 	\begin{tabular}{c||c|ccccccc}
 		code&groups & core & IB & RB & AB & CRP & CR \\
 		\hline
            &1 & 5.90E-04 & 2.95E-03 & 5.60E-03 & 3.30E-03 & 1.96E-03 & 7.98E-03 \\
 		UGKWP&2 & 8.69E-05 & 1.33E-03 & 1.25E-03 & 3.22E-03 & 9.50E-04 & 6.23E-03 \\
 		&3 & 9.78E-04 & 9.22E-04 & 9.29E-04 & 3.30E-03 & 3.86E-04 & 8.69E-03 \\
 		&4 & 1.76E-03 & 1.02E-03 & 1.22E-03 & 2.64E-03 & 8.86E-04 & 5.92E-03 \\
 		\hline
            &1 & 1.38E-03 & 7.60E-03 & 1.39E-02 & 8.53E-03 & 3.53E-03 & 1.55E-02 \\
 		IUGKS&2 & 3.04E-04 & 2.14E-03 & 1.57E-03 & 4.23E-03 & 1.58E-03 & 9.57E-03 \\
 		&3 & 1.15E-03 & 1.27E-03 & 6.60E-04 & 4.67E-03 & 4.58E-04 & 1.27E-02 \\
 		&4 & 2.04E-03 & 1.44E-03 & 2.41E-03 & 4.55E-03 & 9.21E-04 & 9.37E-03 
 	\end{tabular}
 	}
 	\caption{(Four-group axially heterogeneous FBR problem in Sec. \ref{sec:FBR}) The relative error of the region-averaged macroscopic scalar flux of UGKWP and IUGKS compared with MC for different materials. } 	\label{tab:relative_error_region_flux}
 \end{table}

%\input{num_bak.tex}
%-------------------------------------------------
\section{Conclusion}
\label{sec:conclusion}
In this work, the unified gas-kinetic wave-particle method (UGKWP) is utilized to solve the neutron transport equation (NTE). It begins from the unified gas kinetic method (UGKS), and then the unified gas-kinetic particle method (UGKP), which is a particle implementation of UGKS. Finally, UGKWP is presented based on UGKP, where the collisional particles of the re-sampled particles will not be re-sampled any longer, and the numerical flux contributed by these particles is obtained by subtracting the numerical flux contributed by the free transport particles from the total numerical flux in the macroscopic level. This method is also successfully extended to the multi-group model of NTE. The spatially 1D, 2D, and 3D numerical examples are studied to validate UGKWP for NTE.

 \section*{Acknowledgments}
 We thank Prof. Chang Liu from the Institute of Applied Physics and Computational Mathematics for his valuable suggestions. This work of Yanli Wang is partially supported by the Science Challenge program (NO. TZ2025016), the Foundation of the President of China Academy of Engineering Physics (YZJJZQ2022017), and the National Natural Science Foundation of China (Grant No. 12171026, and 12031013). This work of Shuang Tan is partially supported by the National Natural Science Foundation of China (Grant No. 12572329).
\FloatBarrier 
\bibliographystyle{plain}
\bibliography{main.bib}
\appendix
\newpage
\appendix
\section{Appendix}
\label{sec:app}
\subsection{Boundary conditions}
\label{APP:boundary}
In neutron transport simulations, the classical boundary conditions include the inflow boundary, the reflect boundary and the vacuum boundary \cite{lux1991monte}. Appropriate boundary treatment is critical for preserving physical consistency and numerical accuracy. We still consider the spatially 1D single-group model. The spatial domain is \( x \in [x_L, x_R] \), with the angular variable \( \mu \in [-1,1] \) representing the cosine of the particle's travel direction. 
\paragraph{Inflow boundary bondition}
For prescribed inflow boundaries, neutrons entering the domain are sampled from a known angular flux distribution at boundary $x = x_b$ with the form
\begin{equation}
	\label{eq:inflow_boundary}
	\phi(t, x_b, \mu) = 
		\left\{\begin{aligned}
			\phi(t, x_L, \mu), \qquad \mu > 0 , \\
			\phi(t, x_R, \mu), \qquad \mu < 0 .
		\end{aligned}
		\right. 
\end{equation}
At the left boundary \( x_b = x_L \), incoming particles satisfy \( \mu > 0 \), and the total incoming flux per unit area and unit time is given by:
\begin{equation}
	\label{eq:inflow_flux}
	J_\text{in}(t) = \int_{0}^{1} \mu \, \phi(t, x_L, \mu) \, d\mu.
\end{equation}
%\paragraph{Probability Density Function (PDF).}  
To inject particles such that their angular distribution matches the incoming flux \cite{lux1991monte, nance1997role}, the direction \( \mu \) must be sampled according to a probability density function proportional to the incoming particle current at each angle. This leads to:
\begin{equation}
	\label{eq:inflow_pdf}
	P(t, \mu) = \frac{\mu \, \phi(t, x_L, \mu)}{J_\text{in}(t)}, \quad \mu \in [0,1],
\end{equation}
which ensures proper reproduction of the total angular flux during sampling.
The cumulative distribution function corresponding to the above PDF is:
\begin{equation}
	\label{eq:inflow_cdf}
	F(t, \mu) = \int_{0}^{\mu} P(s) \, ds = \frac{1}{J_\text{in}(t)} \int_{0}^{\mu} s \, \phi_\text(t, x_L, s) \, ds.
\end{equation}
To sample \( \mu \), one generates a uniform random number \( r \in [0, 1] \) and computes the inverse transform:
\begin{equation*}
	\mu = F^{-1}(r).
\end{equation*}
If the cumulative distribution function $F(t,\mu)$ is analytically intractable and inverse transform sampling is not feasible, the particle direction $\mu$ can be sampled using the $P(\mu)$ in combination with the acceptance-rejection method. Similarly, the right boundary is treated in the same manner as the left boundary.

\paragraph{Reflective boundary condition}
Reflective boundaries enforce specular reflection, preserving particle energy and group while reversing the normal component of velocity.
At the left boundary \( x_b = x_L \), a particle with $(x_p, v\mu_p, m_e)$ with \( \mu < 0 \) are reflected back into the domain by updating their direction as:
\begin{equation}
	\label{eq:reflect_bd}
	(x_p, v\mu_p, m_e) \rightarrow (x_p, -v\mu_p, m_e).
\end{equation}
The particle's position remains unchanged, and transport continues.
Similarly, at the right boundary \( x_b = x_R \), particles with \( \mu > 0 \) are reflected using the same rule. This condition ensures particle number conservation at the boundary.

\paragraph{Vacuum boundary condition}
Vacuum boundaries assume zero incoming angular flux; particles reaching the boundary exit the system and are not re-injected.
At the left boundary \( x_b = x_L \), particles with \( \mu < 0 \) are absorbed upon contact with the boundary:
\[
\text{If } x_p \leq x_L \text{ and } \mu < 0, \quad \text{delete particle}.
\]
At the right boundary \( x_b = x_R \), the termination condition applies to particles with \( \mu > 0 \). 

\subsection{Parameters of multi-group problem in Sec. \ref{sec:3D_multi_group}}
\label{APP:3D}
\subsubsection{Two-group cross sections and energy group structure for KUCA core 2-group benchmark in Sec. \ref{sec:kuca}}
\label{APP:kuca}
The group-wise cross-section data are provided as follows. $\Sigma_{a,g}, g=1,2$ is the absorption cross section which satisfies 
\[
\begin{aligned}
    \Sigma_{g} = \Sigma_{s,1\rightarrow 1}+\Sigma_{s,1\rightarrow 2}+\Sigma_{a,g}, g = 1, 2. 
\end{aligned}
\]

\begin{table}[!htbp]
	\centering
	\def\arraystretch{1.3}
	{\footnotesize
	\begin{tabular}{ccccccc}
		\toprule
		\textbf{region} & \textbf{group (g)} & $\Sigma_{g}$ & $\nu\Sigma_{f,g}$ & $\Sigma_{s,1\rightarrow g}$ & $\Sigma_{s,2\rightarrow g}$ & $\Sigma_{a,g}$ \\
		\midrule
		Core & 1 & 2.23775E-01 & 9.09319E-03 & 1.92123E-01 & 0.0           & 8.52709E-03 \\
		     & 2 & 1.03864E+00 & 2.90183E-01 & 2.28253E-02 & 8.80439E-01 & 1.58196E-01 \\
		\midrule
		Control Rod & 1 & 8.52325E-02 & 0.0         & 6.77241E-02 & 0.0           & 1.74439E-02 \\
		            & 2 & 2.17460E-01 & 0.0         & 6.45461E-05 & 3.52358E-02 & 1.82224E-01 \\
		\midrule
		Reflector & 1 & 2.50367E-01 & 0.0         & 1.93446E-01 & 0.0           & 4.16392E-04 \\
		          & 2 & 1.64482E+00 & 0.0         & 5.65042E-02 & 1.62452E+00 & 2.02999E-02 \\
		\midrule
		Empty (void) & 1 & 1.28407E-02 & 0.0         & 1.27700E-02 & 0.0           & 4.65132E-05 \\
		             & 2 & 1.20676E-02 & 0.0         & 2.40997E-05 & 1.07387E-02 & 1.32890E-03 \\
		\bottomrule
	\end{tabular}
	}
	\caption{(The KUCA core two-group benchmark problem in Sec. \ref{sec:kuca}) Cross section parameters for each region in different energy groups (Units: $\mathrm{cm}^{-1}$).}
    \label{tab:KUCA_cross}
\end{table}

% \begin{table}[!htbp]
% 	\centering
% 	\def\arraystretch{1.5}
% 	{\footnotesize
% 		\begin{tabular}{c|c|ccccc}
% 		\textbf{region} & \textbf{group (g)} & $\Sigma_{g}$ & $\nu\Sigma_{f,g}$ & $\Sigma_{s,1\rightarrow g}$ & $\Sigma_{s,2\rightarrow g}$ & $\Sigma_{a,g}$ \\
% 		\hline
% 		Core& 1 & 2.23775E-01 & 9.09319E-03 & 1.92123E-01 & 0.0           & 8.52709E-03 \\
% 		& 2 & 1.03864E+00 & 2.90183E-01 & 2.28253E-02 & 8.80439E-01 & 1.58196E-01 \\
% 		\hline
% 		Control Rod& 1 & 8.52325E-02 & 0.0         & 6.77241E-02 & 0.0           & 1.74439E-02 \\
% 		& 2 & 2.17460E-01 & 0.0         & 6.45461E-05 & 3.52358E-02 & 1.82224E-01 \\
% 		\hline
% 		Reflector& 1 & 2.50367E-01 & 0.0         & 1.93446E-01 & 0.0           & 4.16392E-04 \\
% 		& 2 & 1.64482E+00 & 0.0         & 5.65042E-02 & 1.62452E+00 & 2.02999E-02 \\
% 		\hline
% 		Empty (void)& 1 & 1.28407E-02 & 0.0         & 1.27700E-02 & 0.0           & 4.65132E-05 \\
% 		& 2 & 1.20676E-02 & 0.0         & 2.40997E-05 & 1.07387E-02 & 1.32890E-03 
% 	\end{tabular}
%     }
% 	 	\caption{(The KUCA core two-group benchmark problem in Sec. \ref{sec:kuca}) Cross section parameters for each region in different energy group. (Units: $\mathrm{cm}^{-1}$)}
%  \end{table}

\begin{table}[!htbp]
	\centering
	\def\arraystretch{1.3}
	{\footnotesize
 	\begin{tabular}{cccc}
            \toprule
 		\textbf{group (g)} & $E_{\text{sup}}$  & $E_{\text{inf}}$ & $\chi_{g}$ \\
 		\midrule
 		1 & 1.0000E+07 & 6.8256E-01& 1.00000 \\
 		2 & 6.8256E-01 & 1.0000E-05 & 0.0 \\
        \bottomrule
 	\end{tabular}
 	}
 	\caption{(The KUCA core two-group benchmark problem in Sec. \ref{sec:kuca}) Energy group structure.}
        \label{tab:KUCA_energy}
 \end{table}

% \begin{table}[!htbp]
% 	\centering
% 	\def\arraystretch{1.5}
% 	{\footnotesize
%  	\begin{tabular}{cccc}
%  		\textbf{group (g)} & $E_{sup}$ \qquad\qquad $E_{inf}$ &  & $\chi_{g}$ \\
%  		\hline
%  		1 & $1.0000\text{E}+7 \quad 6.8256\text{E}-1$ & & 1.00000 \\
%  		2 & $6.8256\text{E}-1 \quad 1.0000\text{E}-5$ & & 0.00000 \\
%  	\end{tabular}
%  	}
%  	\caption{(The KUCA core two-group benchmark problem in Sec. \ref{sec:kuca}) Energy group structure.}
%  \end{table}
\subsubsection{Cross sections and energy group structure for the four-group axially heterogeneous FBR core problem in Sec. \ref{sec:FBR}}
\label{APP:fbr}
 \begin{table}[htbp!]
 	\centering
 	\renewcommand{\arraystretch}{1.3}
 	{\footnotesize
 	\begin{tabular}{lcccc}
            \toprule
 		\textbf{region} & \textbf{group (g)} & \(\Sigma_{g}\) & \(\Sigma_{a,g}\) & \(\nu\Sigma_{f,g}\) \\
 		\midrule
 		& 1 & 1.4568E-01 & 4.5531E-03 & 2.06063E-02 \\
 		Core& 2 & 2.0517E-01 & 5.2540E-03 & 6.19571E-03 \\
 		& 3 & 3.2931E-01 & 8.0016E-03 & 9.40342E-03 \\
 		& 4 & 3.8910E-01 & 2.7449E-02 & 2.60689E-02 \\
		\midrule
 		& 1 & 1.1968E-01 & 4.3282E-03 & 1.89496E-02 \\
 		Radial \& Inner Blanket& 2 & 2.4219E-01 & 1.9996E-03 & 1.75256E-04 \\
 		& 3 & 5.5647E-01 & 6.7903E-02 & 2.06978E-04 \\
 		& 4 & 3.7943E-01 & 1.5801E-02 & 1.13451E-03 \\
		\midrule
 		& 1 & 1.6439E-01 & 3.5514E-03 & 1.71702E-02 \\
 		Axial Blanket& 2 & 2.2051E-01 & 1.4864E-03 & 1.26064E-02 \\
 		& 3 & 4.5444E-01 & 5.3530E-02 & 1.52784E-04 \\
 		& 4 & 3.8456E-01 & 1.3469E-02 & 7.67302E-04 \\
 		\midrule
 		& 1 & 1.6561E-01 & 6.3915E-04 & 0.0 \\
 		Axial Reflector& 2 & 1.6666E-01 & 1.0667E-04 & 0.0 \\
 		& 3 & 2.6633E-01 & 1.2047E-03 & 0.0 \\
 		& 4 & 8.3491E-02 & 1.3638E-03 & 0.0 \\
		\midrule
 		& 1 & 1.7174E-01 & 1.1330E-03 & 0.0 \\
 		Radial Reflector& 2 & 2.1786E-01 & 9.0792E-04 & 0.0 \\
 		& 3 & 4.6712E-01 & 1.9450E-03 & 0.0 \\
 		& 4 & 7.9519E-01 & 5.7026E-03 & 0.0 \\
		\midrule
 		& 1 & 1.8339E-01 & 5.9764E-03 & 0.0 \\
 		Control Rod& 2 & 3.6612E-01 & 1.7641E-02 & 0.0 \\
 		& 3 & 6.1557E-01 & 8.6274E-02 & 0.0 \\
 		& 4 & 1.0946E+00 & 4.7659E-01 & 0.0 \\
 		\midrule
 		& 1 & 6.50979E-02 & 3.10744E-04 & 0.0 \\
 		No Filled CRP& 2 & 1.09810E-01 & 1.13062E-04 & 0.0 \\
 		& 3 & 1.86765E-01 & 4.40968E-04 & 0.0 \\
 		& 4 & 2.09933E-01 & 1.07518E-03 & 0.0 \\
		\midrule
 		& 1 & 1.36985E-02 & 7.49000E-05 & 0.0 \\
 		Empty Matrix& 2 & 1.69037E-02 & 3.82435E-05 & 0.0 \\
 		& 3 & 3.12271E-02 & 1.39335E-04 & 0.0 \\
 		& 4 & 6.29337E-02 & 4.95515E-04 & 0.0 \\
        \bottomrule
 	\end{tabular}
 	}
 	\caption{(Four-group axially heterogeneous FBR problem in Sec.~\ref{sec:FBR}) Cross section parameters ($\Sigma_{g},\: \Sigma_{a, g},\: \Sigma_{f,g}$) for each region in different energy group (Units: \(\mathrm{cm}^{-1}\)).}
        \label{tab:FBR_cross_1}
 \end{table}

 \begin{table}[!htbp]
	\centering
	\def\arraystretch{1.3}
	{\footnotesize
	 	\begin{tabular}{cccccc}
        \toprule
		\textbf{region} & \textbf{group (g)} & $\Sigma_{s,1\rightarrow g}$ & $\Sigma_{s,2\rightarrow g}$ & $\Sigma_{s,3\rightarrow g}$ & $\Sigma_{s, 4\rightarrow g}$ \\
		\midrule
		& 1 & 7.04236E-02 & 0.0 & 0.0 & 0.0 \\
		Core& 2 & 3.4976E-02 & 1.9544E-01 & 0.0 & 0.0 \\
		& 3 & 1.8822E-02 & 6.2086E-03 & 3.2086E-01 & 0.0 \\
		& 4 & 0.0 & 7.0208E-02 & 9.9275E-04 & 3.6236E-01 \\
		\midrule
		& 1 & 6.9115E-02 & 0.0 & 0.0 & 0.0 \\
		Radial \& Inner Blanket& 2 & 4.0913E-02 & 3.2063E-01 & 0.0 & 0.0 \\
		& 3 & 2.6626E-03 & 9.5702E-03 & 3.4841E-01 & 0.0 \\
		& 4 & 0.0 & 1.9957E-01 & 1.2719E-03 & 3.6631E-01 \\
		\midrule
		& 1 & 7.1604E-02 & 0.0 & 0.0 & 0.0 \\
		Axial Blanket& 2 & 3.2171E-02 & 1.9046E-01 & 0.0 & 0.0 \\
		& 3 & 2.2170E-03 & 8.5985E-03 & 3.7506E-01 & 0.0 \\
		& 4 & 0.0 & 6.8299E-02 & 1.6835E-03 & 3.7486E-01 \\
		\midrule
		& 1 & 1.1565E-02 & 0.0 & 0.0 & 0.0 \\
		Axial Reflector& 2 & 4.8731E-02 & 1.6081E-01 & 0.0 & 0.0 \\
		& 3 & 8.4649E-04 & 6.6409E-03 & 2.6501E-01 & 0.0 \\
		& 4 & 0.0 & 6.5735E-02 & 2.4175E-03 & 3.0547E-01 \\
		\midrule
		& 1 & 1.2335E-02 & 0.0 & 0.0 & 0.0 \\
		Radial Reflector& 2 & 4.6103E-02 & 2.1106E-01 & 0.0 & 0.0 \\
		& 3 & 1.1321E-03 & 6.2710E-03 & 4.4305E-01 & 0.0 \\
		& 4 & 0.0 & 1.0381E-02 & 2.7112E-03 & 7.8947E-01 \\
		\midrule
		& 1 & 1.3437E-02 & 0.0 & 0.0 & 0.0 \\
		Control Rod& 2 & 4.3777E-02 & 3.1858E-01 & 0.0 & 0.0 \\
		& 3 & 2.0655E-02 & 8.9843E-02 & 1.5951E-01 & 0.0 \\
		& 4 & 0.0 & 7.1188E-02 & 7.6620E-03 & 6.1826E-01 \\
		\midrule
		& 1 & 4.74407E-02 & 0.0 & 0.0 & 0.0 \\
		No Filled CRP& 2 & 1.76894E-02 & 1.06142E-01 & 0.0 & 0.0 \\
		& 3 & 4.57012E-04 & 3.55466E-03 & 1.85304E-01 & 0.0\\
		& 4 & 0.0 & 1.77599E-07 & 1.01280E-03 & 2.08858E-01 \\
		\midrule
		& 1 & 9.57999E-03 & 0.0 & 0.0 & 0.0 \\
		Empty Matrix& 2 & 3.95552E-03 & 1.64740E-02 & 0.0 & 0.0 \\
		& 3 & 8.80428E-05 & 3.91394E-04 & 3.09104E-02 & 0.0\\
		& 4 & 0.0 & 7.72254E-08 & 1.77293E-04 & 6.24581E-02 \\
		\bottomrule
	\end{tabular}
	}
 	\caption{(Four-group axially heterogeneous FBR problem in Sec.~\ref{sec:FBR}) Scattering cross section parameters for each region in different energy groups (Units: \(\mathrm{cm}^{-1}\)).}
            \label{tab:FBR_cross_2}
 \end{table}
 
 \begin{table}[!htbp]
 	\centering
	\def\arraystretch{1.3}
        {\footnotesize
        \begin{tabular}{cccc}
            \toprule
 		\textbf{group (g)} &
 		 $E_{\text{sup}}$  & $E_{\text{inf}}$ & $\chi_{g}$ \\
 		\midrule
 		1 & 1.0000E+07 & 1.3534E+06 & 0.583319 \\
 		2 & 1.3534E+06 & 8.6517E+04 & 0.405450 \\
 		3 & 8.6517E+04 & 9.6112E+02 & 0.011231 \\
 		4 & 9.6112E+02 & 1.0000E-05 & 0.0\\
       \bottomrule
 	\end{tabular}
        }
 	\caption{(Four-group axially heterogeneous FBR problem in Sec.~\ref{sec:FBR}) Energy group structure.}
            \label{tab:FBR_energy}
 \end{table}

\end{document}